
\input psfig
\input amssym.def
\input amssym
\magnification=1100
\baselineskip = 0.25truein
\lineskiplimit = 0.01truein
\lineskip = 0.01truein
\vsize = 8.5truein
\voffset = 0.2truein
\parskip = 0.10truein
\parindent = 0.3truein
\settabs 12 \columns
\hsize = 5.4truein
\hoffset = 0.4truein
\font\ninerm=cmr9

\setbox\strutbox=\hbox{%
\vrule height .708\baselineskip
depth .292\baselineskip
width 0pt}
\font\caps=cmcsc10

\def\sqr#1#2{{\vcenter{\vbox{\hrule height.#2pt
\hbox{\vrule width.#2pt height#1pt \kern#1pt
\vrule width.#2pt}
\hrule height.#2pt}}}}
\def\square{\mathchoice\sqr46\sqr46\sqr{3.1}6\sqr{2.3}4}
\def\leaderfill{\leaders\hbox to 1em{\hss.\hss}\hfill}
\font\bigtenrm=cmr10 scaled 1400
\tenrm

\vskip 2in
\centerline{\bf {\bigtenrm EXCEPTIONAL SURGERY CURVES}}
\centerline{\bf {\bigtenrm IN TRIANGULATED 3-MANIFOLDS}}
\tenrm
\vskip 14pt
\centerline{MARC LACKENBY}
\vskip 18pt

\vskip 12pt
\centerline {\caps Abstract}

{\ninerm \smallskip
For the purposes of this paper, Dehn surgery along a curve $K$
in a 3-manifold $M$ with slope $\sigma$ is `exceptional' if
the resulting 3-manifold $M_K(\sigma)$ is reducible or a solid
torus, or the core of the surgery solid torus has finite order
in $\pi_1(M_K(\sigma))$. We show that, providing the exterior
of $K$ is irreducible and atoroidal, and the distance between
$\sigma$ and the meridian slope is more than one, and a homology
condition is satisfied, then there is (up to ambient isotopy) only a finite number
of such exceptional surgery curves in a given compact orientable 
3-manifold $M$, with $\partial M$ a (possibly empty) union of tori.
Moreover, there is a simple algorithm to find all these
surgery curves, which
involves inserting tangles into the 3-simplices of any given triangulation of $M$.
As a consequence, we deduce some results about the finiteness
of certain unknotting operations on knots in the 3-sphere.
\smallskip}

\tenrm

\vskip 12pt
\centerline {\caps 1. Introduction}
\vskip 6pt

Consider the following motivating problem from knot theory.
Let $L$ be a non-trivial knot in $S^3$. If $K$ is an
unknotted curve disjoint from $L$, then Dehn surgery along $K$ with
slope $1/q$ has the effect of adding $\vert q \vert$
full twists to $L$, yielding a knot $L'$, say.
(See Figure 1.2.) Suppose that $L'$ is the unknot, or (more
generally) that $L'$ has smaller genus
than that of $L$. Then, for a given knot $L$, are there only a finitely
many possibilities for $q$ and $K$ (up to ambient isotopy keeping $L$ fixed)?
The following theorem deals with this question.

\noindent {\bf Theorem 1.1.} {\sl Let $L$ be a knot in $S^3$
which is not a non-trivial satellite knot.
Let $K$ be an unknotted curve in $S^3$, disjoint from $L$ and
having zero linking number with $L$. Let $q$ be an integer
with $\vert q \vert > 1$. Suppose that
$1/q$ surgery about $K$ yields a knot $L'$ with
${\rm genus}(L') < {\rm genus}(L)$. Then, for a given knot $L$,
there are only finitely many possibilities for $K$ and $q$ up to
ambient isotopy keeping $L$ fixed,
and there is an algorithm to find them all.}

\vskip 24pt
\centerline{\psfig{figure=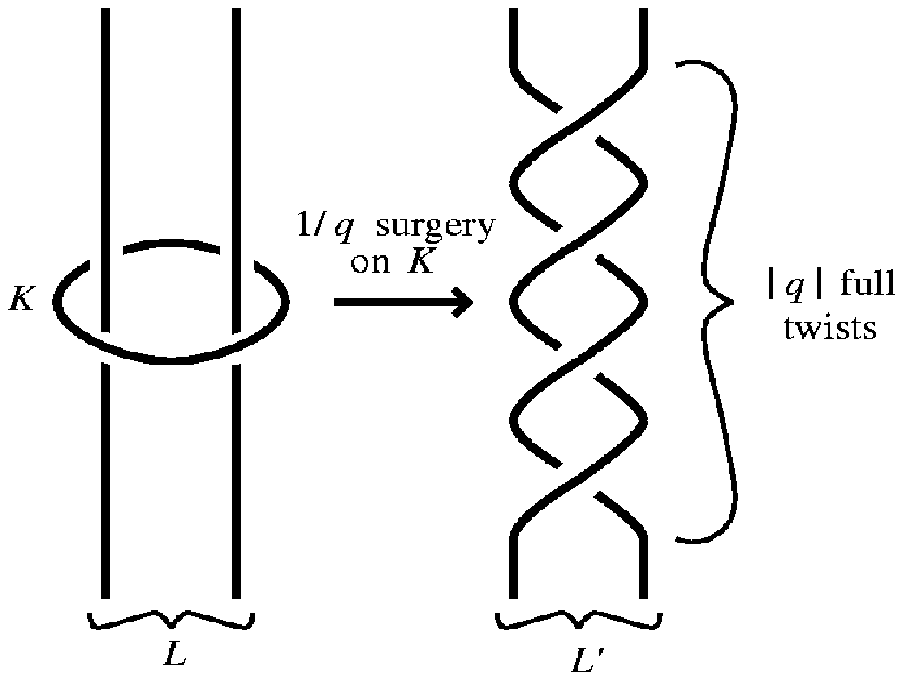}}
\centerline{Figure 1.2.}

Such `unknotting operations' have been the object of
considerable study. For example, the author in [8] dealt with
the case where $K$ bounds a disc intersecting $L$ in
two points of opposite sign, and proved that if such
a surgery reduces the genus of $L$, then there exists
an upper bound on $\vert q \vert$ which depends only on $L$,
not $K$. Theorem 1.1 gives a great deal more than numerical
restrictions on $\vert q \vert$. It provides a
{\sl classification} of all such unknotting operations for a given
knot $L$, when $\vert q \vert > 1$ and the linking number of $K$
and $L$ is zero.

Theorem 1.1 is an almost immediate corollary of new results
on Dehn surgery. Let $M$ be an arbitrary
compact orientable 3-manifold $M$, with $\partial M$
a (possibly empty) union of tori. (In Theorem 1.1,
$M$ is the exterior of the knot $L$.) Our aim is
to find the knots $K$ in $M$ with `exceptional' or `norm-exceptional' surgeries, which
we define as follows.

\noindent {\bf Definition 1.3.} Let $\sigma$ be a slope on
$\partial {\cal N}(K)$ other than the meridional slope $\mu$.
Let $M_K(\sigma)$ be the manifold obtained by Dehn surgery along
$K$ via the slope $\sigma$. Then $\sigma$ is an 
{\sl exceptional slope} and $K$ is an
{\sl exceptional surgery curve} if any of the following holds:
\item{(i)} $M_K(\sigma)$ is reducible,
\item{(ii)} $M_K(\sigma)$ is a solid torus, or
\item{(iii)} the core of the surgery solid torus has
finite order in $\pi_1(M_K(\sigma))$.

\noindent Also, $\sigma$ and $K$ are {\sl norm-exceptional}
if there is some $z \in H_2(M - {\rm int}({\cal N}(K)), \partial M)$
which maps to an element $z_\sigma \in H_2(M_K(\sigma), \partial M_K(\sigma))$,
such that the Thurston norm of $z_\sigma$ less than the Thurston
norm of $z$. (See Section 3 for a definition of the Thurston norm.)

If $K$ satisfies the conditions of Theorem 1.1, then
it is a norm-exceptional surgery curve in $M = S^3 - {\rm int}({\cal N}(L))$.
The reason for distinguishing norm-exceptional surgery curves
from the exceptional case is that, in the former situation, our results 
will be slightly weaker. We restrict attention to knots $K$ with
irreducible atoroidal exteriors. For technical reasons, we 
also have to
assume that $H_2(M  - {\rm int}({\cal N}(K)), \partial M)$
is non-trivial. This implies in particular that the first Betti number of
$M$ must be non-zero. 

We will show that the problem of finding exceptional surgery curves in
a given 3-manifold $M$ falls naturally
into two cases, which depend on $\Delta(\sigma, \mu)$, where
$\Delta(\sigma, \mu)$ is
the intersection number on $\partial {\cal N}(K)$ between the
surgery slope $\sigma$ and the meridian slope $\mu$.
It is not hard to find examples of 3-manifolds $M$
as above containing an infinite
number of pairwise non-isotopic surgery curves $K$ with exceptional 
or norm-exceptional surgery
slopes $\sigma$ satisfying $\Delta(\sigma, \mu) = 1$.
(We will do this in Section 12.) 
However, the main theorem of this paper asserts
that, if $\Delta(\sigma, \mu)> 1$, then
there is only a finite number of possibilities for $K$
and $\sigma$.

\noindent {\bf Theorem 1.4.} {\sl Let $M$ be a compact connected orientable
3-manifold, with $\partial M$ a (possibly empty) union of tori.
Let $K$ be a knot in $M$ such that $M - {\rm int}({\cal N}(K))$
is irreducible and atoroidal, and with 
$H_2(M - {\rm int}({\cal N}(K)), \partial M) \not= 0$.
Let $\sigma$ be an exceptional or norm-exceptional slope on
$\partial {\cal N}(K)$, such that $\Delta(\sigma, \mu) > 1$, where
$\mu$ is the meridian slope on $\partial {\cal N}(K)$. Then,
for a given $M$, there are at most finitely many possibilities for
$K$ and $\sigma$ up to ambient isotopy, and there is an algorithm to find them all.}

The algorithm is surprisingly simple. We describe it in Section 2.
The input into the algorithm is any triangulation of $M$,
or the following generalisation of a triangulation.
Let $P$ be a (possibly empty) collection of components of
$\partial M$. Then a {\sl generalised triangulation}
of $M$ is a representation of $M - P$
as a union of 3-simplices, with some or all of their faces identified in
pairs and then possibly with some subcomplex removed.
For example, an {\sl ideal triangulation} is the
case where $P = \partial M$ and where the subcomplex
removed is the 0-cells. We will also refer to the case
where $P = \emptyset$ as a {\sl genuine triangulation}.

There is a yet simpler algorithm which
deals with the $\sigma$-cable of $K$, which is defined to be
the knot in $M$ lying on $\partial {\cal N}(K)$ having slope $\sigma$.
Recall that a {\sl tangle} is a (possibly empty) collection
of disjoint arcs properly embedded in a 3-ball. Two tangles
are identified if there is an isotopy of the 3-ball which
is fixed on the boundary and which takes one tangle to
the other.

\noindent {\bf Theorem 1.5.} {\sl
There is a finite collection of
tangles, each lying in a 3-simplex and
with the following property. 
Let $M$, $K$ and $\mu$ be as in Theorem 1.4, and
let $\sigma$ be an exceptional slope on $\partial {\cal N}(K)$
with $\Delta(\sigma, \mu) > 1$.
Pick any generalised triangulation of $M$.
Then, we may insert a tangle from this
finite collection into each 3-simplex, in such a way that the tangles
join to form a knot which is ambient isotopic to the
$\sigma$-cable of $K$.
This finite collection of tangles
is constructible and is independent of $M$,
$K$ and $\sigma$.
}

Since these tangles are defined up to isotopy of the
3-simplex $\Delta^3$ which is fixed on $\partial \Delta^3$, 
Theorem 1.5 immediately gives that
there are only finitely many possibilities in $M$ for
the $\sigma$-cable of $K$.

Theorem 1.5 is a very surprising result. If $\Delta(\sigma,\mu)$ is large,
then one would expect the $\sigma$-cable of $K$ to intersect
the triangulation of $M$ in a complicated way. But the above
result asserts that one can control this complexity. It is
also surprising that the same finite collection of tangles
should work for all $M$ and all triangulations.
Note that in Theorem 1.5 we did not assume that
$K$ and $\sigma$ were norm-exceptional. In this case, we have the
following slightly weaker result.

\noindent {\bf Theorem 1.6.} {\sl Let $M$, $K$ and $\mu$ be
as in Theorem 1.4, and let $\sigma$ be a norm-exceptional slope
on $\partial {\cal N}(K)$ with $\Delta(\sigma, \mu) > 1$. If $M$ is closed,
pick any genuine triangulation of $M$.
In the case where $M$ has non-empty boundary,
pick any ideal triangulation of $M$. Then,
we may insert into each 3-simplex a tangle from the
finite collection of Theorem 1.5, in such a way that the tangles
join to form a knot which is ambient isotopic to the
$\sigma$-cable of $K$.}

It is in fact possible to write down explicitly this list
of tangles. We will give an algorithm in Section 11 to do this.
We have not actually run this algorithm, since the
task is fairly lengthy and is more suited to computer
implementation.

\vskip 18pt
\centerline {\caps 2. The algorithm to find all possibilities 
for $K$ and $\sigma$}
\vskip 6pt

In this section, we describe the algorithm for finding, in
a given 3-manifold $M$, all surgery curves $K$ with 
exceptional or norm-exceptional surgery slopes $\sigma$, 
satisfying the conditions of Theorem 1.4.
The first (and most important) step is to construct
a finite list of possibilities for $K$ and $\sigma$,
some of which may be neither norm-exceptional nor exceptional.

We will in Section 11 construct a finite collection
of graphs, each embedded in a 3-simplex $\Delta^3$.
Each graph $G$ meets $\partial \Delta^3$ in a collection
of vertices. These vertices have valence one and lie in
the interior of the 2-simplices of $\partial \Delta^3$. There is also a
specified regular neighbourhood ${\cal N}(G)$ and a
collection of disjoint arcs properly embedded in $\Delta^3$,
lying in $\partial {\cal N}(G)$. Each arc is assigned one of
two labels, $\gamma$ or $\tau$. Each graph $G$ (together
with ${\cal N}(G)$ and the arcs $\gamma$ and $\tau$)
is defined up to isotopy of $\Delta^3$ which is fixed
on $\partial \Delta^3$.

We will show during the course of the paper that it
is possible to ambient isotope $K$ and $\sigma$,
and to find a handle structure ${\cal H}$ on ${\cal N}(K)$
with the following properties.
Each tetrahedron $\Delta^3$ of the generalised triangulation of $M$
intersects the 0-handles and 1-handles of ${\cal H}$ in
${\cal N}(G)$, where $G$ is one of the graphs mentioned above.
The 2-handles of ${\cal H}$ will be
attached to ${\cal N}(G)$ along the arcs $\tau$.
There will also be a curve of slope $\sigma$ on $\partial {\cal N}(K)$
which intersects $\Delta^3$ in the arcs $\gamma$.

The algorithm to find all possibilities for $K$ and $\sigma$
therefore proceeds as follows. We insert one of these graphs 
into each 3-simplex of the generalised triangulation of $M$. 
If $\Delta^2$ is any 2-simplex of $M$ adjacent to two 3-simplices, 
and $G_1$ and $G_2$ are the graphs inserted into these two 3-simplices, 
then we insist that ${\cal N}(G_1) \cap \Delta^2 
= {\cal N}(G_2) \cap \Delta^2$, and
also that the endpoints in $\Delta^2$ of the arcs labelled $\gamma$ 
(respectively, $\tau$) in $G_1$ correspond precisely with
the endpoints in $\Delta^2$ of the arcs labelled $\gamma$
(respectively, $\tau$) in $G_2$. Thus, the graphs $G$
combine to form a graph (which we also
call $G$) embedded in $M$. The collections of arcs
combine to form a collection of curves 
$\gamma$ and $\tau$ properly embedded in $M$ 
and lying in $\partial {\cal N}(G)$.
We insist that each component of $\gamma$ and $\tau$ is a
simple closed curve, and that $\gamma$ is connected. 
Since there are only finitely
many 3-simplices in the representation of $M$ and
there are only finitely many possibilities for
${\cal N}(G) \cap \Delta^3$, $\tau \cap \Delta^3$
and $\gamma \cap \Delta^3$ for each 3-simplex $\Delta^3$
in $M$, there are only finitely many possibilities for
${\cal N}(G)$, $\tau$ and $\gamma$. The handlebody
${\cal N}(G)$ and curves $\tau$ specify a 
handle structure of a 3-manifold $M'$, which
is a candidate for ${\cal N}(K)$. 
At this stage, $M'$ may be something other than a solid torus.

The algorithm proceeds by calculating the
first homology of ${\cal N}(G)$, quotiented by
the subgroup generated by the curves $\tau$.
If the resulting homology $H_1(M')$ is not isomorphic
to ${\Bbb Z}$, we stop. If $H_1(M')$ is isomorphic to ${\Bbb Z}$,
then this implies that $H_1(\partial {\cal N}(M')) \cong {\Bbb Z}
\oplus {\Bbb Z}$. We can algorithmically find generators
$\lambda$ and $\mu$ of $H_1(\partial {\cal N}(M'))$ such
that $\lambda$ maps to $1 \in H_1(M')$, and $\mu$ maps to $0 \in H_1(M')$.
We can construct a simple closed curve representative of $\lambda$
on $\partial {\cal N}(M')$ which avoids the 2-handles of $M'$. 
If $M'$ is ${\cal N}(K)$, then this curve
is ambient isotopic in $M$ to $K$. The simple closed
curve $\gamma$ has slope $\sigma$. Thus, we have
constructed $K$ and $\sigma$.
If we wish, we can also calculate $\Delta(\lambda, \sigma)$
and $\Delta(\mu, \sigma)$. If $K$ is homologically trivial,
this (together with orientation information) gives
the rational number $p/q$ associated with $\sigma$.

The above algorithm constructs a finite number of
possibilities for $K$ and $\sigma$. We now wish to
rule out the cases where $K$ and $\sigma$ are
neither exceptional nor norm-exceptional. We construct the
manifold $M_K(\sigma)$. There is an algorithm to determine
whether $M_K(\sigma)$ is reducible ([5] and [11]), and there
is an algorithm to determine whether $M_K(\sigma)$ is a solid torus 
([5] and [11]). The assumption that $H_2(M - {\rm int}({\cal N}(K)), \partial M)$
is non-trivial implies that $H_1(M_K(\sigma))$ is infinite and hence
that $\pi_1(M_K(\sigma))$ is infinite. If $M_K(\sigma)$ is
irreducible, then according to Corollary 9.9 of [3],
$\pi_1(M_K(\sigma))$ is torsion-free. Thus, if the core
of the surgery solid torus in $M_K(\sigma)$ has finite order
in $\pi_1(M_K(\sigma))$, then it is homotopically trivial. 
There is an algorithm to determine this, since
the word problem is soluble for
the fundamental groups of Haken 3-manifolds [13]. Finally,
there is an algorithm to find the unit ball of the
Thurston norm (Algorithm 5.9 of [12]), 
and so we can determine whether $\sigma$ is norm-exceptional.

\vfill\eject
\centerline {\caps 3. The sutured manifold theory background}
\vskip 6pt

The definition of an exceptional surgery was specifically
designed so that sutured manifold theory can be applied.
Sutured manifolds were defined and studied by Gabai [1]
who used them to construct taut foliations on certain
3-manifolds. In this section, we will outline a version of 
the theory due to Scharlemann [10]. Almost everything in this
section can be found elsewhere, mostly in Scharlemann's
paper [10]. We include it here because it is absolutely
central to our argument, but a reader familiar with the
theory of sutured manifolds may safely skip this section.

Sutured manifold theory is intimately linked to the
Thurston norm. Here, we give a definition of the
Thurston norm and some related definitions of
tautness.

Let $S$ be a compact oriented
surface embedded in a compact oriented 3-manifold $M$.
Let $\chi(S)$ denote its Euler characteristic. If $S$ is connected, define
$$\chi_-(S) = \max \lbrace 0, -\chi(S) \rbrace.$$
When $S$ is not connected, define $\chi_-(S)$ to be the sum of
$\chi_-(S_i)$ over all the components $S_i$ of $S$.

Let $P$ be a subset of $\partial M$, and
let $z$ be an element of $H_2(M, P)$
which is represented by some embedded compact oriented surface. 
Then the {\sl Thurston norm} of $z$ is given by
$$x(z) = \min \lbrace \chi_-(S) :
\hbox{$S$ is an embedded surface representing $z$} \rbrace.$$

Let $(S, \partial S) \subset (M, \partial M)$ be an
oriented compact surface embedded in $M$. 
Let $P$ be a subset of $\partial M$ which
contains $\partial S$.
Then $S$ is {\sl norm-minimising} in $H_2(M, P)$ if
$x([S, \partial S]) = \chi_-(S)$.
In the case where $P = \partial S$,
then $S$ is {\sl taut} if it is incompressible and norm
minimising in $H_2(M, P)$.
This use of the word `taut' is not entirely standard; some
authors (for example, [12]) 
insist only that $S$ be incompressible and norm-minimising
in its class in $H_2(M, \partial M)$. However, our
definition is more suited to sutured manifold theory.

A {\sl sutured manifold} $(M, \gamma)$ is
a compact oriented 3-manifold $M$, with $\partial M$
decomposed into two subsurfaces ${\cal R}_-$ and ${\cal R}_+$, such that
${\cal R}_- \cup {\cal R}_+ = \partial M$ and 
${\cal R}_- \cap {\cal R}_+ = \gamma$, where $\gamma$
is a union of disjoint
simple closed curves, known as the {\sl sutures}. The subsurfaces
${\cal R}_-$ and ${\cal R}_+$ are oriented so that the normal vectors of
${\cal R}_-$ (respectively, ${\cal R}_+$) point into (respectively, out of) $M$.
The symbol ${\cal R}_\pm$ will denote `${\cal R}_-$ or ${\cal R}_+$'.
When we wish to emphasise a particular sutured manifold,
we will use the symbol ${\cal R}_\pm(M)$.

A sutured manifold $(M, \gamma)$ is
{\sl taut} if $M$ is irreducible, and ${\cal R}_-$ and 
${\cal R}_+$ are taut.
For example it is not hard to show the following.
Suppose that $\partial M$ is a (possibly empty)
union of tori, and that ${\cal R}_- = \partial M$ and ${\cal R}_+ = \emptyset$.
Then $(M, \emptyset)$ is taut if and only is $M$
is neither reducible nor a solid torus.

One of the main techniques of the theory
is to decompose a sutured manifold along a properly embedded
surface. If $(M, \gamma)$ is a sutured manifold, and $S$ is an 
oriented surface properly embedded in $M$, with $\partial S$ and
$\gamma$ in general position, then $M_S = M - {\rm int}({\cal N}(S))$
inherits a sutured manifold structure $(M_S, \gamma_S)$.
This is written $(M, \gamma) \buildrel S \over \longrightarrow
(M_S, \gamma_S)$. This decomposition is said to be {\sl taut}
if $(M, \gamma)$ and $(M_S, \gamma_S)$ are both taut.

If $(M, \gamma)$ is a connected taut sutured manifold and $z$ is any non-zero
homology class in $H_2(M, \partial M)$, then (Theorem 2.6 of [10]) there is
a taut decomposition $(M, \gamma) \buildrel S \over \longrightarrow
(M_S, \gamma_S)$ such that
\item{(i)} no curve of $\partial S$ bounds a disc in ${\cal R}_\pm(M)$,
\item{(ii)} no component $X$ of $M_S$ has $\partial X \subset {\cal R}_-(M_S)$
or $\partial X \subset {\cal R}_+(M_S)$, and
\item{(iii)} $[S, \partial S] = z \in H_2(M, \partial M)$.

\noindent When $S$ satisfies (i), we will say that $\partial S$ {\sl has
essential intersection with} ${\cal R}_\pm(M)$. It is not hard to show
that if $(M, \gamma) \buildrel S \over \longrightarrow
(M_S, \gamma_S)$ is a taut decomposition and $\partial S$ 
has essential intersection with ${\cal R}_\pm(M)$,
then $S$ itself is taut.

Thus if $H_2(M, \partial M) \not= 0$, we may perform a
taut sutured manifold decomposition along a taut surface
having essential intersection with ${\cal R}_\pm$. But if $H_2(M, \partial M)$
is trivial,
then it is a classical fact that $\partial M$ is a (possibly empty) union
of 2-spheres. If $M$ is irreducible, then this implies that either
$\partial M = \emptyset$ or $M$ is a collection of 3-balls.
Using this argument, Gabai proved that, if $(M ,\gamma)$ is
a connected taut sutured manifold and $H_2(M, \partial M) \not= 0$, then
there is a sequence of taut decompositions
$$(M, \gamma) = (M_1, \gamma_1) \buildrel S_1 \over \longrightarrow 
(M_2, \gamma_2)
\buildrel S_2 \over \longrightarrow \dots \dots \buildrel
S_{n-1} \over \longrightarrow (M_n, \gamma_n),$$
with $M_n$ a union of 3-balls.

An important step in Gabai's argument is to show that this sequence of
decompositions cannot continue indefinitely. This is
not at all obvious. In the case where $S_i$ is
incompressible and $\partial$-incompressible in
$M_i$ and has no component parallel to a subsurface of
$\partial M_i$, it was proved by Haken [2] that such a sequence
of decompositions must eventually terminate. However,
it is sometimes necessary to consider surfaces $S_i$ which
are $\partial$-compressible. Nevertheless, Gabai constructed
a (complicated) argument which proved that this sequence
of taut sutured manifold decompositions can be guaranteed
to terminate. He did this by defining a complexity of
a sutured manifold and then arguing by induction.
In Section 5, we will offer a new definition of complexity
for a sutured manifold with a given handle decomposition.

There is an extremely useful property of sutured manifold
decompositions, which is summarised in the phrase
`tautness usually pulls back'. It is this property which
makes sutured manifold theory distinctly different from
the theory of Haken manifolds.

\noindent {\bf Theorem 3.1.} (Theorem 3.6 of [10]) {\sl Let $(M, \gamma) \buildrel S \over
\longrightarrow (M_S, \gamma_S)$ be a decomposition,
where $\partial S$ has essential intersection with ${\cal R}_\pm(M)$,
and where no component of $S$ is a compression disc for a torus
component of ${\cal R}_\pm(M)$.
If $(M_S, \gamma_S)$ is taut, then so is $(M, \gamma)$.}

There is a partial converse to this theorem which can
be useful. If $D$ is a disc
properly embedded in $M$ intersecting $\gamma$ transversely
in two points, then $D$ is known as a {\sl product disc}.
If $A$ is an annulus properly embedded in $M$
with one component of $\partial A$ in ${\cal R}_-$
and one in ${\cal R}_+$, then $A$ is known as a
{\sl product annulus}. These surfaces play a useful
r\^ole, since if $(M, \gamma) \buildrel S \over \longrightarrow
(M_S, \gamma_S)$ is a decomposition along a product disc
or an incompressible product annulus, then $(M, \gamma)$
is taut if and only if $(M_S, \gamma_S)$ is taut.

We have now described enough sutured manifold theory
to explain the definition of an exceptional surgery
curve, given in Section 1. The following argument is
well known, and is due to Gabai [1]. Let $K$ be a knot in a compact connected 
orientable 3-manifold $M$, where $\partial M$ is (possibly empty) union of tori.
If $M - {\rm int}({\cal N}(K))$
is neither reducible nor a solid torus, then 
$(M - {\rm int}({\cal N}(K)), \emptyset)$
is a taut sutured manifold, with ${\cal R}_- = \partial M \cup
\partial {\cal N}(K)$. If $H_2(M - {\rm int}({\cal N}(K)), 
\partial M) \not= 0$, then we may perform a 
taut sutured manifold decomposition 
$$(M - {\rm int}({\cal N}(K))) \buildrel S_1 \over \longrightarrow
(M_2 - {\rm int}({\cal N}(K))),$$
such that
\item{$\bullet$} $S_1$ is disjoint from $\partial {\cal N}(K)$,
\item{$\bullet$} no simple closed curve of $\partial S_1$ bounds a disc
in ${\cal R}_\pm(M)$, and
\item{$\bullet$} no component $X$ of $M_2$ has
$\partial X \subset {\cal R}_-(M_2)$ or $\partial X \subset {\cal R}_+(M_2)$. 

\noindent If $K$ is norm-exceptional, we insist that $[S_1, \partial S_1]
= z \in H_2(M - {\rm int}({\cal N}(K)), \partial M)$, where
$z$ is the relevant homology class from Definition 1.3. Repeating
this process, we construct a sequence of taut
sutured manifold decompositions
$$(M - {\rm int}({\cal N}(K)), \emptyset)
\buildrel S_1 \over \longrightarrow \dots
\buildrel S_{n-1} \over \longrightarrow 
(M_n - {\rm int}({\cal N}(K)), \gamma_n),$$
satisfying the above conditions, and 
where $H_2(M_n - {\rm int}({\cal N}(K)), \partial M_n) = 0$.
If $M - {\rm int}({\cal N}(K))$ is atoroidal, then
it is possible to show that this implies that $M_n$
is a solid torus regular neighbourhood of $K$ and possibly some 3-balls.
No component $X$ of $M_n$ has $\partial X \subset {\cal R}_-(M_n)$
or $\partial X \subset {\cal R}_+(M_n)$. 
In particular, if $X$ is the component of $M_n$ containing
$K$, then $\gamma_n \cap X$ is
a collection of essential curves on $\partial X$,
parallel to some slope $\rho$, say, on $\partial {\cal N}(K)$.
If we Dehn fill $M-{\rm int}({\cal N}(K))$ via any slope
$\tau$ on $\partial {\cal N}(K)$, 
then $M_n - {\rm int}({\cal N}(K))$ is filled to
become a 3-manifold $M_n(\tau)$ which is a solid torus and
some 3-balls. Now, $M_n(\tau)$ inherits
a sutured manifold structure $(M_n(\tau), \gamma_n)$
from $M_n - {\rm int}({\cal N}(K))$,
which is taut if the surgery slope $\tau$ is not the 
slope $\rho$ of the sutures. Since tautness pulls
back, this implies that
$$(M_K(\tau), \emptyset) \buildrel S_1 \over \longrightarrow
\dots \buildrel S_{n-1} \over \longrightarrow (M_n(\tau), \gamma_n)$$
is a sequence of taut sutured manifolds, with each $S_i$
taut in $M_i(\tau)$. This implies that
\item{(i)} $M_K(\tau)$ is irreducible,
\item{(ii)} $M_K(\tau)$ is not a solid torus,
\item{(iii)} the core of the surgery solid torus in
$M_K(\tau)$ has infinite order in $\pi_1(M_K(\tau))$, and
\item{(iv)} $S_1$ is taut in $M_K(\tau)$.

\noindent Now, if $\sigma$ is an exceptional or norm-exceptional surgery slope
on $\partial {\cal N}(K)$, then at least one of the
above cannot be true for $M_K(\sigma)$. Thus,
$\sigma$ must be the slope $\rho$ which is parallel to the sutures in
$M_n$. We assume in Theorems 1.4, 1.5 and 1.6 that $\Delta(\sigma, \mu) > 1$
which implies in particular that $\sigma \not= \mu$. Thus,
the facts (i) - (iv) above are true for $\tau = \mu$, and also
$$(M, \emptyset) \buildrel S_1 \over \longrightarrow
\dots \buildrel S_{n-1} \over \longrightarrow (M_n, \gamma_n)$$
is a taut sutured manifold sequence. Each component of
$\gamma_n$ lies inside $M$ as the $\sigma$-cable
of $K$, or as an unknotted curve.
The idea behind Theorems 1.4, 1.5 and 1.6 is (roughly speaking)
inductively to find nice embeddings of $M_i$ in $M$. In particular,
we will show that we can arrange that $\gamma_n \cap \Delta^3$
is one of a finite list of possibilities for each 3-simplex
$\Delta^3$ of $M$. Since one component of $\gamma_n$
is the $\sigma$-cable of $K$, this will establish
Theorem 1.5.

Thus, our definition of an exceptional surgery curve
fits neatly into the sutured manifold setting.
The sutured manifold theory which we have outlined
above formed the basis for a theorem in [6] which
will be an important technical tool in this paper.
This result (Theorem 1.4 of [6]) deals with the
interaction of exceptional surgery curves and
embedded surfaces in a sutured manifold, and is given below.

\noindent {\bf Theorem 3.2.} {\sl Let $(M, \gamma)$ be a
taut sutured manifold, let $K$ be a knot in $M$ and
let $\sigma$ be a slope on $\partial {\cal N}(K)$.
Suppose that at least one of the following is true:
\item{(i)} $\sigma$ is an exceptional surgery slope, or
\item{(ii)} $\sigma$ is a norm-exceptional surgery slope,
$\partial M$ is a (possibly empty) union of tori
and $\gamma = \emptyset$.

\noindent Suppose that $\Delta(\sigma, \mu) > 1$, where
$\mu$ is the meridian slope on $\partial {\cal N}(K)$.
Suppose also that $M - {\rm int}({\cal N}(K))$ is
irreducible and atoroidal and that $H_2(M - {\rm int}({\cal N}(K)),
\partial M) \not= 0$.
Let $F$ be a surface properly embedded in $M$, with
components $F_1, \dots, F_n$, none of which is a
sphere or disc disjoint from $\gamma$. Then there
is an ambient isotopy of $K$ in $M$, after which we have the following
inequality for each $i$:
$$\vert K \cap F_i \vert \leq {-2\chi(F_i) + \vert \gamma \cap F_i \vert
\over 2(\Delta(\sigma, \mu) - 1)}.$$}

The numerator $-2\chi(F_i) + \vert \gamma \cap F_i \vert$
is known as the {\sl index} $I(F_i)$ of $F_i$. 
Note in particular that a product disc
and an annulus disjoint from $\gamma$ both have
zero index. Thus, if $\Delta(\sigma, \mu) > 1$,
Theorem 3.2 implies that we may ambient isotope
$K$ off a collection of product discs and annuli
disjoint from $\gamma$. Given that such surfaces
play a useful technical r\^ole in sutured manifold theory,
this will be very convenient. In fact, this is the only
point in proof of Theorems 1.4, 1.5 and 1.6 where we use
that $\Delta(\sigma, \mu) > 1$.

\vskip 18pt
\centerline {\caps 4. Vertical form and standard form for submanifolds.}
\vskip 6pt

Recall that we are given a generalised triangulation
of $M$. From this, we will construct the dual handle
decomposition, which associates an $i$-handle with each
$(3-i)$-simplex of $M$ not lying entirely in $\partial M$. 
For this dual handle decomposition,
the boundary of each 0-handle has at most four 
discs of intersection with the 1-handles, and each
1-handle has at most three discs of intersection with the
2-handles. An example is given below.

\vskip 18pt
\centerline{\psfig{figure=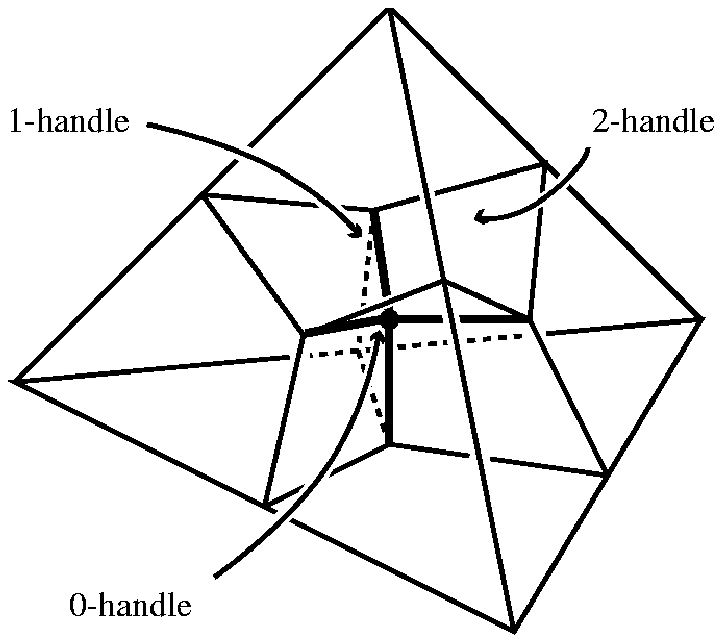}}
\vskip 6pt
\centerline{Figure 4.1.}

We will now give some definitions and conventions regarding
handle decompositions. We will throughout this paper
denote the $i$-handles of a handle decomposition by ${\cal H}^i$.
Henceforth, we will only consider handle decompositions
of $n$-manifolds with the following properties:
\item{$\bullet$} for $i > 0$, the $i$-handles are
attached to $\bigcup_{j < i}{\cal H}^j$, and
\item{$\bullet$} if $H_i = D^{n-i} \times D^i$ 
(respectively, $H_j = D^{n-j} \times
D^j$) is an $i$-handle (respectively, $j$-handle) with $j < i$, then
$H_i \cap H_j = E \times D^j = D^{n-i} \times F$ for some submanifold $E$
(respectively, $F$) of $\partial D^{n-j}$ (respectively, $\partial D^i$).

\noindent In words, the second condition requires that
the attaching map of each handle respects the product
structures of the handles to which it is attached. For a 3-manifold, 
this is relevant only for $j=1$ and $i=2$.
In the case of a handle decomposition of a 3-manifold,
we also insist that

\item{$\bullet$} no 2-handle is disjoint from ${\cal H}^1$.

\noindent We will use the term {\sl handle structure} for
a decomposition satisfying these conditions. Note 
in particular that the dual handle decomposition
of a 3-manifold arising from a generalised triangulation has these
properties. 
We will use ${\cal H}$ to denote a handle structure, but occasionally,
we will also write ${\cal H}(M)$ when we wish to emphasise the
manifold $M$. Note also that a handle structure ${\cal H}$ on a 
3-manifold $M$ induces
a handle structure on $\partial M$, which we will usually write
as ${\cal H}(\partial M)$.

We will fix a handle structure ${\cal H}$ of $M$, and then will consider
embedded submanifolds of $M$. We wish to ensure that
each submanifold lies inside ${\cal H}$ in a manageable way.
The relevant notions are `vertical' and `standard' form,
the first of which we now define.

\noindent {\bf Definition 4.2.} Let $M$ be an $n$-manifold with
a handle structure ${\cal H}$. Let $S$ be an $(n-1)$-manifold
properly embedded in $M$. Then $S$ is in {\sl vertical form}
if, for each $i$-handle $D^{n-i} \times D^i$ of ${\cal H}$, we have
$S \cap (D^{n-i} \times D^i) = E \times D^i$,
where $E$ is a properly embedded submanifold of
$D^{n-i}$. In particular, $S$ is disjoint from ${\cal H}^n$.

\vfill\eject
\centerline{\psfig{figure=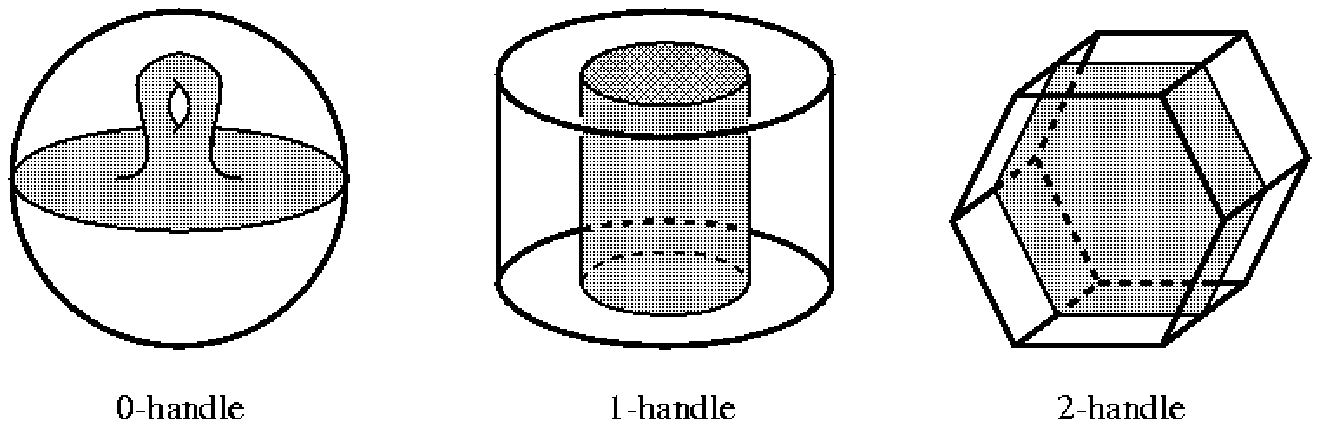}}
\vskip 6pt
\centerline{Figure 4.3.}
\vskip 6pt

The only two cases which we will consider are where $n=2$ or $n=3$.
Examples of 2-manifolds in vertical form in a 3-manifold are given in
Fig. 4.3. The relevance of vertical form is its ubiquity.

\noindent {\bf Lemma 4.4.} {\sl Let $M$ be an $n$-manifold with
a handle structure ${\cal H}$, and let $S$
be an $(n-1)$-manifold properly embedded in $M$.
Then there is an ambient isotopy which takes $S$ into
vertical form with respect to ${\cal H}$.}

\noindent {\sl Proof.} We perform a sequence of ambient
isotopies. The first pulls $S$ off ${\cal H}^n$. The
second places $S$ in vertical form in ${\cal H}^{n-1}$,
and so on. Let $C_i$ be the co-cores of
${\cal H}^i$; thus ${\cal H}^i = C_i \times D^i$. We perform an ambient isotopy
which makes $S$ transverse to $C_i$. By construction,
$S$ is already vertical in ${\cal H}^j$ for $j > i$,
and so we may take this isotopy to be supported
in ${\cal H}^i - \partial {\cal H}^i$.
After the isotopy, we may find a small disc $D_0^i \subset {\rm int}(D^i)$,
such that $S \cap (C_i \times D_0^i) = E_i \times D_0^i$,
for some submanifold $E_i$ of $C_i$. Then we may use the 
product structure on $D^i - {\rm int}(D_0^i) \cong S^{i-1} \times I$
to ambient isotope $C_i \times D_0^i$ onto $C_i \times D^i =
{\cal H}^i$. We can take this isotopy of $M$ to be supported in
an arbitrarily small neighbourhood of ${\cal H}^i$,
and also to leave $S \cap {\cal H}^j$ invariant for
$j > i$. After performing these isotopies for
$i = n, n-1, \dots, 0$, we finish with $S$ in vertical form.
$\square$

If $(M, \gamma)$ is a sutured manifold with a handle structure
${\cal H}$, then by Lemma 4.4 there is an isotopy of
$\partial M$ which takes $\gamma$ into vertical form
(with respect to the induced handle decomposition 
${\cal H}(\partial M)$ on $\partial M$).
This isotopy of $\partial M$ extends to an isotopy
of $M$. We can therefore
assume that $\gamma$ is in vertical form in ${\cal H}(\partial M)$,
and we will henceforth make this assumption.

If $S$ is any surface properly embedded
in $M$, we would like to ensure that we can place $S$ in vertical form,
and still keep $\gamma$ vertical in ${\cal H}(\partial M)$. This is the purpose
of the following lemma.

\noindent {\bf Lemma 4.5.} {\sl Let $(M, \gamma)$ be a
sutured manifold with a handle structure ${\cal H}$,
such that $\gamma$ is vertical in ${\cal H}(\partial M)$.
If $S$ is a surface properly embedded in $M$, in general
position with respect to $\gamma$, then
there is an ambient isotopy which leaves
$\gamma$ invariant and which moves $S$ into vertical
form.}

\noindent {\sl Proof.} The first two steps of the ambient isotopy in
Lemma 4.4 are supported in a small neighbourhood of ${\cal H}^3 \cup {\cal H}^2$.
Hence, we may assume that it leaves $\gamma$ fixed.
Since $S$ and $\gamma$ are in general position,
we may pick the co-cores $C_1$ of the 1-handles so 
that $C_1 \cap S \cap \gamma = \emptyset$. The 
ambient isotopy supported in a neighbourhood of ${\cal H}^1$
can then be taken to leave $\gamma$ invariant.
There is no restriction on $S \cap {\cal H}^0$,
once $S$ lies in the remaining handles in the correct way.
Hence, we have ambient isotoped $S$ into vertical form,
leaving $\gamma$ invariant. $\square$

For inductive purposes, we define a notion of
complexity for surfaces in vertical form in a handle
structure of a 3-manifold.

\noindent {\bf Definition 4.6.} The {\sl complexity} of a vertical surface
$S$ is the ordered pair of integers $(\vert S \cap {\cal H}_2
\vert, \vert \partial S \cap {\cal H}_1 \vert)$.

We order these pairs lexicographically. In other words,
the pairs $(n_1, n_2)$ and $(m_1, m_2)$ satisfy
$(n_1, n_2) > (m_1, m_2)$ precisely when
\item{$\bullet$} $n_1 > m_1$, or
\item{$\bullet$} $n_1 = m_1$ and $n_2 > m_2$.

\noindent It is clear that this ordering is a well-ordering.

In the case of surfaces in 3-manifolds, there is a notion which is
a little stronger than vertical form.

\noindent {\bf Definition 4.7.} Let $S$ be a vertical surface in
a handle structure ${\cal H}$ of a 3-manifold $M$.
Then $S$ is {\sl standard} if $S$ intersects each
handle of ${\cal H}$ in a (possibly empty) collection of discs.

\vskip 24pt
\centerline{\psfig{figure=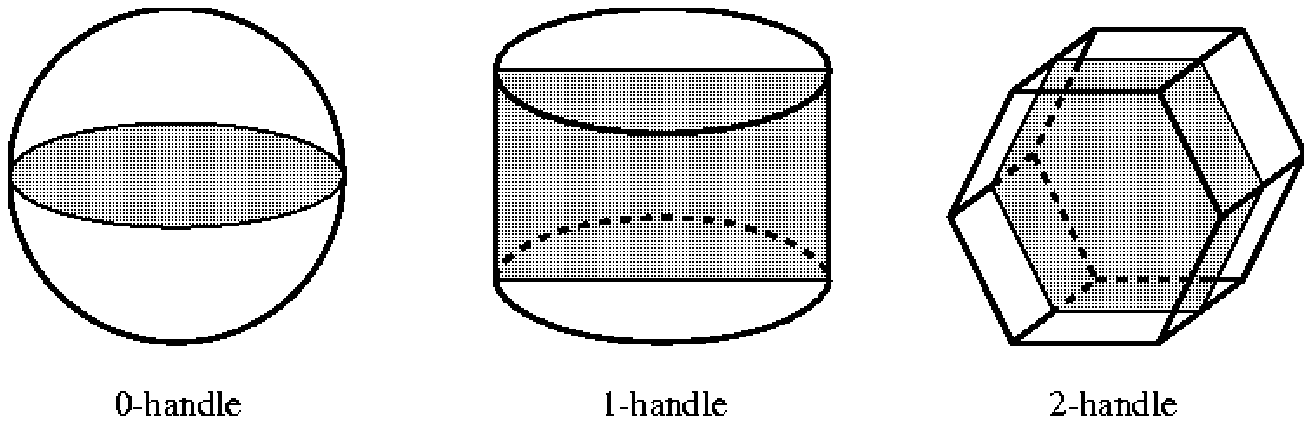}}
\vskip 6pt
\centerline{Figure 4.8.}
\vskip 6pt

Examples of surfaces in standard form are given in
Fig 4.8. A general surface $S$ in $M$ might not have
a representation in standard form, but if $S$ is incompressible
and $M$ is irreducible, then we now show that it can be
ambient isotoped into standard form.

\noindent {\bf Lemma 4.9.} {\sl Let $(M, \gamma)$ be an
irreducible sutured manifold with
a handle structure ${\cal H}$. Let $S$
be a vertical incompressible surface properly embedded in $M$,
with no component of $S$ a 2-sphere.
Then there is an ambient isotopy of $S$ which leaves
$\gamma$ fixed and which takes $S$ into standard
form without increasing its complexity.}

\noindent {\sl Proof.} If $S$
is not in standard form, then it must differ from standard
form in some 1-handle or some 0-handle of ${\cal H}$.
Suppose first that, in some 1-handle $H_1 = D^2 \times D^1$,
there is a component of $S \cap H_1$ which is $\alpha \times D^1$,
for a simple closed curve $\alpha$.  
If both curves of $\alpha \times \partial D^1$ bound 
discs in ${\cal H}^0$, then $S$ has a 2-sphere
component. Hence, we may assume that
$S$ differs from standard form in some 0-handle $H_0$.
That is, suppose that $S \cap H_0$ is not a union of discs. Then, 
since no component of $S$ is a 2-sphere, $S \cap H_0$
is compressible in $H_0$, via a compression disc $D$.
Since $S$ is incompressible, $\partial D$ bounds a disc
$D'$ in $S$. The disc $D'$ does not lie wholly in $H_0$,
and so must intersect ${\cal H}^1$.
As $M$ is irreducible, we may ambient isotope $S$,
taking $D'$ onto $D$. This does not
increase the complexity of $S$, and it
reduces the number of components of $S \cap {\cal H}^1$.
Hence, this process terminates with $S$ in standard form.
The isotopy at each stage leaves $\partial M$ (and hence
$\gamma$) fixed. $\square$

We may therefore assume that if $S$ and $M$ satisfy the
conditions of Lemma 4.9, then $S$ is
in standard form. We will now show that, if
$(M_S, \gamma_S)$ is the sutured manifold resulting from the
decomposition along $S$, then $M_S$ has an induced handle structure
with $\gamma_S$ in vertical form in ${\cal H}(\partial M_S)$.

If $H$ is an $i$-handle $D^{3-i} \times D^{i}$ of ${\cal H}(M)$,
then each component of $H - {\rm int}({\cal N}(S))$
inherits a structure $X \times D^i$,
where $X$ is a $(3-i)$-submanifold of $D^{3-i}$. This is true
because $S$ is vertical. Since $S$ is standard, then
each component of $X$ is a copy of $D^{3-i}$, and so each component of
$H - {\rm int}({\cal N}(S))$
has the structure of an $i$-handle. These handles combine to
give a handle structure on $M_S$.
The curves $\gamma_S$ are a subgraph of the graph
$\partial {\cal N}(\partial S) \cup \gamma$. Since
$\partial S$ and $\gamma$ are both vertical in ${\cal H}(\partial M)$,
the curves $\gamma_S$ are then vertical in ${\cal H}(\partial M_S)$.

It is a very useful property that
$(M_S, \gamma_S)$ inherits a handle structure from that of
$(M, \gamma)$. It is the basis for an
inductive proof of Theorems 1.4, 1.5 and 1.6. However, to
construct such a proof, we need to define a
`complexity' for handle structures. 

\vskip 18pt
\centerline {\caps 5. Complexity of handle structures of sutured manifolds}
\vskip 6pt

We will now define a notion of complexity for a
handle structure ${\cal H}$ of a sutured manifold $(M, \gamma)$.
We will focus on the 2-spheres $\partial {\cal H}^0$.
Lying in these 2-spheres, there is the surface
$\partial {\cal H}^0 \cap ({\cal H}^1
\cup {\cal H}^2)$. We denote this surface by ${\cal F}({\cal H})$,
or sometimes simply ${\cal F}$. 

\vskip 18pt
\centerline{\psfig{figure=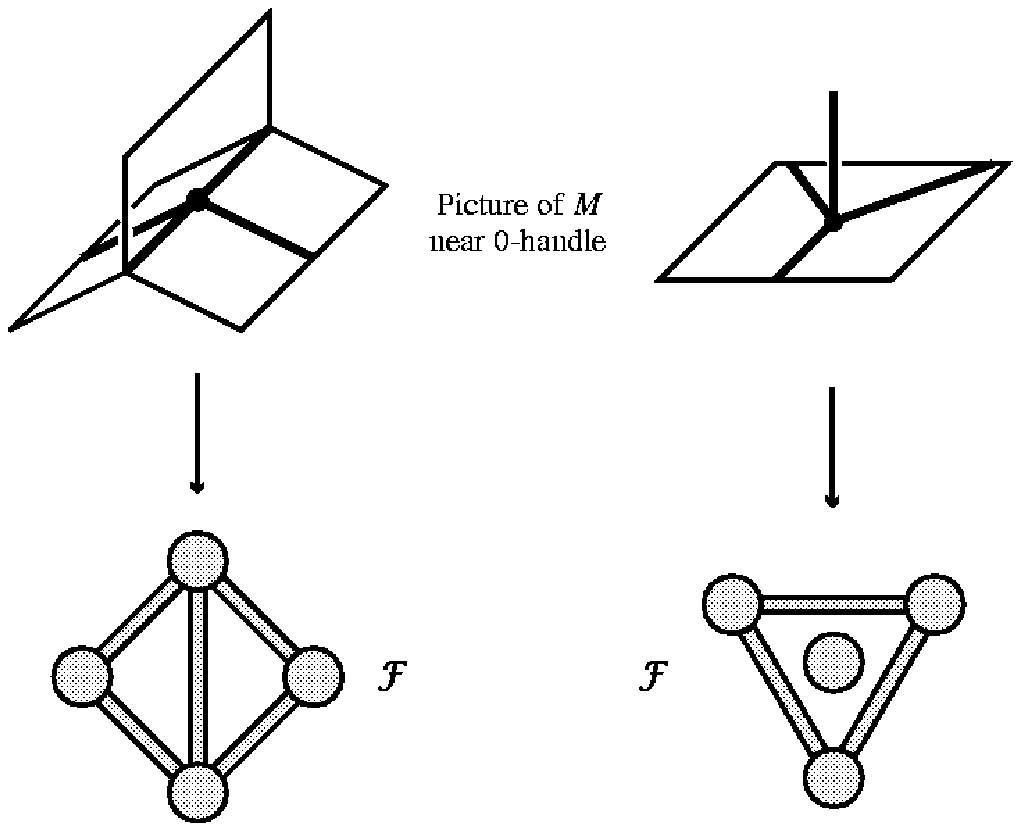}}
\nobreak
\vskip 6pt
\centerline{Figure 5.1.}
\vskip 6pt

Recall from Section 4 that we
insisted that no 2-handle of ${\cal H}$ is disjoint from
${\cal H}^1$. Therefore, ${\cal F}$ inherits a handle structure, with
$\partial {\cal H}^0 \cap {\cal H}^1$ forming the 0-handles of
${\cal F}$ (which we denote by ${\cal F}^0$),
and $\partial {\cal H}^0 \cap {\cal H}^2$
forming the 1-handles of ${\cal F}$ (which we denote by
${\cal F}^1$). Note that
each component of the surface ${\rm cl}(\partial {\cal H}^0 - {\cal F})$
lies either in $\partial M$ or in $\partial {\cal H}^3$,
and the curves $\gamma \cap \partial {\cal H}^0$ are
properly embedded in ${\rm cl} (\partial {\cal H}^0 - {\cal F})$.

If $S$ is in standard form, then the simple closed curves
$S \cap \partial {\cal H}^0$ satisfy the following
(fairly weak) restrictions:
\item{$\bullet$} $S \cap \partial {\cal H}^0$ is disjoint
from ${\cal H}^3$,
\item{$\bullet$} $S \cap {\cal F}$ is vertical in ${\cal F}$,
\item{$\bullet$} no curve of $S \cap \partial {\cal H}^0$ lies
entirely within a 0-handle of ${\cal F}$.

\noindent The nature of ${\cal F}$ will determine the complexity
of ${\cal H}$. One invariant of ${\cal F}$ will be of particular
importance, namely its index. Recall from Section 3
that the index $I(F)$ of a component $F$ of ${\cal F}$ is
defined to be
$$I(F) = -2\chi(F) + \vert F \cap \gamma \vert.$$
If $V$ is a 0-handle of ${\cal F}$, then the {\sl valence} of
$V$ is the number of arcs of $V \cap {\cal F}^1$.
We also define the {\sl index} of $V$ to be
$$I(V) = \vert V \cap {\cal F}^1 \vert + \vert V \cap \gamma \vert -
2.$$
The reason for this terminology is that
$$I(F) = \sum_{V \in F \cap {\cal F}^0} I(V).$$

For each component $F$
of ${\cal F}$, we define the following integers:
$$\eqalign{
C_1(F) &= \vert F \cap {\cal F}^1 \vert + 1,\cr
C_2(F) &= I(F),\cr
C_3(F) &= \vert \partial F \vert.}$$
The ${\cal F}$-{\sl complexity set} 
$C_{\cal F}({\cal H})$ of ${\cal H}$ is defined to
be the set of ordered triples
$$C_{\cal F}({\cal H}) = \lbrace (C_1(F), C_2(F), C_3(F)) :
\hbox{$F$ a component of ${\cal F}$ with $I(F) > 0$} \rbrace,$$
where repetitions are retained. If $X$ is a subset of $M$,
with $X \cap {\cal F}$ a non-empty collection of components
of ${\cal F}$, then we similarly define
$$C_{\cal F}(X) = \lbrace (C_1(F), C_2(F), C_3(F)) :
\hbox{$F$ a component of $X \cap {\cal F}$ with $I(F) > 0$} \rbrace,$$
where again repetitions are retained.

An example is given in Fig. 5.2 of how ${\cal F}$ and its
complexity behave when ${\cal H}$ is decomposed along
a surface $S$.

\vskip 26pt
\centerline{\psfig{figure=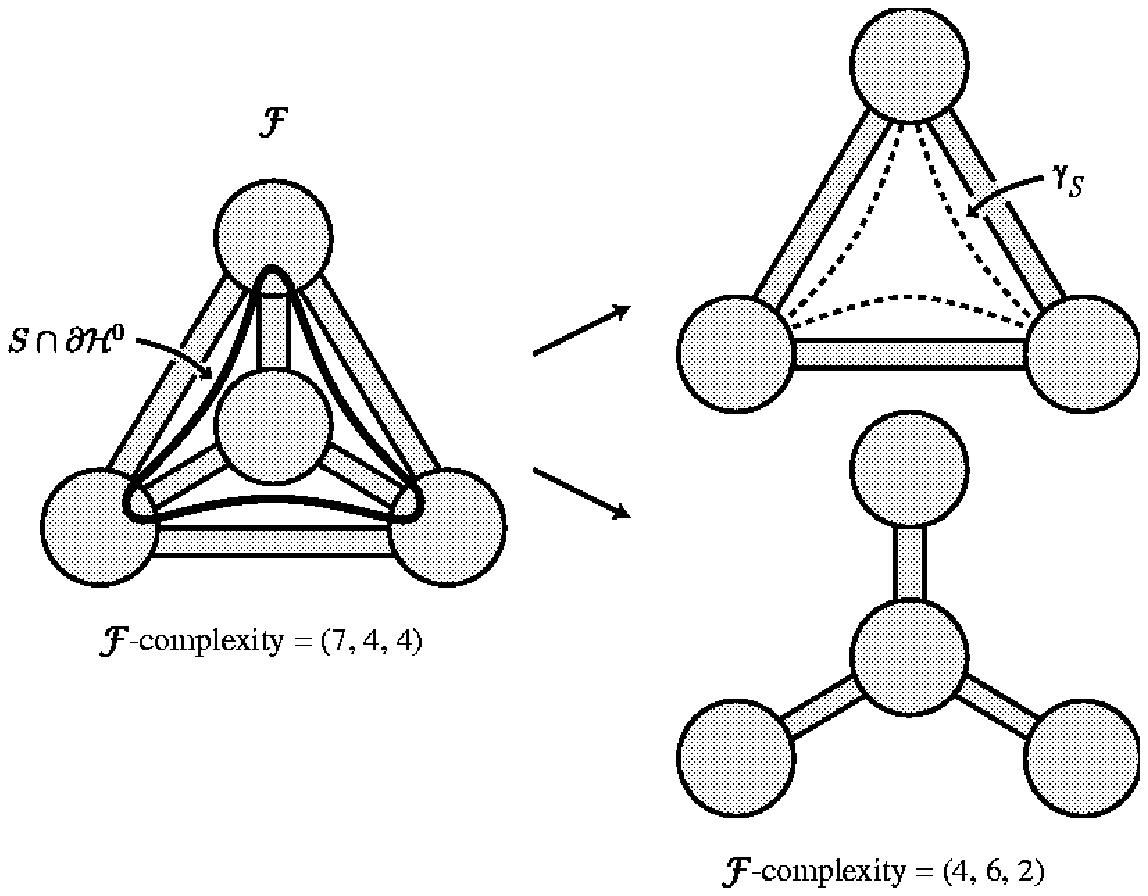}}
\vskip 6pt
\centerline{Figure 5.2.}
\vskip 6pt

We compare the triples
$(C_1(F), C_2(F), C_3(F))$ and $(C_1(F'), C_2(F'), C_3(F'))$
by defining
$(C_1(F), C_2(F), C_3(F)) > (C_1(F'), C_2(F'), C_3(F'))$
if
\item{$\bullet$} $C_1(F) > C_1(F')$, or
\item{$\bullet$} $C_1(F) = C_1(F')$ and $C_2(F) > C_2(F')$, or
\item{$\bullet$} $C_1(F) = C_1(F')$ and $C_2(F) = C_2(F')$
and $C_3(F) > C_3(F')$.

\noindent It is clear that this is a total ordering and
a well-ordering.

We define a total order on the ${\cal F}$-complexity of handle
structures, as follows. 
If ${\cal H}$ and ${\cal H}'$ are two handle structures,
we order their ${\cal F}$-complexity sets 
$C_{\cal F}({\cal H})$ and $C_{\cal F}({\cal H}')$
into two non-increasing sequences of triples. We extend each of
these sequences by concatenating with an infinite
sequence of triples $(0,0,0)$. (Note that always $C_1(F) > 0$,
and so $(C_1(F), C_2(F), C_3(F)) > (0,0,0)$.)
Then, we compare the first (and hence largest) triple
$(C_1(F), C_2(F), C_3(F))$ of
$C_{\cal F}({\cal H})$ with the first (and hence largest) triple
$(C_1(F'), C_2(F'), C_3(F'))$ of $C_{\cal F}({\cal H}')$.
If $(C_1(F), C_2(F), C_3(F)) > (C_1(F'), C_2(F'), C_3(F'))$,
say, then we define $C_{\cal F}({\cal H}) > C_{\cal F}({\cal H}')$.
Otherwise, we pass to the second triples of 
$C_{\cal F}({\cal H})$ and $C_{\cal F}({\cal H}')$.
Continuing in this way, we can compare the
${\cal F}$-complexities of ${\cal H}$ and ${\cal H}'$.

We now define the {\sl complexity} $C({\cal H})$
of a handle structure ${\cal H}$ to be the ordered
pair $(C_{\cal F}({\cal H}), n({\cal H}))$, where $n({\cal H})$
is the number of 0-handles of ${\cal H}$ containing
a component of ${\cal F}({\cal H})$ with positive
index. We compare the complexity of handle
structures ${\cal H}$ and ${\cal H}'$
by asserting that $C({\cal H}) > C({\cal H}')$
if one of the following holds:
\item{$\bullet$} $C_{\cal F}({\cal H}) > C_{\cal F}({\cal H}')$, or
\item{$\bullet$} $C_{\cal F}({\cal H}) = C_{\cal F}({\cal H}')$ and
$n({\cal H}) < n({\cal H}')$.

\noindent {\bf Lemma 5.3.} {\sl This ordering on complexity of handle 
structures is a well-ordering.}

\noindent {\sl Proof.} We need to show that there cannot exist
an infinite strictly decreasing sequence
$\lbrace C({\cal H}_i): i \in {\Bbb N} \rbrace$.
Suppose that there is such a sequence. 
Then $C_{\cal F}({\cal H}_i) \geq C_{\cal F}({\cal H}_{i+1})$ 
for each $i$.
Suppose first that this inequality is strict for
only finitely many $i$. Then we may pass to a
subsequence in which $C_{\cal F}({\cal H}_i)$ is constant.
Then the number of components of ${\cal F}({\cal H}_i)$ with
positive index is constant. However, since
$C({\cal H}_{i}) > C({\cal H}_{i+1})$ for each $i$,
$n({\cal H}_{i}) < n({\cal H}_{i+1})$ for each $i$. This is impossible.

Therefore, we may suppose that
$C_{\cal F}({\cal H}_{i}) > C_{\cal F}({\cal H}_{i+1})$ for infinitely
many $i$. Pass to this subsequence.
Let $T_i^n$ be the $n^{th}$ largest triple of
$C_{\cal F}({\cal H}_i)$. For each $i$,
there is a natural number $N(i)$, such that
\item{$\bullet$} $T_i^n = T_{i+1}^n$ for $n < N(i)$, and
\item{$\bullet$} $T_i^{N(i)} > T_{i+1}^{N(i)}$.

\noindent Define $M(i) = \min_{j \geq i} N(j)$. Then
$\lbrace M(i) : i \in {\Bbb N} \rbrace$ is a non-decreasing
sequence. For all $i$, $M(i) \leq N(i)$, and for infinitely
many $i$, this is an equality. Consider the
sequence of triples $\lbrace T_i^{M(i)} : i \in {\Bbb N} \rbrace$.
Then $T_i^{M(i)} \geq T_{i+1}^{M(i)} \geq T_{i+1}^{M(i+1)}$.
For the infinitely many $i$ when $M(i) = N(i)$, we
have
$$T_i^{M(i)} = T_i^{N(i)} > T_{i+1}^{N(i)} = T_{i+1}^{M(i)}
\geq T_{i+1}^{M(i+1)}.$$
Thus, the infinite sequence of triples
$\lbrace T_i^{M(i)} : i \in {\Bbb N} \rbrace$
contains an infinite strictly decreasing sequence.
This is impossible, since the ordering on
the triples is a well-ordering.  $\square$

By the above lemma, we can use the complexity of handle structures
as the basis for an inductive argument.
We will start with a sutured manifold $(M, \gamma)$
with a handle structure ${\cal H}$. If $H_2(M, \partial M) \not= 0$,
we will perform a taut decomposition $(M, \gamma) \buildrel S
\over \longrightarrow (M_S, \gamma_S)$.
The manifold $M_S$ will inherit a handle structure ${\cal H}'$.
We will try to ensure that the complexity of ${\cal H}'$ is no
more than that of ${\cal H}$ (and preferably, strictly less
than that of ${\cal H}$). The following lemma asserts that,
to guarantee this, we need only restrict attention to
smaller parts of ${\cal H}$. For example, it shows
that we need only check
$C(H_0 \cap {\cal H}') \leq C(H_0)$ for each 0-handle
$H_0$ of ${\cal H}$.

\noindent {\bf Lemma 5.4.} {\sl Let ${\cal H}$ (respectively
${\cal H}'$) be a handle structure for a sutured manifold $(M, \gamma)$
(respectively $(M', \gamma')$). Suppose that the 0-handles
of ${\cal H}$ (respectively ${\cal H}'$) have been partitioned into
$n$ subsets $A_1, \dots, A_n$ (respectively, $A'_1, \dots, A'_n$).
(For example, each $A_i$ may be some 0-handle $H_0$ of ${\cal H}$, 
and $A'_i$ is $H_0 \cap {\cal H'}$.)
Suppose that for each $i$, $C(A'_i) \leq C(A_i)$.
Then $C({\cal H}') \leq C({\cal H})$. Additionally, if
$C(A'_i) < C(A_i)$ for some $i$, then $C({\cal H}') < C({\cal H})$.}

\noindent {\sl Proof.}
Arrange the triples of $C_{\cal F}({\cal H})$ 
into a non-increasing sequence $\lbrace T_j : j \in {\Bbb N} \rbrace$.
Consider the first integer $j$ for which $T_j > T_{j+1}$.
Then the triples $T_1, \dots, T_j$ are all some
fixed triple $T$.
The partitioning of ${\cal H}^0$ gives a partitioning of
$T_1, \dots, T_j$ into $n$ subsets (some of which may be empty).
Say that $k(i)$ of these lie in $A_i$. Since $C(A'_i) \leq
C(A_i)$, we must have $C_{\cal F}(A'_i) \leq
C_{\cal F}(A_i)$. So,
there are at most $k(i)$ copies of $T$ in $A'_i$,
and there are no larger triples. Hence, in $C_{\cal F}({\cal H}')$, there
are at most $j$ copies of $T$ and no larger triples.
If there are fewer than $j$ copies of $T$ in
$C({\cal H}')$, then $C_{\cal F}({\cal H}') < C_{\cal F}({\cal H})$, and the
lemma is proved. Otherwise, we can remove each copy of
$T$ from $C_{\cal F}({\cal H})$ and $C_{\cal F}({\cal H}')$, without affecting any
ordering. Continuing in this fashion with the next largest triples
of $C({\cal H})$, and so on, we see that $C_{\cal F}({\cal H}') 
\leq C_{\cal F}({\cal H})$.
Also, if we have equality, then we must have
had $C_{\cal F}(A'_i) = C_{\cal F}(A_i)$ for each $i$.
Since $C(A'_i) \leq C(A_i)$, the number of 0-handles in
$A'_i$ containing components of ${\cal F}({\cal H}')$
with positive index is at least the number of 0-handles in
$A_i$ containing components of ${\cal F}({\cal H})$
with positive index. Therefore, $n({\cal H}') \geq n({\cal H})$
and so $C({\cal H}') \leq C({\cal H})$. Also, if we have
equality, then we must have had $C(A'_i) = C(A_i)$ for
each $i$. $\square$

To perform an inductive argument we need to ensure that
the complexity of ${\cal H}'$ is less than
that of ${\cal H}$, where ${\cal H}'$ is the induced handle structure on
$(M_S, \gamma_S)$. However, this is not in general true.
To guarantee this, it is important
that each 0-handle of ${\cal F}$ has positive index,
and to ensure this, we may first need to 
decompose $(M, \gamma)$ along some product
discs and annuli, and then simplify the
handle decomposition of the resulting sutured manifold.
Even then, to ensure that complexity is reduced
by decomposition along $S$, we may need to perform
some modifications to $S$.

We will give these procedures in Sections 7-10.
But first we explain the idea behind the above definition
of complexity.
The surface $S$ is in general $\partial$-compressible in
$M$ and in [4] it was shown that there may exist
infinitely long hierarchies of incompressible
$\partial$-compressible surfaces in a 3-manifold.
Thus, it is vital that we use the fact that $M$
has a sutured manifold structure. This is encoded in
the quantity $C_2(F)$ which was defined to
be the index of a component $F$ of ${\cal F}$.
We therefore study how index behaves under decomposition.

Let ${\cal H}$ (respectively, ${\cal H}'$) be the handle
decomposition of $(M, \gamma)$ (respectively, $(M_S, \gamma_S)$).
Let ${\cal F} = {\cal F}({\cal H})$ and let
${\cal F}' = {\cal F}({\cal H}')$.
Let $V$ be a 0-handle of ${\cal F}$ and let
$V'_1, \dots, V'_k$ be the 0-handles $V \cap {\cal F}'$.
Now, $V'_1, \dots, V'_k$ are obtained from $V$ by
cutting along properly embedded arcs.
The endpoint of each arc either lies in
${\cal R}_\pm(M)$ or in ${\cal F}^1({\cal H})$.
Therefore, an elementary counting argument
shows that
$$I(V) = \sum_{i=1}^k I(V'_i).$$
In particular, if $F$ is a component of
${\cal F}$ and $F' = F \cap {\cal F}'$,
then $I(F) = I(F')$.

Hence, we can ensure that the quantity $C_2$ does not increase, 
as long as we create no discs
of ${\cal F}'$ with negative index. Thus, our
goal is to alter $S$ in order to remove these discs.
But, in general, this does not seem to be possible. An example
is given in Figure 5.2. There, a 0-handle of ${\cal H}$
is decomposed into two 0-handles of ${\cal H}'$.
A negative index disc of ${\cal F}'$ is created, but
note that, nevertheless, the complexity of the handle
structure has decreased. 

It is fairly easy to show that, in general, under mild
assumptions on $S$, neither $C_1$ nor $C_3$ can increase.
Our aim is to show that, if $C_1$ is left unchanged, then in
fact no negative index discs of ${\cal F}'$ are
created, and so $C_2$ is not increased. Furthermore, if
$C_1$, $C_2$ and $C_3$ are all left unchanged, then
$[S, \partial S] = 0 \in H_2(M, \partial M)$.

\vskip 18pt
\centerline {\caps 6. Overview of the proof of the main theorems}
\vskip 6pt

We have now developed enough machinery to outline the
proofs of Theorems 1.4, 1.5 and 1.6. We start with a generalised
triangulation of $M$, and using this, we construct the dual handle
structure ${\cal H}(M)$ (which we sometimes abbreviate to
${\cal H}$). Roughly speaking, the idea is to
decompose $M$ along surfaces until we end with
a solid torus neighbourhood of $K$, plus perhaps some 3-balls. At
each stage, we will be examining a 3-manifold $M'$ embedded
in $M$. This manifold $M'$ will have a handle structure
which respects ${\cal H}(M)$, in the following sense.

\noindent {\bf Definition 6.1.} Let $(M, \gamma)$ be a 
sutured manifold with a handle structure ${\cal H}(M)$. 
Let $(M', \gamma')$ be a sutured manifold
lying in $M$ with a handle structure ${\cal H}(M')$.
Then ${\cal H}(M')$ {\sl respects} ${\cal H}(M)$ if
each of the following conditions holds.
\item{$\bullet$} The 0-handles of $M'$ lie in
the 0-handles of $M$.
\item{$\bullet$} The 1-handles of $M'$ lie in
the 1-handles of $M$ in a vertical fashion
and inherit their product structure.
\item{$\bullet$} The surface ${\cal F}(M')$ lies in
${\cal F}(M)$, with the intersection
${\cal F}^1(M') \cap {\cal F}^1(M)$ lying in
${\cal F}^1(M)$ in a vertical fashion.

Note that, if ${\cal H}(M')$ respects ${\cal H}(M)$,
then automatically the arcs $\gamma' \cap {\cal H}^1(M')$
are vertical in ${\cal H}^1(M)$ and the discs
${\cal H}^2(M') \cap {\cal H}^1(M')$ are vertical
in ${\cal H}^1(M)$. Thus, the only restriction on
the 2-handles of $M'$ is a requirement on their attaching maps.
The remainder of each 2-handle may lie inside $M$ in
a complicated way.

Occasionally, the handle structure of $M'$ will resemble
the handle structure of $M$ in some 0-handle, in the
following sense.

\noindent {\bf Definition 6.2.} Suppose that the handle
structure ${\cal H}(M')$ of $(M', \gamma')$ respects the
handle structure ${\cal H}(M)$ of $(M, \gamma)$.
Let $H_0$ be a 0-handle of $M$. Then ${\cal H}(M') \cap H_0$
is obtained from $H_0$ by a {\sl trivial modification} if
each of the following conditions are satisfied:
\item{(i)} each component of ${\cal F}(M') \cap H_0$ with
positive index lies inside a component of ${\cal F}(M)$
with positive index;
\item{(ii)} the components of ${\cal F}(M') \cap H_0$
with positive index all lie in a single 0-handle $H'_0$ of
${\cal H}(M')$;
\item{(iii)} for each component $F$ of ${\cal F}(M) \cap H_0$
with positive index, one of the following holds:
\itemitem{$\bullet$} there is a component $F'$ of ${\cal F}(M')$
such that $F - F'$ is either empty (in which case, $F'$ is a copy of $F$)
or a collar on components of
$\partial F$ which are disjoint from $\gamma$, the collar
respecting the handle structure of $F$, or
\itemitem{$\bullet$} $F$ is the unique component of ${\cal F}(M) \cap H_0$.
It must be a disc intersecting $\gamma$ in precisely four points. These
four points must lie in precisely two 0-handles of $F$,
two points in each 0-handle. The two curves of $\gamma \cap
H_0$ join points in distinct 0-handles of ${\cal F}(M)$.
Also, ${\rm cl}(\partial H_0 - \partial H'_0)$ is a single disc $D$ which
contains the two curves of $\gamma \cap H_0$ and
intersects ${\cal F}(M)$ in precisely two discs, each
lying in ${\cal F}^0(M)$. Additionally, $H'_0 \cap \gamma'$
must be the two arcs of ${\rm cl}(\partial D - {\cal F}(M))$.

Roughly speaking, a trivial modification leaves components of
${\cal F}(M)$ with positive index relatively unaltered. The final case
in the above definition may seem slightly unnatural, but
it will occur later in our argument (in Section 8), and
it is convenient to define it as trivial.

We make the following
definition. If ${\cal H}(M)$ is a handle structure for $M$,
we define the {\sl important} 0-handles ${\cal IH}^0(M)$ to be the
0-handles $H_0$ with $H_0 \cap {\cal F}$ containing
at least one component with positive index.
In Sections 7 and 8, we will prove the following result,
which gives a method of modifying a handle decomposition
so that, afterwards, each 0-handle of ${\cal F}$ has positive index.
Recall from Section 5 that the index of a component of ${\cal F}$
is equal to the sum of the indices of its 0-handles.
So this implies that each component of ${\cal F}$ has positive
index. Therefore, each 0-handle of ${\cal H}$ is
either important or disjoint from the 1-handles and 2-handles.

\vfill\eject

\noindent {\bf Proposition 6.3.} {\sl Let ${\cal H}(M)$
be a handle structure of a taut sutured manifold $(M, \gamma)$.
Suppose that each component of $M$ has non-empty
boundary, and that no component of $M$ is a solid torus.
Suppose also that no component of $M$ is a Seifert fibre space disjoint
from $\gamma$, with base space a disc and having two exceptional fibres.
Then there is a (possibly empty) sequence of
taut decompositions
$$(M, \gamma) \buildrel P_1 \over \longrightarrow
\dots \buildrel P_m \over \longrightarrow (M', \gamma'),$$
where each $P_i$ is either a product disc or an incompressible
annulus disjoint from the sutures. There is a handle structure 
${\cal H}(M')$ of $(M', \gamma')$ and an embedding
of $M'$ in $M$ isotopic to the embedding arising from
the sutured manifold decomposition, with the following properties.
\item{(i)} ${\cal H}(M')$ respects ${\cal H}(M)$.
\item{(ii)} For each 0-handle $H_0$ of ${\cal H}(M)$,
the complexity of $H_0 \cap {\cal H}(M')$ is no more than that of $H_0$.
\item{(iii)} For each 0-handle $H_0$ of ${\cal H}(M)$, the
intersections 
$$\eqalign{& H_0 \cap {\cal IH}^0(M')\cr
& H_0 \cap {\cal IH}^0(M') \cap {\cal F}(M')\cr
& H_0 \cap {\cal IH}^0(M') \cap \gamma' \cr}$$
are each one of a finite number of possibilities (up to trivial
modifications), which depend only on ${\cal F}(M) \cap H_0$ and
$H_0 \cap \gamma$, and are otherwise independent
of $M$ and $M'$.
\item{(iv)} If $H_0$ is a 0-handle of ${\cal H}(M)$
and the complexity of $H_0 \cap {\cal H}(M')$ is
equal to that of $H_0$, then $H_0 \cap {\cal H}(M')$
is obtained from $H_0$ by a trivial modification.
\item{(v)} Each 0-handle of ${\cal F}(M')$ has positive
index.
\item{(vi)} For each 0-handle $H'_0$ of ${\cal H}(M')$,
$H'_0 \cap ({\cal F}(M') \cup \gamma')$ is connected.

}

Once we have such a handle structure, we then perform a
sutured manifold decomposition.

\noindent {\bf Proposition 6.4.} {\sl Let ${\cal H}(M)$ be
a handle structure of a taut sutured manifold $(M, \gamma)$.
Suppose that each 0-handle of ${\cal F}(M)$ has positive index,
and that, for each 0-handle $H_0$ of ${\cal H}(M)$,
$H_0 \cap ({\cal F}(M) \cup \gamma)$ is connected.
Let $(M, \gamma) \buildrel S \over \longrightarrow (M_S, \gamma_S)$
be a taut sutured manifold decomposition,
where $\partial S$ has essential intersection with ${\cal R}_\pm(M)$ 
and $[S, \partial S] \not= 0
\in H_2(M, \partial M)$. Then there is a surface
$S'$ properly embedded in $(M, \gamma)$ and a 
commutative diagram of sutured manifold decompositions
and pull-backs
$$
\matrix{
(M, \gamma) 
&\smash{\mathop{\longrightarrow}\limits^{S}}
&(M_S, \gamma_S) &\smash{\mathop{=}}
& (\hat M_1, \hat \gamma_1) 
&\smash{\mathop{\longleftarrow}\limits^{P_1}}
& \dots
&\smash{\mathop{\longleftarrow}\limits^{P_{r-1}}}
& (\hat M_r, \hat \gamma_r)\cr
\Big\downarrow \rlap{$\vcenter{\hbox{${\scriptstyle S'}$}}$}
&&&&&&&&
\Vert \cr
(M', \gamma') &\smash{\mathop{=}}
& (\hat M_m, \hat \gamma_m) 
&\smash{\mathop{\longleftarrow}\limits^{P_{m-1}}}
& (\hat M_{m-1}, \hat \gamma_{m-1})
&\smash{\mathop{\longleftarrow}\limits^{P_{m-2}}}
& \dots
&\smash{\mathop{\longleftarrow}\limits^{P_{r}}}
& (\hat M_r, \hat \gamma_r).\cr}$$
Each $P_i$ is either a product disc, an incompressible product annulus
or (for $i < r$) a surface parallel to a subsurface $F_i$ of 
${\cal R}_\pm(\hat M_{i+1})$, with the orientations of
$P_i$ and $F_i$ disagreeing near $\partial P_i$.
The induced handle structure ${\cal H}(M')$ on $(M', \gamma')$
satisfies properties (i), (ii), (iii) and (iv) of Proposition 6.3
and also the following:
\item{(v)} For some 0-handle $H_0$ of ${\cal H}(M)$, 
$C({\cal H}(M') \cap H_0) < C(H_0)$.

}

\noindent {\sl Proof of Theorems 1.4, 1.5 and 1.6
using Propositions 6.3 and 6.4.}
Let ${\cal H}$ be the dual handle structure for $M$,
arising from the generalised triangulation of $M$. 
We give $M$ the trivial sutured manifold structure
with ${\cal R}_- = \partial M$ and ${\cal R}_+ = \emptyset$.
If this is not taut, then $M$ is either reducible or
a solid torus. Hence, by Theorem 5.1 of [10], there are no
exceptional or norm-exceptional surgery curves in $M$ satisfying
the hypotheses of Theorem 1.4. Hence we may
assume that $(M, \emptyset)$ is taut.

We will construct a sequence of taut sutured manifolds
$(M_i, \gamma_i)$ where $1 \leq i \leq n$.
The first sutured manifold $(M_1, \gamma_1)$ will be $(M, \emptyset)$.
Each sutured manifold $(M_i, \gamma_i)$ will have a handle structure
${\cal H}_i$, and there will be an embedding of $M_{i}$ in
$M_{i-1}$ having the following properties (some of which are
only relevant for $i>1$):
\item{(i)} ${\cal H}_i$ respects ${\cal H}_{i-1}$.
\item{(ii)} For each 0-handle $H_0$ of ${\cal H}_{i-1}$,
the complexity of $H_0 \cap {\cal H}_{i}$ is no more
than that of $H_0$.
\item{(iii)} For each 0-handle $H_0$ of ${\cal H}_{i-1}$, the
intersections 
$$\eqalign{& H_0 \cap {\cal IH}^0_{i}\cr
& H_0 \cap {\cal IH}^0_{i} \cap {\cal F}({\cal H}_{i})\cr
& H_0 \cap {\cal IH}^0_{i} \cap \gamma_{i} \cr}$$
are each one of a finite number of possibilities (up to
trivial modification), which depend
only on ${\cal F}({\cal H}_{i-1}) \cap H_0$
and $H_0 \cap \gamma_{i-1}$.
\item{(iv)} For a 0-handle $H_0$ of ${\cal H}_{i-1}$, if
the complexity of $H_0 \cap {\cal H}_{i}$ is equal to that
of $H_0$, then $H_0 \cap {\cal H}_i$ is obtained from
$H_0$ by a trivial modification.
\item{(v)} For some 0-handle $H_0$ of ${\cal H}_{i-1}$ ($i > 2$), we
have $C(H_0 \cap {\cal H}_{i}) < C(H_0)$.
\item{(vi)} For $1 < i < n$, each 0-handle of ${\cal F}({\cal H}_{i})$
has positive index.
\item{(vii)} $K$ lies in $M_{i}$.
\item{(viii)} For $1 \leq i < n$,
$H_2(M_i - {\rm int}({\cal N}(K)), \partial M_i) \not=0$.
\item{(ix)} If $M_i(\sigma)$ is the manifold
obtained from $M_i$ by Dehn surgery along $K$ with slope
$\sigma$, then at least one of the following is
true:
\itemitem{$\bullet$} $M_i = M$,
\itemitem{$\bullet$} $(M_i(\sigma), \gamma_i)$ is not taut, or
\itemitem{$\bullet$} the core of the surgery solid torus
has finite order in $\pi_1(M_i(\sigma))$.

\noindent The final manifold $M_n$ of the sequence is a solid torus neighbourhood of $K$,
plus perhaps some 3-balls.
The sequence is constructed using Propositions 6.3 and 6.4 in an
alternating fashion.

We now show how to continue this sequence.
We would like to let $(M, \emptyset) = (M_2, \gamma_2)$, but 
(vi) above need not be satisfied in this case.
Note, however, that each 0-handle of ${\cal F}(M)$ does
indeed have positive index in either of the following
cases:
\item{$\bullet$} $M$ is closed (and so we have a genuine triangulation), or
\item{$\bullet$} $\partial M \not= \emptyset$ and we have an ideal
triangulation.

\noindent In the case where $K$ and $\sigma$ are norm-exceptional,
we would like to ensure that one of the above is true.
In Theorem 1.6, we explicitly make this assumption.
In Theorem 1.4, we alter the given generalised
triangulation of $M$ so that it is either genuine or ideal.
This can be done algorithmically. Hence, in the case where $K$ and
$\sigma$ are norm-exceptional, we let $(M, \emptyset)
= (M_2, \gamma_2)$.

Suppose now that $K$ and $\sigma$ are exceptional and
that some 0-handle of ${\cal F}$ has non-positive
index. Then we use Proposition 6.3 to decompose
$(M, \emptyset)$ along product discs and 
incompressible annuli disjoint from the sutures, resulting in
a taut sutured manifold $(M', \gamma')$ satisfying
(i) - (vi) of 6.3. Since $\Delta(\sigma, \mu) > 1$,
Theorem 3.2 gives that we may ambient
isotope $K$ off each decomposing surface, and
hence $K$ lies in $M'$.

To apply Proposition 6.3, we need to check that
$M$ is not a Seifert fibre space,
with base space a disc, having two exceptional fibres,
and having $\gamma = \emptyset$. We will suppose
it is, and then achieve a contradiction.
Let $\alpha$ be an arc properly embedded
in the base space, separating the two exceptional
points. Then $\alpha$ lifts to an annulus $A$ in $M$.
Using Theorem 3.2, we may ambient isotope $K$ off $A$.
Let $\hat M$ be the solid torus which is the closure of the component of 
$M - A$ containing $K$. Then, $\partial \hat M$
is an incompressible torus in $M - {\rm int}({\cal N}(K))$.
Since $M - {\rm int}({\cal N}(K))$ is atoroidal, 
we deduce that $\partial \hat M$ must be parallel to
$\partial {\cal N}(K)$ and so $K$ is isotopic
to the exceptional fibre lying in $\hat M$. But then
$H_2(M - {\rm int}({\cal N}(K)), \partial M)$ is
trivial, contrary to assumption. Hence, we may
apply Proposition 6.3, and then we
let $(M_2, \gamma_2) = (M', \gamma')$.

We now verify (ix). We only applied Proposition 6.3
in the case where $K$ and $\sigma$ are exceptional. Hence,
either $(M_K(\sigma), \emptyset)$ is not taut or
the core of the surgery solid torus has finite order in 
$\pi_1(M_K(\sigma))$. We shall show that
either $(M'_K(\sigma), \gamma')$ is not taut or
the core of the surgery solid torus has finite
order in $\pi_1(M'_K(\sigma))$. Suppose that $(M'_K(\sigma), \gamma')$ is taut.
Then, by 3.1, the sequence of decompositions
$$(M_K(\sigma), \gamma) \buildrel P_1 \over \longrightarrow
\dots \buildrel P_m \over \longrightarrow (M'_K(\sigma), \gamma'),$$
is taut. In particular, $(M_K(\sigma), \gamma)$ is taut,
and therefore (iii) of 1.3 holds. Also, since
each $P_i$ has essential intersection with ${\cal R}_\pm$, it is therefore
incompressible. Therefore,
the map $\pi_1(M'_K(\sigma)) \rightarrow
\pi_1(M_K(\sigma))$ induced by inclusion is an injection.
Therefore, the core of the surgery solid torus has finite
order in $\pi_1(M'_K(\sigma))$.

Thus, we have now constructed $(M_2, \gamma_2)$ and have
verified that it has the correct properties.
Suppose that we have constructed a sequence 
as far as $(M_i, \gamma_i)$, satisfying (i) - (ix) above.
If $H_2(M_i - {\rm int}({\cal N}(K)), \partial M_i)$ is trivial,
then we stop. If this homology group is non-trivial, then
(see Section 3 or Theorem 2.6 of [10]) we may find a taut decomposition
$$(M_i - {\rm int}({\cal N}(K)), \gamma_i) \buildrel S \over \longrightarrow
(M'_i - {\rm int}({\cal N}(K)), \gamma'_i),$$
such that
\item{$\bullet$} $S$ is disjoint
from $\partial {\cal N}(K)$, 
\item{$\bullet$} no curve of $\partial S$ bounds a disc in ${\cal R}_\pm(M_i)$,
\item{$\bullet$} no component $X$ of $M'_i$ has $\partial X \subset
{\cal R}_-(M'_i)$ or $\partial X \subset {\cal R}_+(M'_i)$, and
\item{$\bullet$} $[S, \partial S] \not=0 \in H_2(M_i
-{\rm int}({\cal N}(K)), \partial M_i)$. 

\noindent This
implies that $[S, \partial S] \not=0 \in H_2(M_i, \partial M_i)$.
In the case where $K$ and $\sigma$ are norm-exceptional
and $M_i = M$, we insist that $[S, \partial S] = z \in 
H_2(M - {\rm int}({\cal N}(K)), \partial M)$, where $z$
is the homology class in Definition 1.3.

Let $M'_i(\sigma)$ be the result of $M'_i$ after Dehn surgery along
$K$ with slope $\sigma$. Since $(M_i(\sigma), \gamma_i)$
is not taut or the core of the surgery solid torus has
finite order in $\pi_1(M_i(\sigma))$, the
argument above gives that $(M'_i(\sigma), \gamma_i)$ is not
taut or the core of the surgery solid torus has finite order
in $\pi_1(M'_i(\sigma))$. 
Using the argument in Section 3 (see also Theorem 1.8 of [1]), 
we deduce that the decomposition
$(M_i, \gamma_i) \buildrel S \over \longrightarrow
(M'_i, \gamma'_i)$ is taut.
Since each 0-handle of ${\cal F}({\cal H}_i)$ has positive index,
we may apply Proposition 6.4 to $(M_i, \gamma_i)$,
and so end with a sutured manifold $(M''_i, \gamma''_i)$
with a handle-decomposition ${\cal H}''_i$, satisfying
(i) - (iv) of 6.3 and (v) of 6.4. Again using 3.2, we may isotope
$K$ off each product disc and incompressible product
annulus, and so we may assume $K$ lies in $M''_i$.
Also, using the commutative diagram in Proposition 6.4,
we deduce that $(M''_i(\sigma), \gamma''_i)$ is not
taut or the core of the surgery solid torus has finite order
in $\pi_1(M''_i(\sigma))$. 

If any component of $M_i''$ is a solid torus disjoint
from $K$, we decompose it along a meridian disc.
If the component $X$ of $M''_i$ containing $K$
is a solid torus, then the atoroidality of 
$M - {\rm int}({\cal N}(K))$ implies that
$K$ is the core of $X$. In this case, we set
$(M_n, \gamma_n) = (X, \gamma''_i \cap X)$,
together with some 3-balls obtained by decomposing
$M''_i - X$. 

We may therefore assume that no component of $M''_i$ is a solid torus,
and so we can apply 6.3 to $(M''_i, \gamma''_i)$ to obtain
a sutured manifold $(M_{i+1}, \gamma_{i+1})$ satisfying
(i) - (ix) above. Note that each component of
${\cal F}({\cal H}_{i+1})$ has positive index,
and therefore, the only 0-handles of ${\cal H}_{i+1}$
not lying in ${\cal IH}^0_{i+1}$ are handles disjoint
from ${\cal F}({\cal H}_{i+1})$.

By (ii), (v) and Lemma 5.4, 
the complexity of ${\cal H}_{i+1}$ is strictly
less than that of ${\cal H}_i$. Hence, eventually, the
sequence terminates with a sutured manifold
$(M_n, \gamma_n)$ such that 
$H_2(M_n - {\rm int}({\cal N}(K)), \partial M_n) =  0$.
Then $M_n$ is a solid torus neighbourhood of $K$, plus perhaps 
some 3-balls. By (ix), the sutures $\gamma_n \cap {\cal N}(K)$
are parallel to $\sigma$.

Note that, for each 0-handle $H_0$
of ${\cal H}(M)$, there are only finitely many possibilities for
${\cal F}(M) \cap H_0$. 
Thus, by induction on complexity using (ii), (iii) and (iv) above, there is in
each 0-handle $H_0$ of ${\cal H}(M)$, only a finite number of possibilities for
$$\eqalign{& H_0 \cap {\cal IH}^0_n\cr
& H_0 \cap {\cal IH}^0_n \cap {\cal F}({\cal H}_n)\cr
& H_0 \cap {\cal IH}^0_n \cap \gamma_n. \cr}$$
Each possibility for $H_0 \cap {\cal IH}^0_n \cap \gamma_n$
gives a tangle in the associated 3-simplex of $M$.
These tangles join to form $\gamma_n$ (with possibly some
unknotted curves removed). Some component of $\gamma_n$
is the $\sigma$-cable of $K$, and hence
the tangles required for Theorems 
1.5 and 1.6 are constructed by taking all possible
subtangles of $H_0 \cap {\cal IH}^0_n \cap \gamma_n$.

Each possibility for 
$H_0 \cap {\cal IH}^0_n$ and
$H_0 \cap {\cal IH}^0_n \cap {\cal F}^0({\cal H}_n)$
gives a graph $G$ in the associated 3-simplex of $M$. When
the collection of these graphs (one in each 3-simplex of 
$M$) are joined, they form the 0-handles and
1-handles of $M_n$. The 2-handles of $M_n$ are attached
along the annuli 
$({\cal H}^0_n \cup {\cal H}^1_n) \cap {\cal H}^2_n$,
which are determined by
${\cal H}^0_n \cap {\cal F}^1({\cal H}_n)$.
Thus, we readily see that the algorithm given in
Section 2 constructs all possibilities for $K$ and $\sigma$.
Hence, Theorem 1.4 is established. $\square$

We end this section with the following.

\noindent {\sl Proof of Theorem 1.1.}
This is an almost immediate corollary of Theorem 1.4,
but there is one complication. If we set
$M = S^3 - {\rm int}({\cal N}(L))$, then
it is not obvious that $M - {\rm int}({\cal N}(K))$
is atoroidal. To establish this, we will use
a modified form of the argument in Proposition 2.3
of [6].

Suppose therefore that $T$ is an incompressible torus in 
$M - {\rm int}({\cal N}(K))$ which is not parallel to
$\partial {\cal N}(K)$ or $\partial {\cal N}(L)$.
Since we are assuming that $K$ is not a
non-trivial satellite knot, then $T$ must be compressible
in $M$ or be parallel to $\partial {\cal N}(L)$.
In the latter case, $K$ lies in the collar between
$T$ and $\partial {\cal N}(L)$, and then it is easy
to see that $L'$ is a winding number one satellite of $L$. 
In particular, ${\rm genus}(L') \geq {\rm genus}(L)$, which is 
contrary to hypothesis. Therefore, $T$ must be compressible in $M$. 
Since $T$ does not lie in a 3-ball in $M$,
it must bound a solid torus $V$ in $M$ which
contains $K$. Let $V'$ be the manifold obtained
from $V$ by $1/q$ Dehn surgery along $K$.

\noindent {\sl Case 1.} $K$ has winding number zero in
$V$. 

Then consider a minimal genus Seifert surface $S$
for the knot $L'$. We may assume that it intersects 
$\partial V'$ in a collection of simple closed curves,
each of which is homologically trivial in $V'$.
Since $K$ has winding number zero in $V$,
these curves are also homologically trivial in
$V$. Hence, we may fill them in with meridian discs in
$V$. This gives a Seifert surface for $L$ with
genus at most that of $S$, which is a contradiction.

\noindent {\sl Case 2.} $K$ has non-zero winding number in $V$.

Since $K$ and $L$ have zero linking number, so do $L$ and the
core of $V$. Therefore, there exists a Seifert surface $S'$
for $L'$ which is disjoint from $V'$.

\noindent {\sl Case 2A.} $\partial V'$ is incompressible
in $V'$. 

Then, by Lemma A.16 of [7], there is a minimal genus
Seifert surface for $L'$ which is disjoint from
$V'$. This gives a Seifert surface for $L$,
and again we reach the contradiction that
${\rm genus}(L) \leq {\rm genus}(L')$. 

\noindent {\sl Case 2B.} $\partial V'$ is compressible in $V'$.

This implies that $V'$ is a solid torus, since it cannot
be reducible, as it lies in $S^3$. Now, $1/q$ surgery along a knot $K$
in the solid torus $V$ never yields another
solid torus if $\vert q \vert > 1$, unless $K$
is a core of $V$ or $K$ lies in a 3-ball in $V$ [9]. 
If $K$ lies in a 3-ball in $V$, then $1/q$ surgery along $K$ does not
alter $L$, which is a contradiction.
If $K$ is a core of $V$, then $T$ is parallel to
$\partial {\cal N}(K)$, contradicting the assumption that
it is essential. 

This proves then that $M - {\rm int}({\cal N}(K))$ is
atoroidal. Theorem 1.1 now follows directly from Theorem 1.4.
$\square$

\vfill\eject
\centerline {\caps 7. Simplifying handle structures}
\vskip 6pt

In the next two sections, we will give a proof of Proposition 6.3.
In particular, we will assume that each component of $M$ has non-empty
boundary, and that no component of $M$ is a solid torus or a Seifert fibre space as
in 6.3. We start by giving various elementary procedures for
simplifying a handle structure ${\cal H}$ of the taut sutured manifold
$(M, \gamma)$. Our aim is to end with a handle structure
in which each 0-handle of ${\cal F}$ has positive index.
Each procedure will satisfy (i) - (iv) of 6.3.
It will not increase the complexity of the handle
structure, but it need not decrease it. To ensure that
these procedures eventually terminate, we therefore
introduce the following definition.

\noindent {\bf Definition 7.1.} Let ${\cal H}$ be a handle
structure for a sutured manifold $(M, \gamma)$.
Define the {\sl extended ${\cal F}$-complexity} for ${\cal H}$ to
be the set of triples
$$C_{\cal F}^+({\cal H}) = \lbrace (C_1(F), C_2(F), C_3(F)) :
\hbox{$F$ a component of ${\cal F}$} \rbrace,$$
where repetitions are retained. Here, $C_1$, $C_2$ and $C_3$
are the integers defined in Section 5. We also define the {\sl extended
complexity} $C^+({\cal H})$ to be the ordered pair
$(C_{\cal F}^+({\cal H}), n({\cal H}))$.

The difference between the extended ${\cal F}$-complexity and
the ${\cal F}$-complexity of a handle structure is that extended
${\cal F}$-complexity also takes into account components of
${\cal F}$ with non-positive index.

We order the extended ${\cal F}$-complexities and
extended complexities as we do the ${\cal F}$-complexities and
complexities of handle structures (see Section 5). As in Lemma 5.3,
this is a well-ordering. The procedures we give in the
next two sections will all reduce the extended complexity
of the handle structure, and so are guaranteed to terminate.

The following lemma will be useful in our verification
that (i) - (iv) of 6.3 holds and that extended complexity
is reduced.

\noindent {\bf Lemma 7.2.} {\sl Let $(M, \gamma)$ be a sutured
manifold with a handle structure ${\cal H}$. Let $(M', \gamma')$
be embedded in $M$, with a handle structure ${\cal H}'$
which respects ${\cal H}$. Suppose that, for each
component $F$ of ${\cal F}$, either $C_{\cal F}^+(F \cap {\cal F}')
< C_{\cal F}^+(F)$ or $F \subset {\cal F}'$. Suppose also
that the former of the above two possibilities holds for
at least one component $F$ of ${\cal F}$.
Suppose also that, if $I(F) \leq 0$, then
the index of each component of $F \cap {\cal F}'$ is non-positive.
Then, $C_{\cal F}^+({\cal H}') < C_{\cal F}^+({\cal H})$, and
(i), (ii) and (iv) of 6.3 are verified.}

\noindent {\sl Proof.} A version of the argument in Lemma
5.4 gives that $C_{\cal F}^+({\cal H}') < C_{\cal F}^+({\cal H})$. 
We now check (ii) and (iv) of 6.3. (Note that (i) of 6.3 is
part of the hypothesis of the lemma.)
For each component $F$ of ${\cal F}$, we have one
of the following possibilities:
\item{(i)} $F \subset {\cal F}'$ and so $F \cap {\cal F}'$ is 
a copy of $F$, or
\item{(ii)} $I(F) > 0$ and 
$C_{\cal F}(F \cap {\cal F}')  \leq
C_{\cal F}^+(F \cap {\cal F}') < C_{\cal F}^+(F) = C_{\cal F}(F)$, or
\item{(iii)} $I(F) \leq 0$ and $C_{\cal F}(F \cap {\cal F}') =
C_{\cal F}(F)$.

In (iii), we are using that if $I(F) \leq 0$, then
the index of each component of $F \cap {\cal F}'$ is non-positive,
and so does not contribute to ${\cal F}$-complexity.
Therefore, for any 0-handle $H_0$ of ${\cal H}$, 
$C_{\cal F}({\cal H}' \cap H_0) \leq C_{\cal F}(H_0)$. 
Also, if we have equality, then (ii) above cannot occur
for any component $F$ of ${\cal F} \cap H_0$, which
implies that components $F$ of ${\cal F} \cap H_0$ with 
positive index remain unchanged and hence that
$n({\cal H}' \cap H_0) \geq n(H_0)$. So,
$C({\cal H}' \cap H_0) \leq C(H_0)$. This verifies (ii) of
6.3. Also, if $C({\cal H}' \cap H_0) = C(H_0)$, then
${\cal H}' \cap H_0$ is obtained from $H_0$ by a trivial
modification, verifying (iv) of 6.3. $\square$

Before we describe the procedures in detail, we mention
that many of them simply remove some handles of ${\cal H}$.
The following lemma will therefore be useful.

\noindent {\bf Lemma 7.3.} {\sl Let ${\cal H}'$ be a collection
of handles of ${\cal H}$ forming a 3-manifold $M'$ embedded
in $M$, with ${\cal H} - {\cal H}'$ containing at least one $i$-handle
for some $i \leq 2$. Suppose that ${\cal H}'$ is a handle structure,
that each handle of ${\cal H} - {\cal H}'$ is disjoint from
$\gamma$ and that $(M', \gamma)$ is a sutured manifold structure.
Then, (i) - (iv) of 6.3 are satisfied, and
extended ${\cal F}$-complexity is reduced.}

\noindent {\sl Proof.} It follows straight from the definition
that ${\cal H}'$ respects ${\cal H}$. 
Let us now check that the hypotheses of 7.2 hold. Let $F$ be some
component of ${\cal F}$ and let $F' = F \cap {\cal F}'$.
The 1-handles of $F$ are either removed or divided up
amongst the components of $F'$. In particular, each component
$X$ of $F'$ has $C_1(X) \leq C_1(F)$. If this inequality is
an equality for some $X$, then in fact $F' = F$.
Hence, either $C_{\cal F}^+(F') < C_{\cal F}^+(F)$ or $F \subset
{\cal F}'$. Also, the former case holds for some component
$F$ of ${\cal F}$.

We now check that if $I(F) \leq 0$, then
each component of $F'$ has non-positive index.
But $F'$ is obtained from $F$ by removing some 1-handles
(or equivalently, cutting $F$ along properly embedded arcs),
then removing some 0-handles disjoint from $\gamma$. Thus,
each component $X$ of $F'$ has $\chi(X) \geq \chi(F)$
and $\gamma \cap X \leq \gamma \cap F$. Hence,
$X$ does not have positive index if $F$ does not have
positive index. 

Thus, by Lemma 7.2, (i), (ii) and (iv) of 6.3 hold
and extended ${\cal F}$-complexity is reduced. Also, (iii) of 6.3
is obvious. $\square$

\noindent {\bf Procedure 1.} Slicing a 0-handle along a disc.

Suppose that there is a disc $D$ properly embedded in
some 0-handle $H_0$ with $D \cap \gamma = D \cap {\cal F} =
\emptyset$, and which separates ${\cal F} \cap H_0$.
Then, $\partial D$ either lies in $\partial {\cal H}^3$ or in
${\cal R}_\pm$. In the former case, $\partial D$ bounds a
disc $D'$ in $\partial {\cal H}^3$, since each component
of $\partial {\cal H}^3$ is a sphere. In the case where
$\partial D$ lies in ${\cal R}_\pm$,
the incompressibility of ${\cal R}_\pm$ implies that
$\partial D$ bounds a disc $D'$ in ${\cal R}_\pm$.
The irreducibility of $M$ implies that $D \cup D'$ bounds a
ball $B$ in $M$. Procedure 1 is the removal of all handles
of ${\cal H}$ which intersect ${\rm int}(B)$, other than
$H_0$. If $D' \subset \partial {\cal H}^3$, we extend
${\cal H}^3$ over $B$. By Lemma 7.3, (i) - (iv) of 6.3 hold,
and extended ${\cal F}$-complexity is reduced. Thus, using this
procedure, we eventually obtain a handle structure ${\cal H}(M')$ on
the resulting sutured manifold $(M', \gamma')$,
with $H'_0 \cap ({\cal F}(M') \cup \gamma')$
connected, for each 0-handle $H'_0$ of ${\cal H}(M')$.
Thus is (vi) of 6.3.

\noindent {\bf Procedure 2.} Collapsing a 2-handle and a 1-handle disjoint
from $\gamma$.

Suppose now that $H_1$ is a 1-handle of $M$ which is disjoint
from $\gamma$ and which intersects ${\cal H}^2$ in a single
disc. Then this disc is contained in a single 2-handle $H_2$.
We may remove $H_1$ and $H_2$, without changing the
homeomorphism type of $(M, \gamma)$. 
Lemma 7.3 gives that (i) - (iv) of 6.3 are satisfied, 
and that extended ${\cal F}$-complexity is reduced.

\noindent {\bf Procedure 3.} Collapsing a 2-handle and a 1-handle
containing an arc of $\gamma$.

Let $H_1$ be a 1-handle of $M$ which intersects
$\gamma$ in a single arc, and which intersects ${\cal H}^2$
in a single disc, lying in a 2-handle $H_2$.
Procedure 3 is the collapsing of $H_1$ and $H_2$.
This moves $\gamma \cap H_1$ onto an arc running along
$\partial ({\cal H}^0 \cup {\cal H}^1 - H_1)$.
Let $(M', \gamma')$ be the new sutured manifold,
with handle structure ${\cal H}'$.

This procedure has the following effect on ${\cal F}$: removing
$H_2 \cap \partial {\cal H}^0$ (which is a collection
of 1-handles of ${\cal F}$), removing $H_1 \cap \partial {\cal H}^0$
(which is precisely two 0-handles of ${\cal F}$) and then replacing
each handle of ${\cal F}$ which we have removed with 
a sub-arc of $\gamma$. Thus, if $F$ is a component of
${\cal F}$, and $F' = F \cap {\cal F}'$, then
each component of $X$ of $F'$ has $C_1(X) \leq C_1(F)$,
and if we have equality for some component $X$, then
in fact $F$ is unchanged by the procedure. This verifies
one of the hypotheses of Lemma 7.2.

We now check that if $F$ has non-positive
index, then each component $X$ of $F'$ has non-positive
index. Suppose therefore $F$ has non-positive index
and that $F$ is changed by the procedure.
It is simple to show that $I(F)$ is the sum of the
indices of $X$, as $X$ ranges over all components of $F'$.
Therefore, the only way that a component $X$ of
$F'$ can have positive index is if another component
of $F'$ has negative index. However, 
since $F$ is changed by the procedure, then
each component of $F'$ touches $\gamma'$, and hence $I(X) \geq 0$
for all components $X$ of $F'$. Lemma 7.2 now gives
us that extended ${\cal F}$-complexity decreases and
that (i), (ii) and (iv) of 6.3 hold. It is straightforward
to verify (iii) of 6.3. 

\noindent {\bf Procedure 4.} Decomposing along a product disc,
then sliding $\gamma$.

Suppose that $H_1$ is a 1-handle of ${\cal H}$ which
is disjoint from ${\cal H}^2$ and which has
$\vert H_1 \cap \gamma \vert = 2$. Let $D$ be one
of the two discs of $H_1 \cap {\cal H}^0$, lying in
some 0-handle $H_0$. Push $D$ a little
into $H_0$. Then $D$ is a product disc, which we
decompose along. This decomposition creates a
new handle decomposition which respects ${\cal H}$,
and leaves both the complexity and extended complexity
of ${\cal H}$ unchanged. But now the two arcs of
$H_1 \cap \gamma$ are joined by an arc of
${\cal H}^0 \cap \gamma$. We may therefore perform an ambient
isotopy which slides $\gamma$ off $H_1$. 
Then, using Procedure 1, we may remove $H_1$.
Again, (i) - (iv) of 6.3 are satisfied and
extended ${\cal F}$-complexity is reduced.

\noindent {\bf Procedure 5.} Collapsing a 3-ball disjoint from $\gamma$.

Suppose that a component of $M$ is a 3-ball disjoint from $\gamma$, comprised
of two 0-handles joined by a 1-handle. Then we may collapse
the 1-handle and one of the 0-handles. This reduces
the extended ${\cal F}$-complexity and (i) - (iv) of 6.3
are satisfied.

\vfill\eject
\noindent {\bf Procedure 6.} Collapsing a 2-handle and a 3-handle.

Let $H_2 = D^1 \times D^2$ be a 2-handle, with one component
of $\partial D^1 \times D^2$ in $\partial M$, and
the other component touching a 3-handle $H_3$. Then we may remove
$H_2$ and $H_3$ without changing the homeomorphism type
of $M$. By Lemma 7.3, this procedure reduces extended ${\cal F}$-complexity and
satisfies (i) - (iv) of 6.3. Note that we are assuming
in 6.3 that each component of $M$ has non-empty boundary. Hence, if
${\cal H}^3$ is non-empty, we may always apply
this procedure somewhere. In this way, we remove all 3-handles
from $M$.

The above six procedures are not enough to ensure
that each 0-handle of ${\cal F}$ has positive index.
To deal with components of ${\cal F}$ which are
annuli disjoint from $\gamma$, we must clump collections
of handles into groups, known as amalgams, which
are defined as follows.

\noindent {\bf Definition 7.4.} An {\sl amalgam} ${\cal A}$
is a connected collection of handles with the following properties:
\item{(i)} ${\cal A}$ is disjoint from $\gamma$,
\item{(ii)} ${\cal A}$ is an $I$-bundle over a connected surface $G$,
\item{(iii)} the $(I - \partial I)$-bundle over $\partial G$ 
is disjoint from $\partial M \cup \partial {\cal H}^3$,
\item{(iv)} the handles of ${\cal A}$ touching the
$(I - \partial I)$-bundle over $\partial G$ are 1-handles
and 2-handles,
\item{(v)} no 2-handle or 3-handle of ${\cal H} - {\cal A}$ touches
${\cal A}$,
\item{(vi)} the $\partial I$-bundle over $G$ lies in
${\cal R}_\pm$, and
\item{(vii)} ${\rm cl}({\cal H} - {\cal A})$ inherits a handle
structure from ${\cal H}$.

\noindent An amalgam is {\sl trivial} if it is a single 2-handle; otherwise
it is {\sl non-trivial}. 

An amalgam ${\cal A}$ behaves in many ways just like a 2-handle.
For example, it is attached onto the 0-handles and 1-handles of 
${\cal H} - {\cal A}$ in a fashion that is very similar
to the attachment of a 2-handle. 

The main example of a non-trivial amalgam ${\cal A}$ is a connected collection of
2-handles and 1-handles disjoint from $\gamma$ and ${\cal H}^3$, such that each
1-handle of ${\cal A}$ intersects ${\cal H}^2$ in precisely
two discs, and these discs lie in 2-handles of ${\cal A}$.
For then the co-core $D^2$
of each 1-handle $H_1 = D^2 \times D^1$ in ${\cal A}$ 
has a product structure as
$I \times I$, in which $H_1 \cap {\cal H}^2 = \partial I \times
I \times D^1$. The product structures on the 1-handles
combine with the product structures on the 2-handles
to form an $I$-bundle structure on ${\cal A}$ with the
required properties. An example is
given in Figure 7.5.

\vskip 18pt
\centerline{\psfig{figure=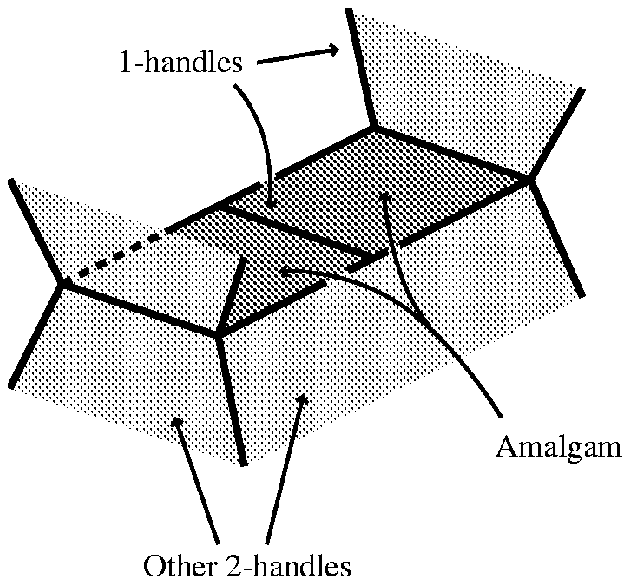}}
\vskip 12pt
\centerline{Figure 7.5.}
\vskip 6pt

In Section 8, we will show how to
remove all non-trivial amalgams. This, together
with Procedures 1 - 5 is enough to ensure that each 0-handle
of ${\cal F}$ has positive index.

\noindent {\bf Lemma 7.6.} {\sl Let ${\cal H}$ be a handle-decomposition
of a connected sutured manifold with non-empty boundary, containing no 
non-trivial amalgams. If some 0-handle of ${\cal F}$ has non-positive index, 
then we may apply one of Procedures 1 - 6.}

\noindent {\sl Proof.} Let $V$ be a 0-handle of ${\cal F}$
with non-positive index. Then, there are a number of cases.
\item{1.} $\vert {\cal H}^2 \cap V \vert = 0$ and
$\vert \gamma \cap V \vert = 0$.
\item{2.} $\vert {\cal H}^2 \cap V \vert = 1$ and
$\vert \gamma \cap V \vert = 0$.
\item{3.} $\vert {\cal H}^2 \cap V \vert = 1$ and
$\vert \gamma \cap V \vert = 1$.
\item{4.} $\vert {\cal H}^2 \cap V \vert = 0$ and
$\vert \gamma \cap V \vert = 2$.
\item{5.} $\vert {\cal H}^2 \cap V \vert = 0$ and
$\vert \gamma \cap V \vert = 1$.
\item{6.} $\vert {\cal H}^2 \cap V \vert = 2$ and
$\vert \gamma \cap V \vert = 0$.

For $i= 2, 3$ and $4$, we may apply Procedure $i$.
In case 1, let $H_1$ be the 1-handle containing
$V$, and let $V'$ be $\partial H_1 \cap \partial {\cal H}^0 - V$.
If we cannot apply Procedure 1 to either of the
discs $V$ or $V'$, then this component of $M$ is a 3-ball disjoint from $\gamma$,
comprised of two 0-handles joined
by a single 1-handle. We may therefore apply Procedure 5. 
Case 5 cannot arise since $\gamma$ separates $\partial M$
into ${\cal R}_-$ and ${\cal R}_+$.
In case 6, the 1-handle of ${\cal H}$ containing
$V$ is part of a non-trivial amalgam, contrary to
assumption. $\square$

\vskip 18pt
\centerline{\caps 8. Removing non-trivial amalgams}
\vskip 6pt

We now give a procedure for removing all non-trivial amalgams,
which will complete the proof of Proposition 6.3.
Suppose that there is a non-trivial amalgam ${\cal A}$
in the handle structure ${\cal H}$ of the taut sutured manifold
$(M, \gamma)$. We will assume that ${\cal A}$ is maximal, in the sense
that if any other handles are added to ${\cal A}$, the
resulting collection of handles does not form an amalgam.
We will also assume that we cannot apply any of
Procedures 1-6 in Section 7. In particular,
due to Procedure 6, this implies that $M$ has no
3-handles. Recall that ${\cal A}$ has
the structure of an $I$-bundle over a connected surface $G$.
The $I$-bundle over $\partial G$ will be denoted by
$\partial_v {\cal A}$.

Note that (iii) of 7.4 implies that $\partial_v {\cal A}$
is a union of intersections between handles of ${\cal H}$.
By (v) of 7.4, only 0-handles and 1-handles of ${\cal H} - {\cal A}$
touch $\partial_v {\cal A}$. By (iv) of 7.4, only 1-handles
and 2-handles of ${\cal A}$ touch $\partial_v {\cal A}$.
Hence, we may define a handle structure
on $\partial_v {\cal A}$ as follows. The 0-handles 
of $\partial_v {\cal A}$
arise from the intersection of 1-handles of ${\cal A}$
with the 0-handles of ${\cal H} - {\cal A}$. The 1-handles of
$\partial_v {\cal A}$ arise from the intersection of 2-handles of
${\cal A}$ with 0-handles and 1-handles of ${\cal H} - {\cal A}$.

\noindent {\bf Lemma 8.1.} {\sl Suppose that we cannot apply
Procedure 2 of Section 7. Then each 0-handle of $\partial_v {\cal A}$
abuts precisely two 1-handles of $\partial_v {\cal A}$.}

\noindent {\sl Proof.} If not, then
some 0-handle of $\partial_v {\cal A}$ abuts precisely one 1-handle
of $\partial_v {\cal A}$. This 0-handle $D$
is a component of ${\cal H}^1 \cap {\cal H}^0$.
Let $H_0$ (respectively, $H_1$) be the 0-handle
(respectively, the 1-handle) of ${\cal H}$ containing $D$. Then
$H_1$ lies in ${\cal A}$, but $H_0$ does not.
Since $D$ abuts precisely one 1-handle of $\partial_v {\cal A}$,
$H_1$ intersects ${\cal A} \cap {\cal H}^2$ in
a single disc, lying in some 2-handle $H_2$.
This in fact must the only intersection between
$H_1$ and ${\cal H}^2$, by (v) of 7.4. Hence, we may
apply Procedure 2 of Section 7 to $H_1$ and $H_2$, 
which is a contradiction. $\square$

The following lemma will also be useful.

\noindent {\bf Lemma 8.2.} {\sl Let ${\cal H}$ be
a handle structure for $(M, \gamma)$ to which
we cannot apply any of Procedures 1 - 6. Let ${\cal A}$ be
a maximal amalgam in ${\cal H}$. Let $F$ be a component
of ${\cal F}$ touching $\partial_v {\cal A}$.
Then $F$ has positive index.}

\noindent {\sl Proof.} Since we cannot apply
any of Procedures 1 - 6, the only 0-handles of
${\cal F}$ with non-positive index have valence two
and are disjoint from $\gamma$ (see the proof of Lemma 7.6).
If $F$ has non-positive index then each 0-handle of
$F$ must be of this form. Hence, $F$ is an annulus
disjoint from $\gamma$. Therefore, $F$ must be the
only component of ${\cal F}$ lying in $H_0$, where
$H_0$ is the 0-handle of ${\cal H}$ containing $F$.
For, otherwise we could apply Procedure 1.

Consider a handle of $F$ lying in $\partial_v {\cal A}$.
This is a component of intersection between $H_0$
and some 1-handle or 2-handle of ${\cal H}$.
By (iv) of 7.4, we must have $H_0 \not\in {\cal A}$.

If $F$ lies entirely in $\partial_v {\cal A}$,
then each handle of ${\cal H}$ touching $H_0$ must
be in ${\cal A}$, and so we may extend ${\cal A}$
over $H_0$, contradicting its maximality. Therefore,
$F \cap \partial_v {\cal A}$ is not the whole of $F$.

If $V$ is a 0-handle of $F$ lying in $\partial_v {\cal A}$, 
then the 1-handle of ${\cal H}$ touching $V$ must lie in
${\cal A}$. Hence, by (v) of 7.4, the 1-handles of $F$ 
touching $V$ also lie in $\partial_v {\cal A}$. Hence, we 
may find a 0-handle $F_0$ of $F$ and a 1-handle
$F_1$ of $F$ which are adjacent, with $F_1$ in
$\partial_v {\cal A}$, but $F_0$ not in $\partial_v {\cal A}$.
Let $H_1$ (respectively, $H_2$) be the 1-handle (respectively,
2-handle) of ${\cal H}$ containing $F_0$ (respectively,
$F_1$). Then, we must have $H_0 \not\in {\cal A}$,
$H_1 \not\in {\cal A}$ and $H_2 \in {\cal A}$.
Let $H'_2$ be the 2-handle other than $H_2$ which
touches $H_1$. (If $H_2$ touches $H_1$ in two discs,
then let $H'_2 = H_2$.) If $H'_2 \in {\cal A}$,
then we may extend ${\cal A}$ over $H_1$. If
$H'_2 \not\in {\cal A}$, we may extend ${\cal A}$
over $H_1 \cup H'_2$. In each case, the maximality
of ${\cal A}$ is contradicted. $\square$

We now consider the various possibilities for 
${\cal A}$ case by case.

\noindent {\sl Case 1.} ${\cal A}$ is an $I$-bundle over a
disc $G$.

In this case, we replace ${\cal A}$ with a single 2-handle $H_2$.
We attach $H_2$ to ${\cal H}^0 \cup {\cal H}^1 - {\cal A}$
using the annulus $\partial_v {\cal A}$.
We now check that ${\cal H}' = ({\cal H} - {\cal A}) \cup H_2$
is a handle structure. By (vii) of 7.4, the only requirement that
is not immediately obvious is that $H_2$ touches some
1-handle. But, if not, then
$\partial_v {\cal A}$ would have been an
annular component of ${\cal F}$, contradicting Lemma 8.2.

We now check that extended ${\cal F}$-complexity is decreased
and that (i) - (iv) of 6.3 are satisfied.
It is clear that ${\cal H}'$ respects ${\cal H}$.
This is because (for $i = 0$ and $1$) each $i$-handle of ${\cal H}'$
is an $i$-handle of ${\cal H}$ and inherits its product structure.
Of course, $H_2$ need not lie in any 2-handle of ${\cal H}$,
but this was not a requirement of Definition 6.1.
This explains why Definition 6.1 did not make more stringent requirements on 
2-handles.

If $F$ is any component of ${\cal F}$ and $F' = F \cap {\cal F}'$,
then $F'$ is either a copy of $F$, or is completely removed,
or is obtained by performing
a sequence of the following operations: remove a 0-handle
of $F$ which abuts precisely two 1-handles of $F$, and amalgamate
these two 1-handles into a single 1-handle of $F'$. 
Hence, Lemma 7.2 ensures that conditions (ii) and (iv) of 6.3 are satisfied
and also that extended ${\cal F}$-complexity is reduced.
Condition (iii) of 6.3 is clear.

We may therefore assume that ${\cal A}$ is an $I$-bundle over a surface
$G$ other than a disc. 

\noindent {\sl Case 2.} $\partial_v {\cal A} = \emptyset$.

In other words, $G$ is a closed surface. If $\partial {\cal A}$
is entirely contained in ${\cal R}_-$ or entirely contained
in ${\cal R}_+$, then it has zero Euler characteristic, since
$(M, \gamma)$ is taut and so $G$ is a torus or Klein bottle. In either case, we pick
a non-separating orientation-preserving curve in $G$, and
perform a decomposition along the $I$-bundle over this
curve. This cuts ${\cal A}$ into a solid torus.
We perform one further taut decomposition along
a product disc, ending with a single 3-ball.
We let this be a 0-handle of ${\cal H}'$.
Suppose now that $\partial {\cal A}$ intersects both ${\cal R}_-$
and ${\cal R}_+$. Then ${\cal A}$ must be a product $G \times I$.
We can then perform a sequence of decompositions along
product annuli and product discs, ending with a 3-ball, which again
we let be a single 0-handle of ${\cal H}'$.

We therefore assume that $\partial_v {\cal A}$ is non-empty.

\vfill\eject
\noindent {\sl Case 3.} Each annulus of $\partial_v {\cal A}$ 
is an incompressible product annulus.

Then by Lemma 4.2 of [10], 
the decomposition $(M, \gamma) \buildrel \partial_v {\cal A} 
\over \longrightarrow
(M', \gamma')$ is taut. We perform this decomposition. 
In other words, we separate off ${\cal A}$ from ${\cal H} - {\cal A}$,
and add sutures $\gamma'$ as appropriate.
By definition, ${\cal H} - {\cal A}$ is a handle
structure. 

The amalgam ${\cal A}$ does not
inherit a handle structure (for example, 1-handles of 
${\cal A}$ need not be attached to 0-handles
of ${\cal A}$). However, since
$\partial_v {\cal A}$ touches both ${\cal R}_-$ and ${\cal R}_+$, the
$\partial I$-bundle over $G$ cannot be connected, and
so ${\cal A}$ must a product $G \times I$. As in Case 2,
we may perform some further decompositions along
product discs, which reduce $G \times I$ to a
ball. We let this be a single 0-handle of ${\cal H}'$.

This whole procedure has the effect of removing
some components of ${\cal F}$ and also
replacing some 0-handles and 1-handles of ${\cal F}$
with arcs of $\gamma' \cap {\cal H}^0(M')$. 
An argument almost identical to that in Procedure 3
of Section 7 establishes that the hypotheses of Lemma 7.2
hold. Therefore, (i) - (iv) of 6.3 hold and
extended ${\cal F}$-complexity has been reduced.

We therefore assume that some annulus of
$\partial_v {\cal A}$ is not an incompressible
product annulus. 

\noindent {\sl Case 4.} Some annulus $A$ of 
$\partial_v {\cal A}$ is compressible in $M$.

Then $A$ compresses in $M$ to two
discs $D'_1$ and $D'_2$ with boundaries in ${\cal R}_\pm$.
Since ${\cal R}_\pm$ is incompressible in $M$
and $M$ is irreducible, $D'_1$ and $D'_2$ are parallel in 
$M$ to discs $D_1$ and $D_2$ in ${\cal R}_\pm$.
We pick $A$ so that the curve $\partial D_1$
is an innermost curve of $\partial_v {\cal A} \cap {\cal R}_\pm$
in ${\cal R}_\pm$. Since ${\cal A}$ is not an $I$-bundle
over a disc, this implies that ${\rm int}(D_1)$ is disjoint
from ${\cal A}$. The parallelity region between
the discs $D_i$ and $D'_i$ is a ball $B_i$.
Then, $B_1$ and $B_2$ are either disjoint or nested.

\noindent {\sl Case 4A.} $B_1$ and $B_2$ are disjoint.

Then, $D_1$ and $D_2$ are disjoint, and
the sphere $D_1 \cup D_2 \cup A$ bounds a ball $B$ in $M$.
Since ${\cal A} \cap B = A$, we can extend the
$I$-bundle structure of ${\cal A}$ over $B$. This contradicts
the maximality of ${\cal A}$.

\vfill\eject
\noindent {\sl Case 4B.} $B_1$ and $B_2$ are nested.

Then $B_1 \subset B_2$ and $D_1 \subset D_2$.
The component $V$ of $M - {\rm int}({\cal N}(A))$
lying wholly within $B_2$ is homeomorphic to the exterior of
a knot in $S^3$. The amalgam ${\cal A}$ lies in $V$, and
we may therefore remove $V$ from $M$ and still retain a handle structure.
This does not change the homeomorphism type of $M$ and Lemma 7.3
gives that extended ${\cal F}$-complexity decreases and that
(i) - (iv) of 6.3 hold.

\noindent {\sl Case 5.} $\partial_v {\cal A}$ is
incompressible, and some component
of $\partial_v {\cal A}$ is not a product annulus.

Now, the $\partial I$-bundle
over $G$ has at most two components. Therefore, if
some component of $\partial_v {\cal A}$ is
not a product annulus, then no component of
$\partial_v {\cal A}$ is a product annulus.
Let us suppose that $\partial_v {\cal A}$
is disjoint from ${\cal R}_-$ (say).

Pick any component $A$ of $\partial_v {\cal A}$.
Then we let $A_1$ and $A_2$ be two parallel copies of $A$,
incoherently oriented in such a way that the parallelity
region $Y$ in $M' = M - {\rm int}({\cal N}(A_1 \cup A_2))$ 
inherits four sutures. Isotope $A_1$ and $A_2$ a little
so that they become standard surfaces.
Consider the decomposition $(M, \gamma) \buildrel A_1 \cup A_2
\over \longrightarrow (M', \gamma')$. 

\noindent {\sl Case 5A.} $(M', \gamma')$ is taut.

Then we perform this decomposition. We now check
the requirements of 6.3 and also that extended
${\cal F}$-complexity has been reduced. We will
use Lemma 7.2 to do this.

Let ${\cal H}'$ be the handle structure which
$(M', \gamma')$ inherits. Consider a component $F$
of ${\cal F}$ which is altered by this decomposition,
and let $F'$ be $F \cap {\cal F}'$.
By Lemma 8.2, $F$ must have had positive index. 

Suppose that the extended ${\cal F}$-complexity of $F'$
is at least that of $F$; we aim to reach a
contradiction. We must have $C_1(X) \geq C_1(F)$ for some
component $X$ of $F'$. But each 1-handle of $F$
gives rise to precisely one 1-handle of $F'$.
Hence, $X$ must have all the 1-handles of 
$F'$. 

Each component of $A \cap F$ yields 
three discs of $F'$. Two of these discs have
no 1-handles and intersect $\gamma'$ four times.
The remaining disc has at least one 1-handle,
and has negative index. Since it has least one 1-handle,
it must be $X$, and therefore $X$ has negative index. 
Therefore $C_2(X) < C_2(F)$. Hence,
$C_{\cal F}^+(F') < C_{\cal F}^+(F)$. Lemma 7.2 now gives that
(i), (ii) and (iv) of 6.3 hold, and that extended complexity
has been reduced. Verifying (iii) of 6.3 is straightforward.

\noindent {\sl Case 5B.} $(M', \gamma')$ is not taut.

Since $M$ is irreducible, so must $M'$ be. Also, 
${\cal R}_\pm(M')$ is norm-minimising in $H_2(M', \gamma')$.
For if $S$ is any surface in $M'$ with $S \cap \partial M' = \gamma'$
and $[S, \partial S] = [{\cal R}_\pm(M'), \gamma']
\in H_2(M', \gamma')$, then $[S - Y, \partial S - Y]
= [{\cal R}_\pm(M') - Y, \gamma' - Y] =
[{\cal R}_\pm(M), \gamma] \in H_2(M, \gamma)$.
So, $\chi_-(S) \geq \chi_-(S - Y) \geq \chi_-({\cal R}_\pm(M))
=\chi_-({\cal R}_\pm(M'))$. Hence, the only way that
$(M', \gamma')$ can fail to be taut is if ${\cal R}_\pm(M')$
is compressible. This compression cannot
reduce $\chi_-({\cal R}_\pm(M'))$, as ${\cal R}_\pm(M')$ is
norm-minimising. Hence, any compressible component of ${\cal R}_\pm(M')$
is a torus or annulus. However, any circle in a compressible
annulus is homotopically trivial in $M$. In particular,
$A$ could not have been incompressible, contrary to assumption.
Thus, if ${\cal R}_\pm(M')$ is not taut, there
are three cases to consider:
\item{(i)} only one of $A_1$ and $A_2$ (say $A_1$) lies in
a compressible torus component of ${\cal R}_\pm(M')$ (called $T_1$, say)
which disjoint from $\gamma'$, or
\item{(ii)} $A_1$ and $A_2$ both lie in the same
compressible torus $T_1$ disjoint from $\gamma'$, or
\item{(iii)} $A_1$ and $A_2$ lie in distinct compressible
tori $T_1$ and $T_2$ disjoint from $\gamma'$.

Since $T_i$ is compressible and $M'$ is irreducible,
$T_i$ bounds a solid torus $V_i$ in $M'$.

In case (iii), the component of $M$ containing $A$ 
is the union of two solid tori, glued
along an essential annulus. Thus, it is a Seifert fibre
space with base space a disc and having at most two
exceptional fibres (which are the cores of the solid tori).
Also, it is disjoint from $\gamma$. Recall that, in
the statement of 6.3, we explicitly
ruled out the case where it is has two exceptional
fibres. If it has at most one exceptional fibre, it
is a solid torus, and again, we ruled this case out.

In case (ii), we pick the $A_i$ which is closest
to ${\cal A}$. Then, the orientation of
$A_i$ and ${\cal R}_+ \cap {\cal A}$ agree
near $\partial A_i$. The decomposition $(M, \gamma)
\buildrel A_i \over \longrightarrow (M_1, \gamma_1)$
is taut. Exactly as in Case 5A, this reduces
extended ${\cal F}$-complexity and (i) - (iv) of 6.3 are satisfied.

In case (i), suppose first that $M - {\rm int}({\cal N}(V_1))$ 
contains ${\cal A}$. As above, the decomposition $(M, \gamma)
\buildrel A_2 \over \longrightarrow (M_1, \gamma_1)$
is taut, and, again, this reduces
extended ${\cal F}$-complexity and (i) - (iv) of 6.3 are satisfied.

Suppose now that $V_1$ contains ${\cal A}$.
Then again the decomposition $(M, \gamma)
\buildrel A_2 \over \longrightarrow (M_1, \gamma_1)$
is taut. This time one must work a little
harder to verify that extended ${\cal F}$-complexity decreases
and that (i) - (iv) of 6.3 are satisfied.
Let $F$ be a component of ${\cal F}$, let $F' = F \cap {\cal F}(M_1)$
and let $H_0$ be the 0-handle of ${\cal H}$ containing $F$. By Lemma 8.2,
if $F$ is altered, then it must have positive
index, and so it contributes towards ${\cal F}$-complexity.
If $C_{\cal F}^+(F') \geq C_{\cal F}^+(F)$, then as in Case 5A, there
must be a single component $X$ of $F'$ containing all
the 1-handles of $F'$. Also, $X$ arises from
a component of $F \cap A$. However, unlike in Case 5A,
$X$ will not have negative index. In fact, it
will be a disc which intersects $\gamma_1$ in
four points. Hence, it has index two. Now,
$F$ has positive index and therefore its index
is at least two. Thus, if $C_{\cal F}^+(F') \geq C_{\cal F}^+(F)$,
the index of $F$ is precisely two. If
$F$ is an annulus intersecting $\gamma$ in two points, 
then $C_3(F) > C_3(X)$ and so $C_{\cal F}^+(F) > C_{\cal F}^+(F')$.
Hence, we may assume that $F$ is
a disc intersecting $\gamma$ in four points.
In this case, all but two 0-handles of $F$ have
valence two and are disjoint from $\gamma$. The
two remaining 0-handles $D_1$ and $D_2$ each have valence one.
These two handles contain a total of four points
of $\gamma \cap F$. Since $\partial_v {\cal A}$ is
disjoint from ${\cal R}_-$, each $D_i$ contains an even number of points
of $\gamma \cap F$. If one of these 0-handles contains
no points of $\gamma$, then it is a compression disc for
${\cal R}_\pm(M_1)$, which contradicts the fact
that $(M_1, \gamma_1)$ is taut. Hence,
each $D_i$ contains precisely two points of $\gamma \cap F$.
If these two points are joined by an arc of $\gamma \cap
\partial H_0$, then again ${\cal R}_\pm(M_1)$ is compressible,
which is a contradiction. Therefore, for each $i$,
the two points $\gamma \cap D_i$ are not joined
by an arc of $\gamma \cap \partial H_0$.

Now, $F'$ is $X$, together with two index zero discs.
We remove the two index zero discs using Procedure 4.
The component $X$ is a copy of $F$, and so 
$C_{\cal F}(F') = C_{\cal F}(F)$. If $F$ was not the
only component of ${\cal F}$ in $H_0$,
then these components of ${\cal F}$ have positive index, since
otherwise we can apply one of Procedures 1-5. Hence,
$n({\cal H}' \cap H_0) > 1 =n(H_0)$, where 
$n({\cal H}' \cap H_0)$ was defined in Section 5 to be the 
number of 0-handles of ${\cal H}' \cap H_0$ containing a
component of ${\cal F}'$ of positive index. This implies
that $C(H_0 \cap {\cal H}') < C(H_0)$. If $F$
was the only component of ${\cal F}$ in $H_0$,
then we have precisely the situation in the final part
of Definition 6.2. Therefore, this is a trivial modification.
This verifies (i) - (iv) of 6.3. 

We now need to check $C^+({\cal H}') < C^+({\cal H})$. But, if it is not,
then the above must happen in every 0-handle of ${\cal H}$
which is altered by the decomposition.
This implies that the component of $M - {\rm int}({\cal N}(V_1))$
containing $A$  is a solid torus 
$V_2$ with $A \cap V_2$ a single annulus in $\partial V_2$
having winding number one. But $V_1$ is a solid torus, and so this
component of $M$ is a solid torus, contradicting one of the assumptions
of 6.3. $\square$

\vskip 18pt
\centerline {\caps 9. Modifications to a decomposing surface}
\vskip 6pt

In the previous two sections, we performed a sequence of
alterations to ${\cal H}$. We are now ready to tackle Proposition 6.4.
Consider the taut decomposition $(M, \gamma) \buildrel S \over \longrightarrow
(M_S, \gamma_S)$, where $S$ is a compact oriented surface properly embedded
in $M$, having essential intersection with ${\cal R}_\pm$.
This implies that $S$ is taut and hence
incompressible. Thus, by Lemmas 4.5 and 4.9, we can assume that $S$ is in standard
form in ${\cal H}$. But, as was remarked in Section 5, 
there is a great deal of freedom over the form of $S \cap \partial {\cal H}^0$.
The aim here is to perform a series of alterations
to $S$, creating a new standard surface $S'$ which has a
considerably more restricted intersection with
$\partial {\cal H}^0$. The sutured manifold obtained by
decomposing $(M, \gamma)$ along $S'$ will be written
as $(M_{S'}, \gamma_{S'})$.

\noindent {\bf Modification 1.} Tubing along an arc.

Suppose that $\alpha$ is an arc in ${\cal R}_\pm$ with
$\alpha \cap S = \partial \alpha$. Then there is an
embedding of $\alpha \times [-1,1]$ in ${\cal R}_\pm$
with $\alpha \times \lbrace 0 \rbrace = \alpha$
and $(\alpha \times [-1,1]) \cap S = \partial \alpha \times
[-1,1]$. Suppose that the orientation that $\alpha \times [-1,1]$ inherits from
${\cal R}_\pm$ agrees with the orientation of $S$
near $\partial \alpha \times [-1,1]$. Then we
call $\alpha$ {\sl a tubing arc}. We
construct a new surface $S'$ as follows. Embed
$\alpha \times [-1,1] \times [0,1]$ in $M$ so that
$(\alpha \times [-1,1] \times [0,1]) \cap \partial M =
\alpha \times [-1,1] \times \lbrace 0 \rbrace
= \alpha \times [-1,1]$ and $(\alpha \times [-1,1] \times [0,1]) \cap S=
\partial \alpha \times [-1,1] \times [0,1]$.
Then let
$$\eqalign {S' = \ &S \cup (\alpha \times [-1,1] \times \lbrace 1 \rbrace)\cr
&\cup (\alpha \times \lbrace -1, 1 \rbrace \times [0,1])\cr
&- (\partial \alpha \times (-1,1) \times [0,1)).\cr}$$
We say that $S'$ is obtained from $S$ by {\sl tubing
along the arc $\alpha$}. Note that
if $\lbrace \ast \rbrace$ is a point in $\alpha - \partial \alpha$,
then $P = \lbrace \ast \rbrace \times [-1,1] \times [0,1]$ is
a product disc in $(M_{S'}, \gamma_{S'})$. There is a commutative diagram:

$$
\matrix{(M, \gamma) &\smash{\mathop{\longrightarrow}\limits^{S'}}
& (M_{S'}, \gamma_{S'})\cr
\Big\downarrow \rlap{$\vcenter{\hbox{$\scriptstyle S$}}$}
&&\ \Big\downarrow \rlap{$\vcenter{\hbox{$\scriptstyle P$}}$}\cr
(M_S, \gamma_S) &\smash{\mathop{=}}
& (M_S, \gamma_S)\cr}$$

\noindent Hence, $(M_S, \gamma_S)$ is taut if and only if
$(M_{S'}, \gamma_{S'})$ is taut. However, $S'$ need
not have essential intersection with ${\cal R}_\pm(M)$.

\vskip 18pt
\centerline{\psfig{figure=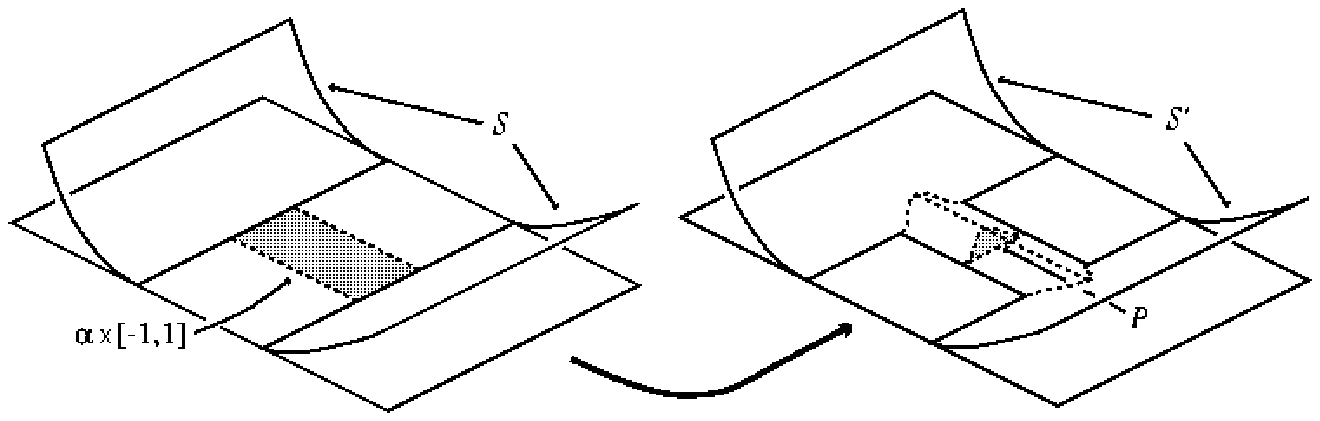}}
\vskip 6pt
\centerline{Figure 9.1.}
\vskip 12pt

\noindent {\bf Modification 2.} Slicing under an incompressible annulus.

Suppose that there is an annulus $A$ in ${\cal R}_\pm$
which is incompressible in $M$ and has $A \cap S = \partial A$.
Let $A \times [0,1]$ be embedded in $M$, so that
$(A \times [0,1]) \cap \partial M =
A \times \lbrace 0 \rbrace  = A$, and so that
$(A \times [0,1]) \cap S = \partial A \times [0,1]$.
If the orientation of $A$ agrees with that of $S$
near $\partial A$, then we construct
a new surface $S' = S  \cup (A \times \lbrace 1 \rbrace)
- (\partial A \times [0,1))$ by {\sl slicing under the incompressible annulus $A$}.
Let $C$ be a core circle of $A$. If we give $C \times [0,1]$
any orientation, then we have a commutative diagram:

$$
\matrix{(M, \gamma) &\smash{\mathop{\longrightarrow}\limits^{S'}}
& (M_{S'}, \gamma_{S'})\cr
\Big\downarrow \rlap{$\vcenter{\hbox{$\scriptstyle S$}}$}
&&\ \Big\downarrow \rlap{$\vcenter{\hbox{$\scriptstyle C
\times [0,1]$}}$}\cr
(M_S, \gamma_S) &\smash{\mathop{=}}
& (M_S, \gamma_S)\cr}$$

Since $A$ is incompressible in $M$, the annulus $C \times [0,1]$ is
incompressible in $(M_{S'}, \gamma_{S'})$,
and so by 4.2 of [10], $(M_S, \gamma_S)$ is taut if
and only if $(M_{S'}, \gamma_{S'})$ is taut.

\noindent {\bf Modification 3.} Sliding $\partial S$ across $\gamma$.

Suppose that $D$ is a disc in $\partial M$, such that
$D \cap S$ is an arc $\alpha$ in $\partial D$ and $D \cap \gamma$
is an arc properly embedded in $D$ disjoint from $\alpha$. Suppose that
the orientation of $D$ near $\partial S$ agrees with
the orientation of $S$. Let $D \times [0,1]$ be embedded in
$M$ so that $(D \times [0,1]) \cap \partial M = D \times
\lbrace 0 \rbrace = D$
and $(D \times [0,1]) \cap S = \alpha \times [0,1]$.
Let 
$$\eqalign{S' = \ &S \cup \partial (D \times [0,1])\cr
 &- (D -\partial D) \times \lbrace 0 \rbrace\cr
 &- (\alpha -\partial \alpha) \times [0,1). \cr}$$
Then $(M_{S}, \gamma_S)$ and $(M_{S'}, \gamma_{S'})$
are homeomorphic.

\vskip 18pt
\centerline{\psfig{figure=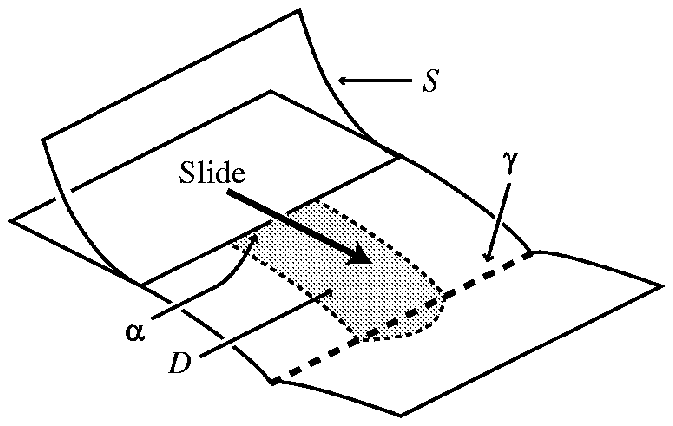}}
\vskip 6pt
\centerline{Figure 9.2.}
\vskip 12pt

\noindent {\bf Modification 4.} Slicing under a disc of contact.

Suppose that there is a disc $D$
in ${\cal R}_\pm$ with $D \cap S = \partial D$, with the
orientation of $D$ matching that of $S$ near $\partial D$.
Then $D$ is known as a {\sl disc of contact}.
Embed $D \times [0,1]$ in $M$, so that
$(D \times [0,1]) \cap \partial M = 
D \times \lbrace 0 \rbrace  = D$, and so that
$(D \times [0,1]) \cap S = \partial D \times [0,1]$.
The surface $S' = S \cup (D \times \lbrace 1 \rbrace)
- (\partial D \times [0,1))$ is obtained from $S$ by 
{\sl slicing under the disc of contact}.

\vskip 18pt
\centerline{\psfig{figure=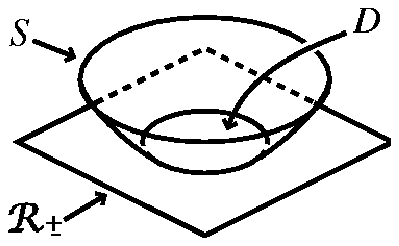}}
\vskip 6pt
\centerline{Figure 9.3.}
\vskip 12pt

Unfortunately, in this case, we have no guarantee that $(M_{S'}, \gamma_{S'})$
is taut. But the following sequence of lemmas circumvents this.
We introduce a temporary definition.

\noindent {\bf Definition 9.4.} The surface $S$ in $(M, \gamma)$
is {\sl mountainous} if some curve of $\partial S$
bounds a disc $D$ in ${\cal R}_\pm$, such that the
orientation of $D$ disagrees with that of $S$ near $\partial D$.
The disc $D$ may also intersect $S$ away from $\partial D$.

\noindent {\bf Lemma 9.5.} {\sl Suppose that $S$ is not mountainous.
Let $S'$ be obtained from $S$ by a sequence of modifications 1, 2 and 3.
Then $S'$ is not mountainous.}

\noindent {\sl Proof.} It suffices to consider the case where $S'$ is obtained
from $S$ by a single modification of type 1, 2 or 3. Consider
first a modification of type 1. Suppose that some
curve of $\partial S'$ bounds a disc $D$ in ${\cal R}_\pm$,
with orientation on $D$ disagreeing with that on $S'$. Then
at least one of the curves of $\alpha \times \lbrace -1,1 \rbrace$
must lie in $\partial D$, since $S$ was not mountainous.

If both curves of
$\alpha \times \lbrace -1, 1 \rbrace$ lie in $\partial D$,
then the curves $\partial D  \cup (\partial \alpha \times
[-1,1]) - ((\alpha - \partial \alpha)
\times \lbrace -1 ,1 \rbrace)$ would have bounded the discs $D -
((\alpha - \partial \alpha) \times [-1,1])$. Then
$S$ would have been mountainous.

Suppose now that only one curve of $\alpha \times \lbrace -1, 1 \rbrace$
lies in $\partial D$, say $\alpha \times \lbrace -1 \rbrace$.
Then $\alpha \times \lbrace 1 \rbrace$
lies in $D - \partial D$, and so is part of a
component of $\partial S'$ bounding a disc $D'$ in ${\cal R}_\pm$.
Then, $D - (D' -\partial D') - ((\alpha - \partial \alpha)
\times [-1, 1])$ is a disc which would have made $S$ mountainous.

Now consider the case where $S'$ is obtained from $S$
by slicing under an incompressible annulus $A$. This has the
effect of removing two curves of $\partial S$, neither
of which bounded discs in ${\cal R}_\pm$. Hence, in this
case, $S'$ is mountainous if and only if $S$ is mountainous.

Finally, consider the case where $S'$ is obtained from $S$
by sliding an arc of $\partial S$ across $\gamma$. Then, this
only creates new intersection points between the surface
and $\gamma$, and so a curve of $\partial S'$ disjoint from 
$\gamma$ is a copy of a curve of $\partial S$ disjoint from
$\gamma$. Thus, if $S'$ is mountainous, then so is $S$.
$\square$

\noindent {\bf Lemma 9.6.} {\sl Suppose that no component
of $S$ is a disc disjoint from $\gamma$. Let $S'$
be a surface obtained from $S$ by modifications 1, 2
and 3. Then no component of $S'$ is a disc disjoint
from $\gamma$.}

\noindent {\sl Proof.} It suffices to consider a
single modification of type 1, 2 or 3. If a
component of $S'$ which is a disc disjoint from $\gamma$
arises by tubing along an arc $\alpha$,
then the components of $S$ containing
$\partial \alpha$ were both discs disjoint from $\gamma$,
contrary to assumption. If a disc component 
of $S'$ arises by slicing under an annulus, then
that annulus could not have been incompressible in $M$. 
A type 3 modification cannot create components of $S'$
disjoint from $\gamma$. $\square$

\noindent {\bf Lemma 9.7.} {\sl Let $S$ be a
surface in a taut sutured manifold $(M, \gamma)$
having essential intersection with ${\cal R}_\pm$.
Let $S_2$ be obtained from $S$ by a sequence
of modifications 1, 2 and 3, and let $S_3$ be obtained from $S_2$
by slicing under a disc of contact $D$. Then
$S_3$ is in fact obtained from $S$ by a sequence
of modifications 1, 2 and 3.}

\noindent {\sl Proof.} We shall prove this by induction
on the number $n$ of type 1, 2 and 3 modifications
from $S$ to $S_2$. For $n = 0$, the statement of
the lemma is empty, since $S = S_2$ has essential intersection with
${\cal R}_\pm(M)$ and so
has no discs of contact.

Suppose now the lemma is true for sequences of length
at most $(n-1)$. Let $S_1$ be the surface obtained
from $S$ by the first $(n-1)$ modifications. Then
$S_2$ is obtained from $S_1$ by a modification of type
1, 2 or 3, and $S_3$ is obtained from $S_2$ by slicing
under a disc of contact $D$.

Suppose that $S_2$ is obtained from $S_1$ by tubing along
an arc $\alpha$. Then, $D$ is disjoint from 
$\alpha \times (-1,1)$. If the disc $D$ is disjoint from
$\alpha \times \lbrace -1, 1 \rbrace$, then we can
obtain $S_3$ from $S_1$ by slicing under $D$, which gives a surface $S_4$ say,
then tubing along $\alpha$. Inductively, $S_4$ is obtained
from $S$ by modifications 1, 2 and 3, and so the lemma is
true in this case. Hence, we may assume that $D$
touches at least one of the arcs $\alpha \times \lbrace -1, 1 \rbrace$.
If it touches only one arc, then $S_3$ is ambient isotopic to $S_1$, and
the lemma is true. If $D$ touches both arcs, then
$D \cup (\alpha \times [-1,1])$ is an annulus $A$
in ${\cal R}_\pm$. The two curves of $\partial A$ are
boundary components of $S_1$. If $A$ is compressible in $M$, then
both curves of $\partial A$ bound discs in ${\cal R}_\pm$,
since $(M, \gamma)$ is taut.
One of these discs has an orientation disagreeing
with that of $S_1$ near the boundary of the disc, and so $S_1$
is mountainous, contrary to Lemma 9.5.
Hence, $A$ is incompressible in $M$. We may slice under
$A$ to obtain $S_3$ from $S_1$. This proves the lemma in this case.

Suppose that $S_2$ is obtained from $S_1$ by slicing
under an incompressible annulus $A$. Then $A$ cannot lie
in $D$, since $A$ is incompressible. Hence, we can obtain
$S_3$ from $S_1$ by slicing under $D$, then slicing
under $A$. The inductive hypothesis proves the
lemma. 

Similarly, if $S_2$ is obtained from $S_1$
by sliding an arc of $\partial S_1$ across $\gamma$, then
the relevant component of $\partial S_2$ is disjoint
from $D$, and therefore we may obtain $S_3$ from $S_1$
by slicing under $D$, then performing the type 3 modification.
Again, the inductive hypothesis proves the lemma.
$\square$

\noindent {\bf Corollary 9.8.} {\sl Let $S$ be a
surface in the taut sutured manifold $(M, \gamma)$
having essential intersection with ${\cal R}_\pm$.
Then any surface obtained from $S$ by modifications 1, 2, 3 and 4
is in fact obtained from $S$ by modifications 1, 2 and 3.}

Unfortunately, if $S'$ is a surface created from $S$ by
modifications 1, 2 and 3, then $S'$ need not be incompressible.
The incompressibility of $S$ was useful in showing
that $S$ can be isotoped into standard form.
We therefore need the following lemma.

\noindent {\bf Lemma 9.9.} {\sl Let $S$ be a surface
in $(M, \gamma)$, having essential intersection with ${\cal R}_\pm$.
Suppose that $(M, \gamma) \buildrel S \over
\longrightarrow (M_S, \gamma_S)$ is taut. Let $S'$
be a surface obtained from $S$ by a sequence of modifications
1, 2 and 3. If $S'$ is in vertical form with respect
to some handle structure on $(M, \gamma)$, then we may perform an
ambient isotopy of $S'$ and perhaps some type 4 modifications,
taking $S'$ into standard form, without increasing its complexity.}

\noindent {\sl Proof.} Consider again the proof of
Lemma 4.9. The crucial property of $S$ was that
if $D$ is any disc in $M - \partial M$ with $D \cap S = \partial D$,
then $\partial D$ bounds a disc $D'$ in $S$, and
we may ambient isotope $D'$ onto $D$. In fact,
we need only restrict attention to discs $D$
lying in a single 0-handle of $M$. In the case of
$S'$, the circle $\partial D$ need not bound
a disc in $S'$, since $S'$ might not be incompressible.
But, since $(M_{S'}, \gamma_{S'})$
is taut, $\partial D$ bounds a disc $D'$ in
${\cal R}_\pm (M_{S'})$. Consider a circle $C$
of $D' \cap \partial S'$ innermost in $D'$.
By Lemma 9.6, this cannot bound a disc of
$S'$. Hence, it bounds a disc of contact in ${\cal R}_\pm(M)$.
Slice under this disc of contact.
By Corollary 9.8, the new surface (also called $S'$) is in
fact obtained from $S$ by a sequence of modifications 1, 2 and 3.
Hence, the new $(M_{S'}, \gamma_{S'})$ is taut.
So, we may repeat this process and, in this way,
we may remove all circles of $D' \cap \partial S'$.
But then $\partial D$ bounds a disc $D'$ in $S'$, and we
may ambient isotope $D'$ onto $D$. The new surface $S'$
is obtained from the old $S'$ by removing $S' \cap (D' - \partial D')$
and gluing in $D$. Thus, the complexity of the new $S'$
is no more than the complexity of the old $S'$.
Continuing in this fashion, we may get $S'$
into standard form. $\square$

We will alter $S$, using modifications 1, 2, 3 and 4, until $S$
has become a standard surface satisfying
each of the following three conditions:
\item{1.} Each curve of $S \cap \partial {\cal H}^0$ meets
any 1-handle of ${\cal F}$ in at most one arc.
\item{2.} There exists no tubing arc in $\partial {\cal F}^0
\cap {\cal R}_\pm$.
\item{3.} Suppose that $D$ is a disc in ${\cal F}^0$
with $\partial D$ the union of two
arcs $\alpha$ and $\beta$, where $\alpha = S \cap \partial D$
and $\beta = D \cap \partial {\cal F}^0$.
Suppose that one endpoint of $\alpha$ lies in ${\cal R}_\pm$
and one endpoint lies in ${\cal F}^1$. Then at least one
of the following must happen:
\itemitem{$\bullet$} $\beta$ touches at least two components
of $\partial {\cal F}^0 \cap \partial {\cal F}^1$,
\itemitem{$\bullet$} $\beta$ touches $\gamma$
and the orientation of $\alpha$ and $\beta$ disagree
locally near $\alpha \cap \beta \cap {\cal R}_\pm$, or
\itemitem{$\bullet$} $\beta$ touches $\gamma$ at
least twice, and the orientation of $\alpha$ and $\beta$ agree
locally near $\alpha \cap \beta \cap {\cal R}_\pm$.

Diagrams clarifying Conditions 1, 2 and 3 are given in
Figures 9.10, 9.11 and 9.12.
The alterations to $S$ which ensure that Conditions 1, 2 and
3 hold will reduce its complexity, and
so they are guaranteed to terminate.
The above three conditions are not quite sufficient
for our purposes. We also wish to ensure that the following 
two conditions hold:
\item{4.} Each curve of $S \cap \partial {\cal H}^0$ meets
any component of ${\cal R}_\pm \cap \partial {\cal H}^0$ in
at most one arc.
\item{5.} If $\alpha$ is an arc of $S \cap {\cal F}^0$ 
with both endpoints in ${\cal R}_\pm$,
then each of the two arcs in $\partial {\cal F}^0$
joining $\partial \alpha$ must either touch
$\partial {\cal F}^1$ or hit $\gamma_S$ more than
twice.

A diagram clarifying Condition 5 is given in Figure 9.17.
To achieve Conditions 4 and 5, we will need two further
types of modification to $S$, which we will describe later.
We now show that Conditions 1 - 3 can be achieved. By Lemma 4.9,
we may assume that $S$ is in standard form. Each alteration
to $S$ leaves it in vertical form, but not necessarily standard
form. However, we can then use Lemma 9.9 to get $S$ into standard
form, since the alterations to $S$ used there result in the
removal of some components of $S \cap \partial {\cal H}^0$,
and hence the new $S$ also satisfies Conditions 1-3.

\noindent {\bf Condition 1.} Each curve of $S \cap \partial {\cal H}^0$ meets
any 1-handle of ${\cal F}$ in at most one arc.

\vskip 20pt
\centerline{\psfig{figure=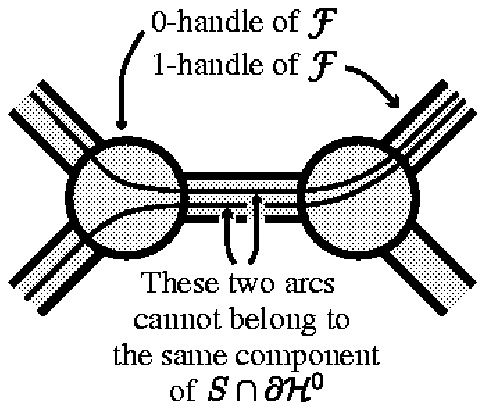}}
\vskip 6pt
\centerline{Figure 9.10.}
\vskip 10pt

Suppose that this condition does not hold. We will construct
a ball $B$ lying in $M - \partial M$, such that $B \cap S$ is a disc
in $\partial B$. We will then ambient isotope $S$ across $B$, 
and in doing so, reduce the complexity of $S$.

By assumption, there is a curve $C$ of $S \cap \partial {\cal H}^0$
containing two sub-arcs $\alpha_1$ and $\alpha_2$, which
are both properly embedded in the same 1-handle
$D_1$ of ${\cal F}$.
Pick $C$ to be a curve of $S \cap \partial {\cal H}^0$
which is innermost in $\partial {\cal H}^0$ amongst
all curves with this property.
Since $S$ is standard, there is a disc $D_2$ of
$S \cap {\cal H}^0$ with $\partial D_2 = C$.
Let $H_0$ (respectively, $H_2$) be the 0-handle (respectively, 2-handle)
containing $\alpha_1$ and $\alpha_2$, and let $E_1$ and $E_2$ be the discs of
$H_2 \cap S$ containing $\alpha_1$ and $\alpha_2$.
By the `innermost' assumption on $C$, the two arcs $\alpha_1$ 
and $\alpha_2$ are adjacent in $D_1$,
in the sense that no other arcs of $S \cap D_1$ lie between them.
Let $B'$ be the closure of the component of $H_2 - (E_1 \cup E_2)$ lying
between $E_1$ and $E_2$. The ball $B'$ will be part of $B$.

Let $\beta_1$ be an arc in the interior of $D_1$ which runs
from $\alpha_1$ to $\alpha_2$, but which intersects $S$ in
no other points. Let $\beta_2$ be an arc properly embedded
in $D_2$, with $\partial \beta_2 = \partial \beta_1$.
Then $\beta_1 \cup \beta_2$ bounds a disc $D_3$ in $H_0$.
We can assume that $D_3 - \partial D_3$ misses $S$ and $\partial H_0$.
Let $B$ be a small neighbourhood of $D_3 \cup B'$ in $M$. 
Then an ambient isotopy of $S$ across
$B$ has the effect of reducing $\vert S \cap {\cal H}^2 \vert$,
by removing the discs $E_1$ and $E_2$. The new surface
is a vertical surface with lower complexity than that of $S$.

\vfill\eject

\noindent {\bf Condition 2.}
There exists no tubing arc in $\partial {\cal F}^0 \cap {\cal R}_\pm$.

\vskip 20pt
\centerline{\psfig{figure=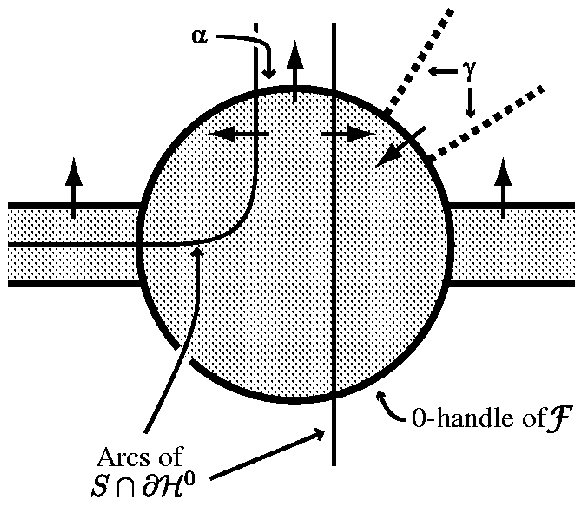}}
\vskip 6pt
\centerline{Figure 9.11.}
\vskip10pt

Suppose that $\alpha$ is such an arc. We will tube $S$
along $\alpha$. Let $H_1$
be the 1-handle of ${\cal H}$ containing $\alpha$. We may pick $\alpha \times
[-1,1]$ so that $(\alpha \times [-1,1]) \cap H_1$ is
vertical in $H_1$, and so that ${\rm cl}((\alpha \times [-1,1]) - H_1)$
is two small discs in ${\cal H}^0$. Then the surface $S'$
constructed from $S$ by tubing along $\alpha$ has
lower complexity than that of $S$, since
$\vert S \cap {\cal H}^2 \vert = \vert S' \cap {\cal H}^2 \vert$,
and $\vert \partial S' \cap {\cal H}^1 \vert
= \vert \partial S \cap {\cal H}^1 \vert - 2$.

\noindent {\bf Condition 3.}
Suppose that $D$ is a disc in ${\cal F}^0$ 
with $\partial D$ the union of two
arcs $\alpha$ and $\beta$, where $\alpha = S \cap \partial D$
and $\beta = D \cap \partial {\cal F}^0$.
Suppose that one endpoint of $\alpha$ lies in ${\cal R}_\pm$
and one endpoint lies in ${\cal F}^1$. Then at least one
of the following must happen:
\item{$\bullet$} $\beta$ touches at least two components
of $\partial {\cal F}^0 \cap \partial {\cal F}^1$,
\item{$\bullet$} $\beta$ touches $\gamma$
and the orientation of $\alpha$ and $\beta$ disagree
locally near $\alpha \cap \beta \cap {\cal R}_\pm$, or
\item{$\bullet$} $\beta$ touches $\gamma$ at
least twice, and the orientation of $\alpha$ and $\beta$ agree
locally near $\alpha \cap \beta \cap {\cal R}_\pm$.

\vskip 24pt
\centerline{\psfig{figure=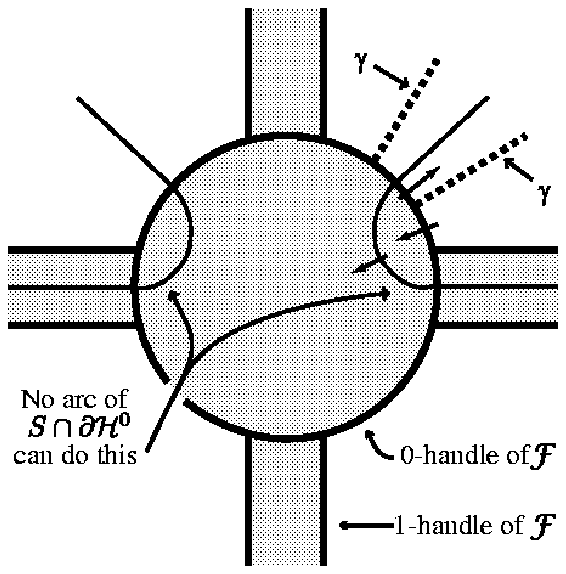}}
\vskip 8pt
\centerline{Figure 9.12.}
\vskip 4pt

Suppose that $D$ is such a disc but that it fails to satisfy
each of the three alternatives of Condition 3.
In particular, $\beta$ touches only
one component of $\partial {\cal F}^0 \cap \partial {\cal F}^1$ 
(which is therefore the component
of $\partial {\cal F}^0 \cap \partial {\cal F}^1$ 
which contains an endpoint of $\alpha$).
Suppose also that there is no sub-arc of $\beta$ which
violates Condition 2. Let $H_1 = D^2 \times I$ be the 1-handle of 
${\cal H}$ containing $D$.

There are a number of cases to consider. Suppose first
that $\beta \cap \gamma = \emptyset$.
Then, let $B$ be the vertical ball $D \times I$ in $H_1$.
Let $D'$ be the disc of $S \cap {\cal H}^2$
which touches $\alpha$ at a single point.
Let $B'$ be the closure of the component of ${\cal H}^2 - D'$ which has
non-empty intersection with $\beta - \partial \beta$.
There is an ambient isotopy of $S$ across the ball
$B \cup B'$, leaving $S$ in vertical form, and reducing
$\vert S \cap {\cal H}^2 \vert$.
This isotopy will also move parts of $S$ lying in
$D \times I$, but this causes no problems.

Suppose now that $\beta \cap \gamma \not= \emptyset$.
Then, by assumption, the orientations of $\alpha$ and
$\beta$ agree near $\alpha \cap \beta \cap {\cal R}_\pm$,
and also $\beta \cap \gamma$ is a single point.
Suppose first that there is no arc of $S \cap {\cal F}^0$ 
other than $\alpha$ lying in $D$.
Then we perform a type 3 modification, supported in a small
neighbourhood of $H_1$, which slides
the vertical arc $(\beta \cap \alpha \cap {\cal R}_\pm) \times I$ in $H_1$
across the vertical arc $(\beta \cap \gamma) \times I$.
Then, we can perform the ambient isotopy described above
to reduce the complexity of $S$.

Suppose now that there exists some arc of $S \cap \partial {\cal F}^0$
other than $\alpha$ lying in $D$.
Let $\alpha_1$ be the arc adjacent to $\alpha$.
If the sub-arc $\beta_1$ of $\beta$ lying between $\alpha$
and $\alpha_1$ is a tubing arc, then Condition 2 is violated.
If $\beta_1$ is not a tubing arc, then there are two possibilities:
$\beta_1$ touches $\gamma$, or
the orientations of $\alpha_1$ and $\beta_1$
disagree near $\alpha_1 \cap \beta_1$.
Applying this argument to each arc of $S \cap D$
with an endpoint lying between $\beta \cap \alpha \cap {\cal R}_\pm$
and $\beta \cap \gamma$, we see that these arcs
are all coherently oriented. In particular, we can
slide each of these arcs across $\gamma$.
Then, we can apply the ambient isotopy described above.

Thus, we may ensure that $S$ satisfies Conditions 1, 2 and 3.
We now ensure that $S$ also satisfies Condition 4 and 5. To do this,
we will need two further modifications. 

\noindent {\bf Modification 5.} Surgery along a product disc.

This is defined to be the reverse of a type 1 modification.

\noindent {\bf Modification 6.} Removal of a product region.

Suppose that a component $F_1$ of $S$ is parallel to
a surface $F_2$ in ${\cal R}_\pm$, and the orientations
of $F_1$ and $F_2$ disagree near $\partial F_1 = \partial F_2$.
Suppose also that the interior of the parallelity 
region between $F_1$ and $F_2$ is disjoint from
$S$. If we remove $F_1$ from $S$, creating a new surface
$S'$, then $(M_{S}, \gamma_{S})$ is homeomorphic
to $(M_{S'}, \gamma_{S'})$, plus a product component
$(F_1 \times [0,1], \partial F_1 \times \lbrace 0 \rbrace)$.
Hence, $(M_{S'}, \gamma_{S'})$ is taut if and only if
$(M_S, \gamma_S)$ is taut.

We need some lemmas to ensure that modifications 5 and 6
are well behaved.

\noindent {\bf Lemma 9.13.} {\sl Let $S$ be a surface in $(M, \gamma)$,
and let $S_2$ be a surface obtained
from $S$ by a sequence of $n$ type 5 modifications.
Let $S_3$ be obtained from $S_2$
by slicing under a disc of contact. Then in fact, $S_3$
is obtained from $S$ by a sequence of type 4 modifications,
and then at most $n$ type 5 modifications.}

\noindent {\sl Proof.} This is by induction on $n$. For $n=0$,
the lemma is trivial. Therefore, assume that the
lemma is true for less than $n$ type 5 modifications,
and let $S_1$ be the surface obtained from $S$
by the first $(n-1)$ type 5 modifications.
Let $\alpha$ be the tubing arc for $S_2$. If
the disc of contact $D$ is disjoint from $\alpha$,
then we may slice under $D$ before doing the
type 5 move, and so the lemma is true by induction.
If $D$ is not disjoint from $\alpha$, then
$D_1 \cup D_2 = D - (\alpha \times (-1,1))$ is two discs of contact
for $S_1$, and $S_3$ is obtained by slicing under
both $D_1$ and $D_2$. Applying the inductive hypothesis
twice proves the lemma.
$\square$

\noindent {\bf Lemma 9.14.} {\sl Let $S$ be a surface in $(M, \gamma)$.
Let $S_2$ be a surface obtained from $S$ by a 
sequence of type 5 modifications,
and let $S_3$ be a surface obtained from $S_2$ by a type 6 modification.
Then $S_3$ is obtained from $S$ by at most one type 6 modification,
then perhaps some type 5 modifications.}

\noindent {\sl Proof.} We will prove this by induction on the
number $n$ of type 5 modifications from $S$ to $S_2$.
For $n = 0$, the statement of the lemma is empty.
So, assume that the statement is true for sequences
of length at most $(n-1)$. Let $S_1$ be the surface
obtained from $S$ by the first $(n-1)$ type 5 modifications.
The surface $S_1$ is obtained from $S_2$ by tubing along an arc $\alpha$.
Let $F_1$ be the component of $S_2$ which we remove in
the type 6 modification. If neither component of $\partial \alpha$
lies in $F_1$, then we may obtain $S_3$ from $S_1$ by 
performing the type 6 modification, then the type 5 modification.
The lemma is then true by induction. If both components
of $\partial \alpha$ lie in $F_1$, then a single type 6
modification takes $S_1$ to $S_3$, and again the
lemma is true by induction. If a single component of
$\partial \alpha$ lies in $F_1$, then we find a (possibly empty) collection of
properly embedded arcs in $F_1$ which cut it to a disc.
These arcs define type 5 moves which can be applied to
$S_1$, at the end of which we obtain a surface ambient
isotopic to $S_3$. Hence, in this
case also, the lemma is true. $\square$

\noindent {\bf Lemma 9.15.} {\sl Let $S$ be a surface in $(M, \gamma)$
which is not mountainous.
Let $S'$ be obtained from $S$ by a sequence of type 6
modifications, then by slicing under a disc of contact $D$. Then $S'$
is obtained from $S$ by a type 4 modification, then some
type 6 modifications.}

\noindent {\sl Proof.} Suppose that some component $F_1$ of
$S$ has $\partial F_1$ lying in ${\rm int}(D)$. Then it
must be removed by some type 6 modification. In particular,
it must be parallel to a subsurface $F_2$ of $D$.
The outermost component of $\partial F_2$ in $D$ is
a component of $\partial S$ which makes $S$ mountainous,
contrary to assumption. Hence, each component of $S$ is disjoint
from ${\rm int}(D)$, and we may therefore slice under $D$
before performing the type 6 modifications. $\square$

The above lemmas all give the following proposition.

\noindent {\bf Proposition 9.16.} {\sl Let $S$ be a taut surface
in the taut sutured manifold $(M, \gamma)$, with $S$
having essential intersection with ${\cal R}_\pm$. Let $S'$
be obtained from $S$ firstly by a (possibly empty) sequence of modifications
1, 2, 3 and 4, then by a (possibly empty) sequence of 
modifications 3, 4, 5 and 6.
Then $S'$ is in fact obtained from $S$ by a (possibly empty)
sequence of modifications
1, 2 and 3, then possibly by some type 6 modifications, then
possibly some type 3 and 5 modifications. In particular, no type
4 modifications are needed.}

\noindent {\sl Proof.} Consider the sequence
of numbers from 1 to 6 which are the type of
each modification. Ignore repetitions; for
example, if we perform two type 3 modifications
in a row, then only write down one 3.
Lemma 9.13 implies that if we write down 54, we
may replace this with 45 or 4. Lemma 9.14 implies
that if we write down 56, we may instead
write down 65, 5, 6 or nothing. Also, we may replace
34 with 43, since each slide across $\gamma$ creates
a component of $\partial S$ touching $\gamma$,
whereas each slice under a disc of contact deals with
a component of $\partial S$ disjoint from $\gamma$. Hence,
we can perform the type 4 modifications before the type 3
modifications. Similarly, we can replace 36 with 63.
Hence, in the sequence of 3, 4, 5 and 6 modifications,
we can arrange to do
all the type 5 and 3 modifications last.
Corollary 9.8 asserts that we may replace the initial
sequence of 1, 2, 3 and 4 with just a sequence of 1, 2 and 3.
Let $S_1$ be the surface obtained from $S$ after this
initial sequence. Then, the sequence of numbers
after $S_1$ is a (possibly empty) sequence of 4 and 6
(starting with 6), and then possibly a sequence of 5 and 3. 
If the sequence of 6 and 4 is empty or a single
6, the proposition is proved. Otherwise, the sequence of 6 and 4
starts with 64. By Lemma 9.5, $S_1$ is not mountainous
and so, by Lemma 9.15, we may replace the
64 with 46. We may include the 4 in the initial sequence
of 1, 2 and 3. Proceeding in this way, we prove the proposition.
$\square$

We have so far modified $S$ using modifications 1, 2, 3 and 4,
resulting in a surface satisfying Conditions 1, 2 and 3. We are
now going to make some further alterations, using modifications
3, 4, 5 and 6, resulting in a surface $S'$ which also satisfies
Conditions 4 and 5. The point behind the above proposition
is that we can in fact obtain $S'$ from $S$ without slicing
under discs of contact.
 
\noindent {\bf Condition 4.} Each curve of $S \cap \partial {\cal H}^0$ meets
any component of ${\cal R}_\pm \cap \partial {\cal H}^0$ in
at most one arc.

Suppose that, on the contrary, there are two arcs $\alpha_1$ and
$\alpha_2$ of $S \cap {\cal R}_\pm \cap \partial {\cal H}^0$
properly embedded in a component of ${\cal R}_\pm \cap
\partial {\cal H}^0$, such that $\alpha_1$ and $\alpha_2$ belong to
the same component of $S \cap \partial {\cal H}^0$.
We may assume that there is an arc $\beta$ in ${\cal R}_\pm \cap
\partial {\cal H}^0$ with one endpoint in $\alpha_1$,
one endpoint in $\alpha_2$, and the remainder of $\beta$
disjoint from $S$. There is also a disc $D$ in ${\cal H}^0$
with $\partial D$ containing $\beta$, and $D \cap S =
{\rm cl}(\partial D - \beta)$.

Suppose first that $\beta$ is a tubing arc and hence
that $D$ is disjoint from $\gamma_S$.
Since ${\cal R}_\pm(M_S)$ is incompressible, $\partial D$ bounds a disc
$D'$ in ${\cal R}_\pm(M_S)$. If $\partial S \cap D'$ is
a single arc, then we may ambient isotope $S \cap D'$ onto
$D$. It is straightforward to verify that the resulting surface
$S'$ is standard and still satisfies Conditions 1 - 3.
This procedure does not increase $\vert S \cap {\cal H}^2 \vert$
and it decreases $\vert \partial S \cap {\cal H}^1 \vert$.
Hence, it decreases the complexity of $S$.

We must deal with the case where $\partial S \cap D'$ contains
some circles. Pick one $C$ innermost in $D'$, bounding a 
disc $E$ in $D'$. Then $E$ is either a disc of contact
or a disc component of $S$. In the former case, we slice
under the disc of contact. We now give a procedure for
dealing with the latter case. The curve $C$ bounds a 
disc $E'$ in ${\cal R}_\pm(M)$, and $E \cup E'$ bounds
a ball $B$ in $M$. Pick a curve $C'$ of $E' \cap \partial S$
innermost in $E'$. Then we may apply either modification 4
or modification 6 to the component of $S$ containing $C'$.
In this way, we eventually remove all components of
$S$ lying in $B$. We can then apply modification 6
to $E$. Continuing
in this fashion, we eventually remove all circles of
$\partial S \cap D'$. Then we may ambient isotope $S \cap D'$ onto $D$.

Suppose now that
$\beta$ is not a tubing arc, in which case the disc $D$ is
a product disc in $M_S$. Let $D \times [-1,1]$ be
embedded in ${\cal H}^0$, so that
\item{$\bullet$} $D \times \lbrace 0 \rbrace = D$,
\item{$\bullet$} $(D \times [-1,1]) \cap \gamma = \emptyset$,
\item{$\bullet$} $(D \times [-1,1]) \cap {\cal R}_\pm(M)
= \beta \times [-1,1]$,
\item{$\bullet$} $(D \times [-1,1]) \cap S = {\rm cl}(\partial D - \beta)
\times [-1,1]$.

\noindent Let $S'$ be $S - ({\rm cl}(\partial D -\beta) \times (-1,1))
\cup (D \times \lbrace -1,1 \rbrace)$. Let $\lbrace \ast \rbrace$
be a point in $\beta - \partial \beta$. Then $S$ is obtained
from $S'$ by tubing along the arc $\lbrace \ast \rbrace \times
[-1,1]$. Hence, $S'$ is obtained from $S$ by a type $5$ modification.
It is straightforward to check
that $S'$ still satisfies Conditions 1, 2 and 3.
This procedure leaves the complexity of $S$ unchanged.
It also creates two discs of $S' \cap {\cal H}^0$
from one disc of $S \cap {\cal H}^0$. Each of the new discs
either touches $\gamma$ or touches ${\cal H}^1$. 
But $S \cap \gamma = S' \cap \gamma$ and
$S \cap {\cal H}^1 = S' \cap {\cal H}^1$. Hence, eventually,
this process terminates.

\noindent {\bf Condition 5.} 
If $\alpha$ is an arc of $S \cap {\cal F}^0$ 
with both endpoints in ${\cal R}_\pm$,
then each of the two arcs in $\partial {\cal F}^0$
joining $\partial \alpha$ must either touch
$\partial {\cal F}^1$ or hit $\gamma_S$ more than
twice.

\vskip 18pt
\centerline{\psfig{figure=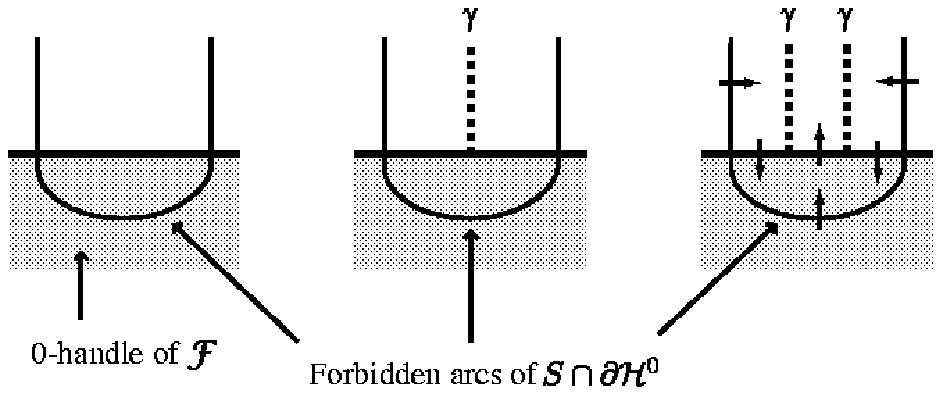}}
\vskip 6pt
\centerline{Figure 9.17.}
\vskip12pt

Let $\alpha$ be such an arc. Suppose that there is an arc $\beta$
in $\partial {\cal F}^0$
joining the endpoints of $\alpha$, such that
$\beta$ is disjoint from ${\cal F}^1$ and touches
$\gamma_S$ at most twice. Let $D$ be the disc
of ${\cal F}^0$ containing $\alpha$. Then, since Condition 2
holds, we may take $\alpha$ to be extrememost in $D$,
separating off a disc $D'$ from $D$ with $D' \cap S = \alpha$.
Then $D'$ is a disc properly embedded in $M_S$.
Hence, $\partial D'$ hits $\gamma_S$ either twice
or not at all.  If $\partial D'$ is disjoint from 
$\gamma_S$, then $\beta$ is a tubing arc, contrary to
Condition 2. 

Suppose therefore that $\partial D'$
hits $\gamma_S$ at precisely two points. If
both of these points lie at the endpoints of $\alpha$,
then we will apply modification 5. If $H_1 = D \times D^1$
is the 1-handle containing $\alpha$, then we will take
the tubing arc to be a slight extension of $\lbrace \ast
\rbrace \times D^1$, where $\lbrace \ast \rbrace$ is a
point in $\beta$. The result of this type 5 modification
is to leave $\vert S \cap {\cal H}^2 \vert$ unchanged and
to reduce $\vert \partial S \cap {\cal H}^1 \vert$ by $2$.
Hence, it reduces the complexity of $S$.

If there is a point $P$ of 
$\gamma_S \cap \partial D'$
not lying at an endpoint of $\alpha$, then 
$P$ lies on $\gamma$, and we
can perform a type 3 modification sliding 
$\alpha \times D^1$ across $P \times D^1$. 
This slide leaves the complexity of $S$ unchanged.
After possibly performing this operation once more,
we end with a situation where both points of
$\gamma_S \cap \partial D'$ lie at the endpoints of $\alpha$. Hence,
we may apply modification 5 to reduce the complexity of $S$.
It is clear that, in the above procedure, we have not violated
Conditions 1 - 3.

\vfill\eject
\centerline {\caps 10. Behaviour of handle complexity under decomposition}
\vskip 6pt

The aim of this section is to complete the proof of Proposition 6.4.
Recall that we are given a taut decomposition
$(M, \gamma) \buildrel S \over \longrightarrow (M_S, \gamma_S)$,
where $S$ has essential intersection with ${\cal R}_\pm(M)$. Also, $(M, \gamma)$ is
equipped with a handle structure ${\cal H}$,
for which each 0-handle of ${\cal F} = {\cal F}({\cal H})$
has positive index.

In the previous section, we performed a sequence of
alterations to $S$, resulting in a new standard surface (called $S'$, say)
satisfying Conditions 1 - 5. Let $(M, \gamma) \buildrel S' \over
\longrightarrow (M', \gamma')$ be the decomposition along
$S'$ and let ${\cal H}'$ be the induced handle structure
on $M'$. Note that $[S, \partial S] = [S', \partial S']
\in H_2(M, \partial M)$.

Proposition 9.16 asserted that
it sufficed to perform a sequence of modifications 1, 2 and 3, then
some type 6 modifications, then some type 3 and 5 modifications.
If a type 1 or 2 modification to $S$ results in a surface $S_1$,
then there is a pull-back $(M_S, \gamma_S) \buildrel P \over
\longleftarrow (M_{S_1}, \gamma_{S_1})$, where $P$ is a
product disc or incompressible product annulus. Hence, the
sequence of modifications 1, 2, 3 and 6 gives rise to the sequence of
pull-backs
$$(M, \gamma) = (\hat M_1, \hat \gamma_1) \buildrel P_1 \over
\longleftarrow \dots \buildrel P_{r-1} \over \longleftarrow
(\hat M_r, \hat \gamma_r)$$
which was mentioned in Proposition 6.4.
Then modifications 3 and 5 give the sequence of decompositions 
$$(\hat M_r, \hat \gamma_r) \buildrel P_r \over \longrightarrow
\dots \buildrel P_{m-1} \over \longrightarrow (\hat M_m, \hat \gamma_m)
= (M', \gamma').$$
In the case of a type 5 modification, the relevant
decomposing surface is a product disc.

We just have to check that conditions (i) - (v) of 6.4
hold if $S'$ satisfies Conditions 1 - 5. Since $S'$
is standard, (i) is trivially true. We claim that
Conditions 1, 2 and 4 guarantee (iii). Conditions 1 and 4
ensure that there is only a finite number of possibilities for each
curve of $S' \cap \partial {\cal H}^0$, up to ambient isotopy
which keeps $\gamma$ and each handle invariant.

\noindent {\bf Lemma 10.1.} {\sl Let $C$ and $C'$ be two disjoint simple
closed curves of $S' \cap \partial {\cal H}^0$ where $S'$ satisfies
Conditions 1, 2 and 4. If there is an isotopy of $\partial {\cal H}^0$
which leaves $\gamma$ and each handle of ${\cal F}$ invariant and
which takes $C$ onto $C'$, then $C$ and $C'$ are parallel in a way
which respects $\gamma$ and the handle structure on ${\cal F}$.}

\noindent {\sl Proof.} Let $\alpha$ be an arc of $C \cap
{\cal F}^0$, $C \cap {\cal F}^1$ or $C \cap {\cal R}_\pm$, and let
$\alpha'$ be the image of $\alpha$ after the isotopy taking
$C$ to $C'$. It suffices to show that no arc of $C$ or $C'$ lies
between $\alpha$ and $\alpha'$. Suppose that there is such an
arc. If $\alpha$ lies in ${\cal F}^1$, then this means that Condition 1
is violated. If $\alpha$ lies in ${\cal R}_\pm$, then Condition 4
is violated. If $\alpha$ lies in ${\cal F}^0$, then by Condition 4,
$C$ must intersect ${\cal F}^0$ in two arcs which are
joined by two arcs in ${\cal R}_\pm$. If $V$ is the 0-handle of
${\cal F}$ containing $\alpha \cup \alpha'$, then some sub-arc of
$\partial V$ is a tubing arc for $S'$, contrary to Condition 2.
$\square$

Hence, there is only a finite number of possible arrangements
for $S' \cap \partial {\cal H}^0$, up to possibly taking multiple
parallel copies of each curve and performing an ambient
isotopy which leaves $H_0$, $H_0 \cap {\cal F}$ and
$H_0 \cap \gamma$ invariant. Consider therefore 
a collection $C$ of $n$ parallel curves of $S' \cap \partial {\cal H}^0$, 
the parallelity regions respecting ${\cal F}$ and $\gamma$.
Let $H'_1, \dots, H'_{n-1}$
be the associated 0-handles of ${\cal H}'$ lying between them.
There are two possibilities: either each curve of $C$
misses $\partial {\cal F}$ or each curve of $C$ hits
$\partial {\cal F}$. In the former case, ${\cal F}' \cap H'_i$
is an annulus disjoint from $\gamma'$ for each $i$, and so
none of the $H'_i$ lie in ${\cal IH}^0(M')$, and hence can be ignored.
In the latter case, we claim that at most one
$H'_i$ lies in ${\cal IH}^0(M')$. For if two adjacent
curves of $C$ are coherently oriented, then ${\cal F}'$ intersects
the 0-handle between them in a collection of product discs.
Hence, this 0-handle does not lie in ${\cal IH}^0(M')$.
If two adjacent curves of $C$ are incompatibly oriented
(say they point towards each other), then the arcs
of $\partial {\cal F}$ lying between them must
all point out of ${\cal F}$. Otherwise, Condition 2 is
violated. Hence, at most one pair of adjacent curves
of $C$ can be incoherently oriented. In particular,
at most one $H'_i$ can lie in ${\cal IH}^0(M')$.
Therefore, (iii) of 6.4 is established.

We will now focus on a component $F$ of ${\cal F}$,
and will compare its complexity with the complexity
of $F' = F \cap {\cal F}'$.

\noindent {\bf Lemma 10.2.} {\sl Let $S'$ be a standard surface
satisfying Condition 1. Then no simple closed curve of $S' \cap F$ bounds
a disc in $F$.}

\noindent {\sl Proof.} We may pick such a simple closed curve
to be innermost in $F$, bounding a disc $D$, which
inherits a handle structure
from $F$. Since $D$ is a disc, there is a 0-handle of $D$ with valence
at most one. If this 0-handle has valence zero, then
$S'$ is not standard. If this 0-handle has valence one, then
Condition 1 is violated. $\square$

The following corollary of Lemma 10.2 is a simple
property of planar surfaces.

\noindent {\bf Corollary 10.3.} {\sl Let $S'$ be a standard 
surface satisfying Condition 1. Then one of the following holds:
\item{$\bullet$} each component $X$ of $F'$ has 
$\vert \partial X \vert < \vert \partial F \vert$, or
\item{$\bullet$} $F'$ is obtained from $F$ by cutting
along arcs and circles which are parallel to arcs and
circles in $\partial F$.

}

Condition 1 also has the following implication.

\noindent {\bf Lemma 10.4.} {\sl Let $S'$ be a standard surface
satisfying Condition 1. Then any component of $F'$
meets any 1-handle of $F$ in at most one disc.}

\noindent {\sl Proof.} Suppose, on the contrary, that there
is a component $X$ of $F'$ meeting a 1-handle $D$
of ${\cal F}$ in more
than one disc. Let $D_1$ and $D_2$ be two discs of
$X \cap D$, and let $\alpha_1$ be an arc in $X$ joining
$D_1$ to $D_2$. Let $\alpha_2$ be an arc in $D$ transverse
to $S' \cap D$, joining the endpoints of $\alpha_1$, in
such a way that $\alpha_1 \cup \alpha_2$ forms a
simple closed curve, which bounds a disc $D'$ in
$\partial {\cal H}^0(M)$. There exists at least one
circle $C$ of $S' \cap \partial {\cal H}^0$ entering $D'$ through $\alpha_2$.
This arc cannot leave $D'$ through $\alpha_1$ and so must
leave $D'$ through $\alpha_2$. Hence, $C$ violates Condition 1.
$\square$

\noindent {\bf Lemma 10.5.} {\sl Let $S'$ be a standard surface
satisfying Conditions 1, 2 and 3. Let $D$ be a component
of $F'$ with a negative index 0-handle. Then
there is a 1-handle of $F$ which touches $D$
but no other component of $F'$.}

\noindent {\sl Proof.}
The 0-handle $V$ of $D$ must be disjoint from $\gamma'$ and
have valence at most one. The boundary of $V$
is divided into $\partial V \cap \partial F$,
$\partial V \cap S'$ and at most one arc $V \cap {\cal F}^1(M')$.
If $\alpha$ is an arc of $\partial V \cap \partial F$
with both endpoints lying in $S'$, then $\alpha$ is
a tubing arc, contrary to Condition 2. 

Suppose first that $V$ has zero valence.
Then, $\partial V$ is divided into $\partial V \cap \partial F$
and $\partial V \cap S'$. However, we cannot have
an arc of $\partial V \cap \partial F$, since
its endpoints would lie in $S'$ and so would be
a tubing arc. Hence, $\partial V$ lies wholly in $S'$.
But this violates the assumption that $S'$ is standard,
which is a contradiction.

Now suppose that $V$ has valence one,
with a single 1-handle $E$ of ${\cal F}'$ attached to it.
Let $\beta_1$ and $\beta_2$ be the two
arcs of $\partial E \cap \partial D$. Then each
$\beta_i$ originally came from $\partial F$
or from $S' \cap F$.

If both $\beta_1$ and $\beta_2$ lie inside $S' \cap F$,
then the arc $\partial V - \partial E$ also lies
inside $S' \cap F$, for otherwise there would an arc of
$\partial V \cap \partial F$ with endpoints in $S'$,
violating Condition 2. But if $\partial V - \partial E$
lies wholly in $S'$, then Condition 1 is violated.
Similarly, if $\beta_1$ lies inside $S' \cap
F$ and $\beta_2$ lies inside $\partial F$, then
Condition 3 is violated. If both $\beta_1$ and $\beta_2$
lie inside $\partial F$, then $E$ is the required 1-handle of $F$
lying solely in $D$. $\square$

An example of a component $D$ of $F'$ with $I(D) < 0$
is given in Fig. 5.2.

\noindent {\bf Proposition 10.6.} {\sl Let $F$ be a
component of ${\cal F}$, and let $F' = F \cap {\cal F}'$.
Suppose that every 0-handle of $F$ has positive index.
If $S$ is a standard surface
satisfying Conditions 1 - 5, then
$C_{\cal F}(F') \leq C_{\cal F}(F)$. Also, if
we have equality, then each component of $S' \cap F$ is
a circle parallel to a component of $\partial F$ disjoint from
$\gamma$. The parallelity region inherits a handle structure
from ${\cal F}$ in which each 0-handle has valence two.}

\noindent {\sl Proof.} Suppose
that $C_{\cal F}(F') \geq C_{\cal F}(F)$. By Lemma 10.4, each
component $Y$ of $F'$ has $C_1(Y) \leq C_1(F)$.
If some component $D$ of $F'$ has negative index,
then some 0-handle of $D$ has negative index
and so by Lemma 10.5, all remaining components $Y$ of $F'$
have $C_1(Y) < C_1(F)$. But by definition $D$ does not
contribute to the ${\cal F}$-complexity of $F'$.
Hence $C_{\cal F}(F') < C_{\cal F}(F)$, which is contrary to assumption.
Thus no component of $F'$ has negative index.
Hence, the index of $F$ is shared among the components
of $F'$. Since $C_{\cal F}(F') \geq C_{\cal F}(F)$, then one
component $X$ of $F'$ has $C_1(X) = C_1(F)$ and $I(X) = I(F)$.
All other components $Y$ of $F'$ have zero index,
and so, by definition, they do not contribute to the
${\cal F}$-complexity of $F'$.
By Corollary 10.3, $\vert \partial X \vert \leq \vert \partial F
\vert$. Hence, $C_{\cal F}(F') \leq C_{\cal F}(F)$, 
and so $C_{\cal F}(F') = C_{\cal F}(F)$.

We now wish to examine further the case when
$C_{\cal F}(F') = C_{\cal F}(F)$. Since $\vert \partial X \vert
= \vert \partial F \vert$, Corollary 10.3 implies that $F'$ is obtained
from $F$ by cutting along arcs and circles which are
parallel to arcs and circles in $\partial F$.
Each component of $F' - X$ has index zero.
If $V$ is a 0-handle of $F'$ not lying in $X$, then
$V$ cannot have negative index. For otherwise,
Lemma 10.5 would imply that $C_1(X) < C_1(F)$. Thus,
each 0-handle $V$ of $F' - X$ must have zero index. 

We will now show that in fact there are no index zero
discs of $F'$. If there is such a disc, then there is
an arc of $F \cap S'$ extrememost in $F$, parallel to an
arc in $\partial F$ via a parallelity disc $D$. If this
disc $D$ does not have zero index, then $D = X$ and hence
$F$ is a disc. In this case, we may find an arc of
$F \cap S'$ which is extrememost in $F$ and which
does separate off an index zero disc. Thus, we may
assume that $D$ has zero index. Let $V$ be a 0-handle
of $D$ with valence at most one. Since $V$ has
zero index, there are two cases to consider.
If $V$ has valence zero and hits $\gamma'$
twice, then Condition 5 is violated.
If $V$ has valence one and hits $\gamma'$
once, then Condition 3 is violated. 

Hence, each component of $F \cap S'$ is a simple closed curve
parallel to a curve of $\partial F$. An
extrememost component of $F \cap S'$ separates off an annulus $A$.
If $A = X$, then $F$ is an annulus, in which case we may find an
arc of $F \cap S'$ extrememost in $F$ which is parallel to
a component of $\partial F$ via a component of $F'$ other than $X$.
Hence, we may assume $A \not= X$. Therefore, $A$ has zero index,
and so each 0-handle of $A$ has valence two and is disjoint from
$\gamma$. Repeating this argument for each component of $F \cap S'$
proves the proposition. $\square$

The following verifies (ii), (iv) and (v) of Proposition 6.4
and completes the proof of that proposition and hence of
Theorems 1.4, 1.5 and 1.6.

\noindent {\bf Proposition 10.7.} {\sl Suppose that every
0-handle of ${\cal F}$ has positive index.
Suppose also that $H_0 \cap ({\cal F}(M) \cup \gamma)$ is
connected for each 0-handle $H_0$ of ${\cal H}(M)$.
Let $S'$ be a standard surface satisfying Conditions 1 - 5, 
with $[S', \partial S'] \not= 0 \in H_2(M, \partial M)$.
Then $C(H_0 \cap {\cal H}') \leq C(H_0)$ for each 0-handle 
$H_0$ of ${\cal H}$. Also, this inequality is strict for
some 0-handle $H_0$. If this inequality is an equality for some
0-handle $H_0$, then $H_0 \cap {\cal H}'$ is obtained from $H_0$
by a trivial modification.}

\noindent {\sl Proof.} Suppose that $C(H_0 \cap {\cal H}') \geq C(H_0)$
for some 0-handle $H_0$ of ${\cal H}$.
Then $C_{\cal F}(H_0 \cap {\cal H}') \geq C_{\cal F}(H_0)$.
But, by Proposition 10.6, each component $F$
of ${\cal F}$ has $C_{\cal F}(F \cap {\cal F}') \leq C_{\cal F}(F)$. 
Hence $C_{\cal F}(H_0 \cap {\cal H}') \leq C_{\cal F}(H_0)$.
Therefore, for each component $F$ of ${\cal F} \cap H_0$,
we must have $C_{\cal F}(F \cap {\cal F}') = C_{\cal F}(F)$.
Proposition 10.6 then implies
that each component of $S' \cap {\cal F} \cap H_0$ is
a circle parallel to a component of $\partial {\cal F}$
disjoint from $\gamma$. Therefore, $n({\cal H}' \cap H_0) \geq
n(H_0) = 1$. Hence, $C({\cal H}' \cap H_0) \leq C(H_0)$.

Suppose that this is an equality for some 0-handle $H_0$
of ${\cal H}$. Then, as above, this implies
that $n({\cal H}' \cap H_0) = n(H_0) = 1$.
Also, the argument above gives that each component 
of $S' \cap {\cal F} \cap H_0$ is
a circle parallel to a component of $\partial {\cal F}$
disjoint from $\gamma$. This component of
$\partial {\cal F}$ bounds a disc in $\partial H_0$
with interior disjoint from ${\cal F}$,
since $H_0 \cap ({\cal F}(M) \cup \gamma)$ is connected.
Hence, $H_0 \cap {\cal H}'$ is obtained from $H_0$ by a trivial
modification.

Suppose now that $C(H_0 \cap {\cal H}') = C(H_0)$ for every 0-handle
$H_0$ of ${\cal H}$. We aim to achieve a contradiction.
Let $C$ be the collection of circles extrememost in ${\cal F}$.
Then, there is a collection of annuli $A$ in ${\cal F}$ which
is disjoint from $\gamma$ and with $A \cap S' = C$. Let $D$ be the
collection of discs of $S' \cap {\cal H}^0$
which $C$ bounds. Then $A \cup D$ is
a collection of discs properly embedded in $M$
which are parallel to discs in ${\cal R}_\pm$ via balls $B_0$.
These balls lie in ${\cal H}^0$ since $H_0 \cap ({\cal F}(M) \cap 
\gamma)$ is connected for each 0-handle $H_0$ of ${\cal H}(M)$.

For each 1-handle $H_1 = D^2 \times [0,1]$, the discs 
$D^2 \times \lbrace 0 \rbrace$ and $D^2 \times
\lbrace 1 \rbrace$ are each divided up by the decomposition
along $S'$. For $i=0$ and $1$, all but one 0-handle of $D^2 \times 
\lbrace i \rbrace  - {\rm int}({\cal N}(S'))$
has index zero. The remaining component
has index equal to the index of $D^2 \times \lbrace i \rbrace$. But the index
of $D^2 \times \lbrace i \rbrace$ is positive, since we are assuming that
the index of each 0-handle of ${\cal F}$ is positive.
Hence, the product structure on $H_1$ matches $A \cap 
(D^2 \times \lbrace 0 \rbrace)$ with $A \cap (D^2 \times \lbrace 1 \rbrace)$. 
We may therefore unambiguously define $A \cap D^2$.
Let $B_1$ be the union (over all 1-handles) of
the balls $(A \cap D^2) \times [0,1]$.
Similarly, we may find a collection $B_2$ of
components of ${\cal H}^2 - {\rm int}({\cal N}(S'))$,
and such that $B_2 \cap {\cal H}^0 = {\cal H}^2 \cap A$.

Then $B_0 \cup B_1 \cup B_2$ is a parallelity region
between some closed components of $S'$ and a subsurface of
${\cal R}_\pm$. If we remove these components, we may
repeat the argument, and show eventually that each component
of $S'$ which touches ${\cal F}$ is closed and parallel to 
some component of ${\cal R}_\pm$. This does not quite show
that $[S', \partial S'] = 0 \in H_2(M, \partial M)$,
since there may be components of $S'$ which are disjoint
from ${\cal F}$. Such a component $X$ lies entirely in a
0-handle $H_0$ of ${\cal H}$. But recall from above that
$n({\cal H}' \cap H_0) = n(H_0) = 1$. Hence, $\partial X$
cannot separate components of $H_0 \cap {\cal F}$.
In particular, $X$ is parallel to a disc in $\partial M$.
Therefore, $[S', \partial S'] = 0 \in H_2(M, \partial M)$,
contrary to assumption. $\square$

\vfill\eject
\centerline{\caps 11. The algorithm to construct the tangles}
\vskip 6pt

In this section, we demonstrate how to construct the graphs
$G$ required for algorithm of Theorem 1.4 which we outlined
in Section 2. Recall that each graph $G$ is embedded in
a 3-simplex $\Delta^3$ and comes with a regular neighbourhood
${\cal N}(G)$ and arcs labelled $\gamma$ and $\tau$ in
$\partial {\cal N}(G)$. Recall that the arcs $\gamma$ form
the tangles required for Theorems 1.5 and 1.6.

In line with the rest of this paper, we work with the
handle structure ${\cal H}$ arising by dualising the
given generalised triangulation of $M$. We will focus on
a single 0-handle $H_0$ of ${\cal H}$. The algorithm starts
with the 0-handle $H_0$ and the surface ${\cal F}(M) \cap H_0$.
This surface is one of finitely many possibilities, but, for
the moment, we will assume that ${\cal F}(M) \cap H_0$ is as
in Figure 11.1. In general, it may be a subsurface of this; we will
explain later how to cope with this eventuality. 

\vskip 18pt
\centerline{\psfig{figure=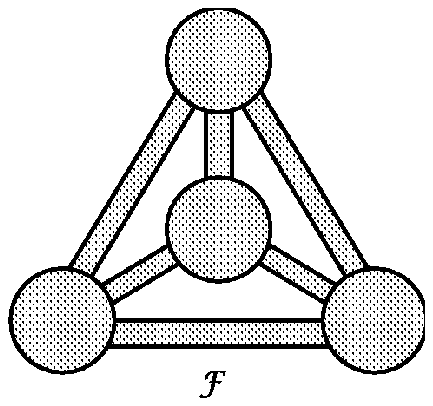}}
\vskip 6pt
\centerline{Figure 11.1.}
\vskip12pt

At each stage $j \in {\Bbb N}$ of the algorithm, we will
construct a finite list of possibilities for
the following objects lying in $H_0$:
\item{$\bullet$} a subset ${\cal H}_j^0$ of $H_0$, which is
a union of 3-balls embedded in $H_0$,
\item{$\bullet$} a subsurface ${\cal F}({\cal H}_j)$ of 
${\cal H}_j^0 \cap {\cal F}(M)$, and
\item{$\bullet$} arcs $\gamma_j$ properly embedded in
${\rm cl}(\partial {\cal H}^0_j - {\cal F}({\cal H}_j))$.

\noindent Each component of $\partial {\cal H}_j^0 -
({\cal F}({\cal H}_j) \cup \gamma_j)$ will have a
specified orientation, pointing into or out of
${\cal H}_j^0$. When we wish to refer to the above
data, we will denote it simply by ${\cal H}_j$.

For $j=1$, we take ${\cal H}_1^0 = H_0$,
${\cal F}({\cal H}_1) = H_0 \cap {\cal F}(M)$,
and $\gamma_1 = \emptyset$. We consider all
possible orientations for $\partial {\cal H}_1^0
-{\cal F}({\cal H}_1)$. If (as we supposed above)
$H_0 \cap {\cal F}(M)$ is as in Figure 11.1, then
there are four components of $\partial {\cal H}_1^0
-{\cal F}({\cal H}_1)$, and so there are 16 possible
orientations.

The algorithm constructs the list of possibilities
for ${\cal H}_{j+1}$ by considering each
possibility for ${\cal H}_j$ in turn,
and performing some modifications to it, which
we describe below. These modifications have the
property that, if ${\cal H}_j$ is some fixed
possibility at the $j^{\rm th}$ stage, then
each possibility for ${\cal H}_{j+1}$ to which it
gives rise satisfies one of the following:
\item{$\bullet$} $C({\cal H}_{j+1}) < C({\cal H}_j)$ or
\item{$\bullet$} $C({\cal H}_{j+1}) = C({\cal H}_j)$ and
$C^+({\cal H}_{j+1}) < C^+({\cal H}_j)$.

\noindent Thus, by Lemma 5.3, the algorithm will terminate
at stage $m$, say. However, we do not know this
value of $m$ until we run the algorithm.

It should be clear that this algorithm is modelling
within the single 0-handle $H_0$ what is happening
in the proof of Theorem 1.4. Recall that, in that proof,
we constructed a sequence of sutured manifolds embedded
within $M$, and examined how each sutured manifold
$(M_i, \gamma_i)$ intersected any given 0-handle $H_0$.
However, the intersection $M_i \cap H_0$ does
not necessarily correspond precisely with the $i^{\rm th}$ stage
of the algorithm we are about to outline. This
is because, when passing from a single possibility 
for ${\cal H}_j$ to several possibilities for ${\cal H}_{j+1}$,
we insist that complexity or extended complexity strictly decreases 
within our given 0-handle. However, at each stage in the induction of
Theorem 1.4, we merely insisted that complexity decreased
within {\sl some} 0-handle (not necessarily the one
we are examining). Therefore, in order to determine
how the final sutured manifold $M_n$ lies in $H_0$,
we must consider every possibility for ${\cal H}_j$,
where $1 \leq j \leq m$. Given one such possibility ${\cal H}_j$,
we construct the graph $G$ by associating a
vertex of $G$ with each component of ${\cal H}_j^0$;
we associate an edge of $G$ with each component of
${\cal F}^0({\cal H}_j)$; the curves $\tau$ are
specified by ${\cal F}^1({\cal H}_j)$;
the arcs $\gamma$ are formed by taking all
possible subtangles of $\gamma_j$.

We now give the heart of the algorithm, namely
the procedure which constructs each possibility
for ${\cal H}_{j+1}$ arising from a single possibility for
${\cal H}_j$. We apply one of the following procedures
to ${\cal H}_j$ (and we give the points in Sections 7 - 9
where we applied them).

\noindent 1. Removal of a component of ${\cal H}_j^0$.

This can occur in Procedures 1 and 5 of Section 7.
It can also occur in Cases 1, 2, 3 and 4B of Section 8.
When a component of ${\cal H}^0_j$ is removed, so are the
components of ${\cal F}({\cal H}_j)$ and arcs of $\gamma_j$
which it contains.

\noindent 2. Removing some handles of ${\cal F}({\cal H}_j)$
disjoint from $\gamma_j$.

In order that the new surface inherits a handle structure,
we insist that if a 0-handle of ${\cal F}({\cal H}_j)$
is removed, then so are the 1-handles of ${\cal F}({\cal H}_j)$
which abut it. Also, we may only perform this procedure
if the components of $\partial {\cal H}_j^0
- ({\cal F}({\cal H}_j) \cup \gamma_j)$ which touch any
removed handle have orientations which agree.
This operation can occur in Procedures 1, 2, 5 and 6
of Section 7, and Case 4B of Section 8.

Note that, in general, ${\cal F}(M) \cap H_0$ is obtained
from the surface in Figure 11.1 by removing some handles.
Therefore, by applying this procedure at the first stage
$j=1$, we can incorporate all possibilities for
${\cal F}(M) \cap H_0$ into this algorithm.

\noindent 3. Replacing handles of ${\cal F}({\cal H}_j)$ with a
sub-arc of $\gamma_{j+1}$.

Here, we may replace a 1-handle of ${\cal F}({\cal H}_j)$ with an
arc of $\gamma_{j+1}$, providing that the components of 
$\partial {\cal H}_j^0 - ({\cal F}({\cal H}_j) \cup \gamma_j)$ 
which touch this handle have orientations which disagree.
We may also remove a 0-handle of ${\cal F}({\cal H}_j)$
which has valence one and which intersects $\gamma_j$
in a single point, providing that we also remove the
1-handle of ${\cal F}({\cal H}_j)$ which it abuts
and we then replace these handles with a sub-arc of
$\gamma_{j+1}$. This occurs in Procedure 3 of Section 7
and in Case 3 of Section 8.

\noindent 4. Removal of a product disc component of 
${\cal F}({\cal H}_j)$.

If $F$ is a disc component of ${\cal F}({\cal H}_j)$
intersecting $\gamma_j$ twice, we may replace $F$ with
an arc of $\gamma_{j+1}$ joining the two points
of $F \cap \gamma_j$. This occurs in Procedure 4
of Section 7.

\noindent 5. Removal of a valence two 0-handle of ${\cal F}({\cal H}_j)$

If $V$ is a 0-handle of ${\cal F}({\cal H}_j)$ which is disjoint 
from $\gamma_j$ and which abuts two distinct 1-handles of
${\cal F}({\cal H}_j)$, then we may combine $V$ and the
two 1-handles into a single 1-handle of ${\cal F}({\cal H}_{j+1})$.
This occurs in Case 1 of Section 8.

\noindent 6. Decomposition along a surface.

This step models the sutured manifold decomposition outlined
in Section 9. We only perform this operation providing each
0-handle of ${\cal F}({\cal H}_j)$ has positive index
and providing $H_0 \cap ({\cal F}({\cal H}_j) \cup \gamma_j)$
is connected for each 0-handle $H_0$ of ${\cal H}_j$.
We construct all possible oriented curves $C$ which satisfy
Conditions 1 - 5 of Section 8 (viewing $C$ as a
possibility for $S' \cap \partial {\cal H}^0$).
There is only a finite number of possibilities
($C_1, \dots, C_t$, say) for $C$. We then let
$C'$ be a collection of disjoint simple closed curves in
$\partial {\cal H}^0_j$, each curve being a copy of one of the
$C_i$'s, and with no two components of $C'$ representing
the same $C_i$ (although, two components of $C'$ may
be the same underlying curve, but have opposite orientations).
We insist that $C'$ also satisfies Condition 2 of Section 8.
We then extend $C'$ to a collection of disjoint
discs properly embedded in ${\cal H}^0_j$. We
then decompose ${\cal H}^0_j$ along these discs,
creating a new collection of 0-handles ${\cal H}^0_{j+1}$,
which naturally inherit ${\cal F}({\cal H}_{j+1})$
and sutures $\gamma_{j+1}$. By the argument of Propositions 10.6 and 10.7, 
$C({\cal H}_{j+1}) < C({\cal H}_j)$, unless each component
of $C'$ is a curve lying entirely in ${\cal F}({\cal H}_j)$
parallel to some component of $\partial {\cal F}({\cal H}_j)$
disjoint from $\gamma_j$,
the parallelity region respecting the handle structure
of ${\cal F}({\cal H}_j)$. In this case, the modification
is trivial. We therefore do not include this case as
a possibility for ${\cal H}_{j+1}$. However, the modification
may alter the orientations of some components of
$\partial {\cal H}^0_j -{\cal F}({\cal H}_j)$ disjoint
from $\gamma_j$. We therefore have to consider all
possible orientations for these components as giving distinct possibilities
for ${\cal H}_j$.

\noindent 7. Amalgam removal

We only perform this operation when each component of 
${\cal F}({\cal H}_j)$ has positive index and each 0-handle
of ${\cal F}({\cal H}_j)$ with non-positive index has
valence two and is disjoint from $\gamma_j$.
Suppose that $D$ is a union of handles of ${\cal F}({\cal H}_j)$
which forms a disc disjoint from $\gamma$. Suppose also
that if $D$ has any 0-handles, then each such 0-handle 
abuts precisely two 1-handles
of ${\cal F}({\cal H}_j)$, both of which lie in $D$.
Suppose also that the two components of
$\partial {\cal H}^0_j - ({\cal F}({\cal H}_j) \cup \gamma_j)$
which touch $D$ have the same orientation.
We take one or two copies of $\partial D$ and move
them a little, creating a curve $C_1$ (and possibly $C_2$)
which intersect ${\cal F}({\cal H}_j)$
in a collection of arcs lying in 0-handles of ${\cal F}({\cal H}_j)$.
Extend each $C_i$ to a disc $D_i$ properly embedded in ${\cal H}^0_j$.
If we have both $D_1$ and $D_2$, we orient them inconsistently,
in a way which gives the parallelity region between them
some sutures. If we are just dealing with $D_1$, we consider
both possible orientations, providing that we are not dealing with the
final case of Definition 6.2. We then decompose ${\cal H}^0_j$
along $D_1$ (and possibly $D_2$), as outlined in operation 6 above.
This occurs in Case 5 of Section 8.
Note that we cannot necessarily include this case here
in operation 6, since the curves $C_1$ and $C_2$ might fail
Conditions 2, 4 or 5 of Section 9. 

The following figure gives a concrete example
of some of the above operations.
It is clear that these procedures may implemented algorithmically,
although they may pose some challenges for a computer programmer.

\vskip 24pt
\centerline{\psfig{figure=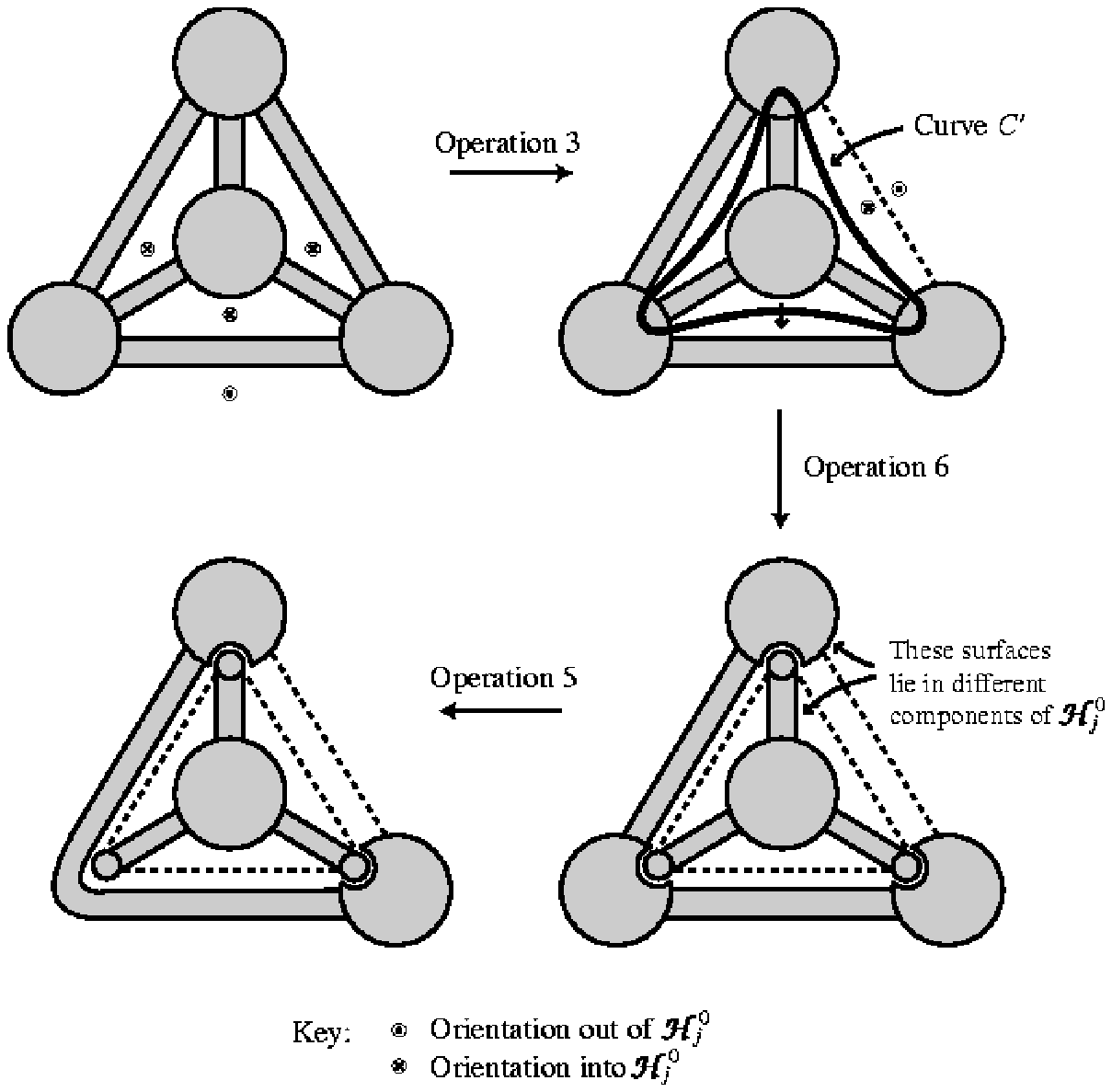}} 
\vskip 18pt
\centerline{Figure 11.2.}

\vfill\eject
\centerline{\caps 12. Exceptional and norm-exceptional surgeries 
with $\Delta(\sigma, \mu) = 1$.}
\vskip 6pt

We now give examples which demonstrate that the
restriction on $\Delta(\sigma, \mu)$ in Theorems 1.4, 1.5 and 1.6
is necessary.
We give a method of constructing in a 3-manifold $M$ (satisfying
certain conditions) an infinite number of surgery curves $K$
with exceptional or norm-exceptional surgery slopes $\sigma$ satisfying
$\Delta(\sigma, \mu) = 1$, where $\mu$ is the
meridian slope on $\partial {\cal N}(K)$. 

Let $M$ be a compact orientable 3-manifold with $\partial M$
a (possibly empty) union of tori. Suppose also
that $M$ is irreducible, atoroidal and has incompressible boundary.
Let $S$ be a connected oriented surface properly embedded
in $M$ with $[S, \partial S] \not= 0 \in H_2(M, \partial M)$,
and so that $S$ is incompressible and norm-minimising in its homology class.
Then $S$ is neither a sphere nor a disc.
Let $K$ be any essential simple closed curve on $S$ disjoint from $\partial S$.
Let $\sigma$ be the slope of the curves $\partial {\cal N}(K) \cap S$,
which is known as `surface framing'.

\noindent {\bf Proposition 12.1.} {\sl The slope
$\sigma$ is exceptional or norm-exceptional.}

\noindent {\sl Proof.} The surface $S$ determines
a class $z \in H_2(M - {\rm int}({\cal N}(K)), \partial M)$
as follows. The two curves $S \cap \partial {\cal N}(K)$ divide
$\partial {\cal N}(K)$ into two annuli. Attach either
of these annuli to $S - {\rm int}({\cal N}(K))$ and
let $S'$ be the
resulting surface. Then $z = [S', \partial S'] \in 
H_2(M - {\rm int}({\cal N}(K)), \partial M)$ is
independent of the choice of annulus. In fact,
$S'$ is norm-minimising in its class in
$H_2(M - {\rm int}({\cal N}(K)), \partial M)$,
since $\chi_-(S') = \chi_-(S) = x([S, \partial S]) 
\leq x([S', \partial S'])$.
Let $z_\sigma \in H_2(M_K(\sigma), \partial M_K(\sigma))$ be the image 
of $z$ under the map induced by inclusion.

We may construct a surface $S_\sigma$ in
$M_K(\sigma)$ by starting with $S - {\rm int}({\cal N}(K))$
and attaching a disc to each curve of $S \cap \partial {\cal N}(K)$,
the discs being meridian discs in the surgery solid torus.
Then $[S_\sigma, \partial S_\sigma] = z_\sigma \in 
H_2(M_K(\sigma), \partial M_K(\sigma))$. Also,
$-\chi(S_\sigma) = -\chi(S') - 2$. Since we assumed that
$S$ was connected, there are two possibilities: 
\item{(i)} $\chi_-(S_\sigma) < \chi_-(S')$, or 
\item{(ii)} $\chi_-(S_\sigma) = \chi_-(S') = 0$.

\noindent In case (i), $x(z_\sigma) < x(z)$, and therefore
$K$ and $\sigma$ are norm-exceptional.
In case (ii), $S_\sigma$ is a non-separating
sphere or two non-separating discs in $M_K(\sigma)$. If $S_\sigma$ is a
sphere, then $M_K(\sigma)$ is reducible. If $S_\sigma$
is two discs, then $M_K(\sigma)$ has compressible boundary,
which implies that either $M_K(\sigma)$ is a solid
torus or it is reducible. Thus, in this case, $K$ and $\sigma$ are
exceptional. $\square$

We now show that one may find an infinite number of
such knots $K$ on a given $S$ (satisfying some conditions),
such that no two knots in this collection are ambient
isotopic to each other in $M$.

\noindent {\bf Proposition 12.2.} {\sl Let $M$ be a
compact 3-manifold with $\partial M$ a (possibly empty)
union of tori and with $H_1(M)$ torsion free. Let $S$ be a compact 
connected oriented surface properly embedded in $M$ which has positive
genus and which is norm-minimising in its class in
$H_2(M, \partial M)$. Then, we may find an infinite
collection of knots, each essential curves on $S$, no two of which
are ambient isotopic in $M$.}

\noindent {\sl Proof.} We may find two simple closed
curves $C_1$ and $C_2$ on $S$ which intersect each other
precisely once. If one of these curves has infinite order
in $H_1(M)$ ($C_1$, say), then, for any integer $n$,
consider the curve $n C_1 + C_2$, which is constructed by taking
$n$ (coherently oriented) parallel copies of $C_1$, together
with $C_2$ and smoothing off the double-points. 
This is the required collection of knots on $S$.

Suppose therefore that $C_1$ and $C_2$ have finite order
in $H_1(M)$. Since $H_1(M)$ is torsion free,
this implies that $C_1$ and $C_2$ are homologically trivial.
We will construct our collection of knots
by analysing the `Seifert form' on $S$. Given two disjoint
homologically trivial closed curves $\alpha_1$ and $\alpha_2$
in $M$, define their linking number $lk(\alpha_1, \alpha_2)$ to be
the signed intersection number
between $\alpha_2$ and a (not necessarily embedded) Seifert surface for $\alpha_1$.
This is independent of the choice of Seifert surface
for $\alpha_1$, since any two Seifert surfaces can be
glued to form a closed (not necessarily embedded) surface, with which
$C_2$ has zero intersection, since it is homologically trivial. 
Also, it is symmetric: $lk(\alpha_1, \alpha_2) = lk(\alpha_2,
\alpha_1)$. Given any curve $C$ on $S$, define $C^+$ to be the
push-off of $C$ from $S$ in some specified normal direction.
Define the framing $fr(C)$ of any curve $C$ on $S$ which
is homologically trivial in $M$ as $lk(C^+, C)$.
Now, 
$$lk(C^+_1, C_2) - lk(C_1, C^+_2) = \pm 1,$$
since $C_1$ and $C_2$ intersect in one point on $S$.
This implies that
$$lk(C^+_1, C_2) + lk(C_1, C^+_2) \not= 0.$$
Let $n_1$ be an arbitrary integer. Then
$$\eqalign{& fr(n_1 C_1 + C_2) =
lk((n_1 C_1 + C_2)^+, n_1 C_1 + C_2) \cr
&=n_1^2 lk(C^+_1, C_1) +
n_1(lk(C^+_2, C_1) + lk(C^+_1, C_2))
+ lk(C^+_2, C_2).\cr
&= n_1^2 k_1 + n_1 k_2 + k_3,\cr}$$
for integers $k_1$, $k_2$ and $k_3$, where
$k_2 = lk(C_1, C^+_2) + lk(C^+_1, C_2) \not= 0$.
Hence, $fr(n_1 C_1 + C_2)$ takes infinitely many values.

We now claim that if $C$ and $C'$ on $S$ are two closed curves
on $S$ which are homologically trivial in $M$ and freely
homotopic in $M$, then $fr(C) = fr(C')$. A
free homotopy is realised by a map $f: A \rightarrow M$,
where $A$ is an annulus and where $f(\partial A) = C \cup C'$.
We may ensure that $f^{-1} (S)$ is $\partial A$, together
with some properly embedded arcs and circles in $A$.
We may also ensure that no region of $A - {\rm int}({\cal N}(f^{-1}(S)))$
is a disc. Therefore (using the fact that $A$ is an
annulus) we can guarantee that $f^{-1}(S)$ is a
collection $\alpha_0, \alpha_1, \dots, \alpha_n$ of disjoint essential simple 
closed curves in $A$, where $\partial A = \alpha_0 \cup \alpha_n$. 
Since the image of the annulus lying between $\alpha_i$ and $\alpha_{i+1}$
is disjoint from $S$, then $fr(\alpha_i) = fr(\alpha_{i+1})$.
Therefore, $fr(C) = fr(C')$.

Hence, we have constructed the required infinite
collection of knots. $\square$

We now show that we may ensure that each knot $K$ in
this infinite collection has $M - {\rm int}({\cal N}(K))$
irreducible and atoroidal.

\noindent {\bf Proposition 12.3.} {\sl Let $M$ and $S$ be
as above. Then each essential simple closed curve $K$ on $S$ has
$M - {\rm int}({\cal N}(K))$ irreducible. Also, there are (up to ambient
isotopy in $M$) at most finitely many knots $K$ on $S$ for which
$M - {\rm int}({\cal N}(K))$ is toroidal.}

\noindent {\sl Proof.} Let $K$ be an essential simple closed curve
on $S$. If $M - {\rm int}({\cal N}(K))$ contains
a reducing sphere, then this bounds a ball in $M$.
By assumption, $M$ is irreducible, and so the knot $K$
must lie in this 3-ball, and is therefore homotopically
trivial in $M$. However, $K$ is essential on $S$
and $S$ is $\pi_1$-injective, since it is incompressible.
This is a contradiction.

Suppose now that $M - {\rm int}({\cal N}(K))$ is toroidal, and let $T$
be an essential torus in $M - {\rm int}({\cal N}(K))$. 
Since $M$ is atoroidal, $T$ either is parallel in $M$ to a component of 
$\partial M$ or is compressible in $M$. Consider the former case, and
let $T^2 \times I$ be the parallelity region between $T$ and
a component $T'$ of $\partial M$. The intersection 
$S \cap (T^2 \times I)$ is a collection of discs and annuli,
with $K$ being a core of one of these annuli. Hence, $K$ is parallel to
a curve on $T'$. It is not hard to show that if $K_1$ and $K_2$ are two curves on
$S$ both parallel to curves in $T'$, then either $K_1$ and $K_2$ are ambient isotopic
in $M$ or $S$ contains a component parallel to $T'$.
However, $S$ is connected and non-trivial in $H_2(M, \partial M)$,
which gives a contradiction.

Hence, we may restrict attention to the case where $T$ is compressible
in $M$. Then $T$ bounds a solid torus $V$ in $M$, since $K$
does not lie in a 3-ball. We may assume that
the surface $V \cap S$ is incompressible in $V$
and so is a collection of discs and annuli.
The knot $K$ lies on one such annulus $A$. If $A$ has
winding number one in $V$, then $K$ is a core of $V$
and so $T$ is parallel to $\partial {\cal N}(K)$,
contradicting the assumption that $T$ is essential in
$M - {\rm int}({\cal N}(K))$.
If $A$ has winding number greater than one in $V$,
then a cabling annulus for $K$ is constructed by gluing
$A - {\rm int}({\cal N}(K))$ to the closure of one
of the components of $T - A$. It is now not hard to
show that $K$ is ambient isotopic to the core $K'$ of
an annular component $A'$ of $V \cap S$, where $K'$ has
a cabling annulus disjoint from $S$.

For the purposes of the proof of Proposition 12.3,
we may consider the knot $K'$ instead of $K$.
Suppose therefore that the cabling annulus for $K$
is disjoint from $S$.

\noindent {\sl Claim.} Let $K_1$ and $K_2$ be two essential simple
closed curves on $S$, each cabled with cabling annulus
disjoint from $S$. Suppose that the cabling annuli both
lie on the same side of $S$. Then there is an isotopy
of $S$ which takes $K_1$ off $K_2$.

We may take the cabling annulus $A_i$ for $K_i$
to be properly embedded in $M - {\rm int}({\cal N}(S))$.
Then a component of $M - {\rm int}({\cal N}(S\ \cup A_i))$ is
a solid torus $V_i$. A simple examination
of the intersection between the annulus $A_1$ and
the solid torus $V_2$ establishes that we may isotope
$A_1$ off $V_2$ unless the winding number of $A_2$
along is $V_2$ is two. Thus, the claim is proved unless
the winding number of $A_i$ in $V_i$ is two, for both
$i = 1$ and $2$. In this case $V_i$ is an $I$-bundle
over a M\"obius band, with $V_i \cap \partial {\cal N}(S)$
being precisely the $\partial I$-bundle. Therefore,
if $K_1$ and $K_2$ cannot be homotoped off
each other, $V_1 \cup V_2$ is an $I$-bundle over a connected
non-orientable surface $G$ other than a M\"obius band.
The $I$-bundle over $\partial G$ is a collection of
annuli. If any of these annuli are compressible,
we may extend the $I$-bundle. Thus, we may construct
an $I$-bundle $X$ over a compact connected non-orientable surface $G'$, such that
the $I$-bundle over $\partial G'$ is incompressible in $M$,
and so that $X \cap \partial {\cal N}(S)$ is the
$\partial I$-bundle over $G'$. If $G'$ is a M\"obius band,
then there is an isotopy of $S$ taking $K_1$ off $K_2$.
Suppose therefore that $G'$ has negative Euler characteristic.
Expand the $I$-bundle a little, so that the $\partial I$-bundle lies in $S$.
If we remove the $\partial I$-bundle from $S$, and
attach the $I$-bundle over $\partial G'$, we
create a surface $S'$ with $[S', \partial S'] =
[S, \partial S] \in H_2(M, \partial M)$, and with
$\chi(S') > \chi(S)$. This contradicts the
assumption that $S$ is norm-minimising and
incompressible. 

There are at most finitely many disjoint essential non-parallel
simple closed curves on $S$. This proves the proposition. $\square$

Propositions 12.1, 12.2 and 12.3 give the following result.

\noindent {\bf Theorem 12.4.} {\sl Let $M$ be
a compact orientable 3-manifold with $\partial M$
a (possibly empty) union of tori. Suppose
that $M$ is irreducible and atoroidal, and 
has incompressible boundary. Suppose also that
$H_1(M)$ is torsion free and that some
non-trivial element of $H_2(M, \partial M)$
is represented by a norm-minimising incompressible
surface with positive genus.
Then (up to ambient isotopy) there is an infinite number of surgery curves $K$
in $M$, with exceptional or norm-exceptional surgery slopes $\sigma$ satisfying
$\Delta(\sigma, \mu) = 1$, where $\mu$ is the
meridian slope on $\partial {\cal N}(K)$. We may
ensure that each knot $K$ in this collection
has $M - {\rm int}({\cal N}(K))$ irreducible
and atoroidal, and has $H_2(M - {\rm int}({\cal N}(K)), \partial M) \not= 0$.}

\vskip 18pt
\centerline{\caps 13. References}
\vskip 6pt

\item{1.} D. GABAI, {\sl Foliations and the topology of 3-manifolds II},
J. Differ. Geom. {\bf 26} (1987) 461-478.
\item{2.} W. HAKEN, {\sl Some results on surfaces in 3-manifolds},
Studies in Modern Topology (1968) 39-98.
\item{3.} J. HEMPEL, {\sl 3-Manifolds}, Ann. of Math. Studies, No. 86,
Princeton Univ. Press, Princeton, N. J. (1976)
\item{4.} W. JACO, {\sl Lectures on Three-Manifold Topology},
Regional Conference Series in Mathematics, No. 43, Providence (1980), 
A. M. S.
\item{5.} W. JACO and U. OERTEL, {\sl An algorithm to decide if a 3-manifold
is Haken}, Topology {\bf 23} (1984) 195-209.
\item{6.} M. LACKENBY, {\sl Surfaces, surgery and unknotting operations},
Math. Ann. {\bf 308} (1997) 615-632.
\item{7.} M. LACKENBY, {\sl Dehn surgery on knots in 3-manifolds}, J.
Amer. Math. Soc. {\bf 10} (1997) 835-864.
\item{8.} M. LACKENBY, {\sl Upper bounds in the theory of unknotting
operations}, \hfill\break Topology {\bf 37} (1998) 63-73.
\item{9.} Y. MATHIEU, {\sl Unknotting, knotting by twists and property
(P) for knots in $S^3$}, Knots 90 (ed. A. Kawauchi) Walter de
Gruyter \& Co (1992) 93-102.
\item{10.} M. SCHARLEMANN, {\sl Sutured manifolds and generalized Thurston
norms}, J. Differ. Geom. {\bf 29} (1989) 557-614.
\item{11.} A. THOMPSON, {\sl Thin position and the recognition problem
for $S^3$}, Math. Res. Lett. {\bf 1} (1994) 613-630.
\item{12.} J. TOLLEFSON and N. WANG, {\sl Taut normal surfaces},
Topology {\bf 35} (1996) 55-75.
\item{13.} F. WALDHAUSEN, {\sl The word problem in fundamental groups of
sufficiently large irreducible 3-manifolds}, Ann. Math. {\bf 88} (1968)
272-280.

\vskip 12pt
\+ DPMMS \cr
\+ University of Cambridge \cr
\+ 16 Mill Lane \cr
\+ Cambridge CB2 1SB\cr
\+ England\cr

\end